\documentclass[a4paper]{article}
\usepackage[utf8]{inputenc}
\usepackage[a4paper,left=3cm,right=3cm,top=2cm,bottom=3cm,bindingoffset=5mm]{geometry}
\usepackage[english]{babel}
\usepackage[T1]{fontenc}
\usepackage{graphicx, subfig}
\graphicspath{{img/}}
\usepackage{fancyhdr}
\usepackage{lmodern}
\usepackage{color}

\usepackage{amsfonts}
\usepackage{amsmath}
\usepackage{MnSymbol}
\usepackage{hyperref}

\title{$m$-Axial and $m$-Circular $3m$-Polygons}
\author{Rolf Haag \\
\href{rhaag.98@gmail.com}{rhaag.98@gmail.com}}
\date{}

\begin{document}
\maketitle
\begin{abstract}
The present article includes the enumeration of $n$-polygons with two certain symmetry properties: For a number $3m$ of vertices, we count the $3m$-polygons with $m$ symmetry axes and the $3m$-polygons, that match after three elementary rotations, but have no symmetry axes. For those polygons we give complete lists of representatives of their equivalence-classes and closed formulas for their number.\\

\begin{center}
\textbf{Keywords}
\end{center}
Hamiltonian cycles $\cdot$ Polygons $\cdot$ Symmetry-classes $\cdot$ Euler's $\varphi$-function
\end{abstract}

\tableofcontents

\section{Introduction}
\label{sec:introduction}
\subsection{Definition of a $n$-polygon}
\label{subsec:definition_of_a_n_polygon}
$n$ vertices are regularly distributed in a circle. We consider the Hamiltonian cycles through the $n$ vertices. \cite{Herman2019} In this paper such Hamiltonian cycles are called $n$-polygons. The usual polygons are the special case where all edges have minimal length.\\

Let $ n $ be a natural number $ n\geq3 $ and $ S ^ 1 \subset \mathbb{R} ^ 2 = \mathbb{C} $ the unit circle in the Euclidean plane. The finite subset
\begin{center}
$V_n:=\{v_k:=e^{2 \pi i k/n} \mid k = 0,1,\ldots,n-1 \} \subset S^1$
\end{center}
represents the vertices of an $n$-polygon. \\

To describe the $n$ polygons we use the $n$-cycles $\sigma = (\sigma_1, \sigma_2, \cdots, \sigma_n)$ consisting of the $n$ numbers $\lbrace 0,1, \cdots , n-1 \rbrace$ in any order. The associated $n$-polygon $P(\sigma)$ is given by the path $\overline{v_{\sigma_1}v_{\sigma_2} \cdots v_{\sigma_n}v_{\sigma_1}}$ or more precisely by combining the links $\overline{v_{\sigma_i}v_{\sigma_{i+1}}}$, $i= 1,2, \cdots , n$, where $\sigma_{n+1}= \sigma_1$. \\

Each of the n edges is assigned its ``length'' $e_i$. $e_i = 1$ means, that the $i$-th edge runs counterclockwise from the vertex $V_i$ to the following vertex $V_{i + 1}$. $e_i = 2$ means that the $i$-th edge runs counterclockwise from the vertex $V_i$  to the vertex $V_{i + 2}$ and so on. $e_i = n$ is not possible, since this would mean the connection of the vertex $V_i$ to itself. Therefore, only the numbers between $1$ and $n-1$ are allowed to describe the ``length'' of the edges. The ``length'' of an edge $e_i$ is referred to briefly as a side of the n-polygon. Therefore, an n-polygon can also get described by the $n$-cycle of its sides:
$(e_1, e_2, \ldots, e_i, \ldots, e_n)$.\cite{Brueckner1900}\\

\textbf{About the sums of the sides and the number $u$ of revolutions of any $n$-polygon:}

So that a $n$-polygon does not close prematurely, that is, before all $n-1$ other vertices are passed, no sums of $1, 2, 3, \ldots, n-1$ consecutive sides may be divisible by $n$, the sum $s_n$ of all $n$ sides on the other hand, it must be a multiple of $n$. So if $s_n = u\cdot n$, the integer $u$ is the number of the revolutions in the circle made by the $n$-polygon during its construction.\\

\textbf{Example: Cycle of vertices, cycle of sides, sums of sides and number of revolutions}\\
\begin{figure}[!htp]
\begin{center}
\includegraphics[width=0.4\textwidth]{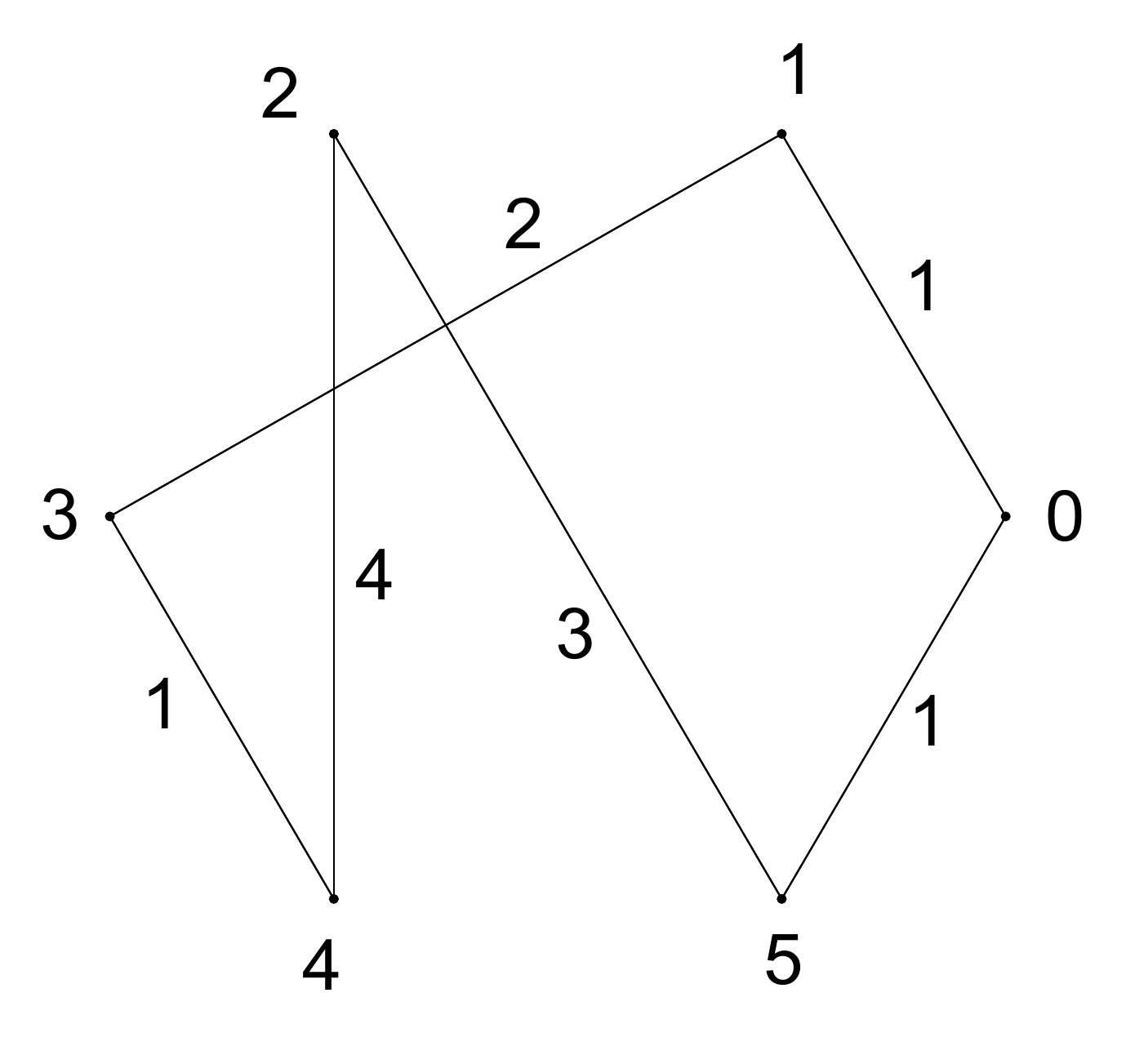}
\caption{n = 6: Cycle of vertices (0 1 3 4 2 5) and cycle of sides (1 2 1 4 3 1)}
\end{center}
\end{figure}

Previous sums $s_1, s_2, \cdots, s_{n-1}$ are not divisible by $n$:\newline
$s_1=1; s_2=1+2=3; s_3=1+2+1=4; s_4=1+2+1+4=8; s_5=1+2+1+4+3=11;$\\

The last sum $s_n$ is divisible by $n$:\newline
$s_6=s_n=1+2+1+4+3+1=12=2 \cdot n$\\

The number of revolutions: $\Rightarrow u=2$

\subsection{Definition of the basic equivalence relation}
\label{subsec:definition_of_the_basic_equivalence_relation}
We denote by $C(n)$ the set of all $n$-polygons and define the following equivalence relation on $C(n)$:\\

Two $n$-polygons $P_1(n)$ and $P_2(n)$ are said to be equivalent, denoted $P_1(n)\stackrel{\equiv}{_E}P_2(n)$, if they are obtainable from one another by a rotation, but not by a reflection.\\

The question of the number $\vert C(n)\stackrel{\equiv}{_E}\vert$ of equivalence classes of the equivalent $n$-polygons has been answered by Golomb/Welch \cite{Golomb1960} and proved by Herman/Poimenidou \cite{Herman2019} in a second way. We now go further into the details:

\subsection{$m$-axial $3m$-polygons - Definition and question}
\label{subsec:m-axial_3m-polygon_-_definition_and_question}
Let $m>2$ be an integer and $n=3m$. A $m$-axial $3m$-polygon is a $3m$-polygon with $m$ symmetry-axes.\\

Those polygons are interesting, because they are ``high'' symmetric. Only the completely regular $n$-polygons with $n$ axes and, provided that $n$ is an even number, the $n$-polygons with $\dfrac{n}{2}$ axes, which were the object of our previous paper \cite{2019arXiv190904124H}, have a ``higher regularity''.\\

In this paper we will show several examples of $m$-axial $3m$-polygons. We ask for the number $\vert P_m(3m)\vert$ of their equivalence classes.

\subsection{$m$-circular $3m$-polygons - Definition and question}
\label{subsec:m-circular_3m-polygon_-_definition_and_question}
Let $m>2$ be an integer and $n=3m$. A $m$-circular $3m$-polygon is a $3m$-polygon without symmetry-axes and the property that it can be made to coincide with itself by $m$ turns through the angles $\dfrac{1 \cdot 3\cdot 2\pi}{n},\dfrac{2\cdot 3 \cdot 2\pi}{n},\cdots,\dfrac{i\cdot 3 \cdot 2\pi}{n},\cdots,\dfrac{\dfrac{n}{3} \cdot 3 \cdot 2\pi}{n}$.\\

We will show several examples of $m$-circular $3m$-polygons. We ask for the number $\vert Q_m(3m)\vert$ of their equivalence classes.

\subsection{Group theoretical aspect of the questions asked}
\label{subsec:group_theoretical_aspect_of_the_questions_asked}

The completely regular $n$-polygons with $n$ axes belong to the complete dihedral group $D_{2n}$, i.e. on completely regular $n$-polygons all $2n$ operations of the dihedral group can be applied. These operations merge each regular $n$-polygon into itself.\\

The $n$-polygons that are the subject of this investigation, on the other hand, belong to certain subgroups of the dihedral group $D_{2n}$: Let $m>2$ be an integer and $n=3m$.
\begin{enumerate}
\item The $m$-axial $3m$-polygons belong to the subgroup with $2m$ elements. These are $m$ rotations and $m$ reflections.
\item The $m$-circular $3m$-polygons belong to the cyclic subgroup with $m$ elements, namely the $m$ rotations by the angles, mentioned in the definition of the $m$-circular $Q_m(3m)$ $3m$-polygons.
\end{enumerate}

\section{Results}
\label{sec:results}
\subsection{$m$-axial $3m$-polygons}
\label{subsec:m-axial_3m-polygons}
\subsubsection{Main-theorem 1}
\label{subsubsec:main_theorem_1}
Let $m>2$ be an integer and $n=3m$.\newline
The different equivalence classes of $n$-polygons with $m$ axes are represented by the $n$-tuples $(a, b, a, a, b, a, \ldots, a, b, a)$ of their sides, if $a$ and $b$ have the following six properties:
\begin{enumerate}
\item $a\in \mathbb{N}$ with $a\equiv $1 mod 3,
\item $b\in \mathbb{N}$ with $b\equiv $1 mod 3,
\item $gcd\left(u,m\right)=1$,
\item $1 \leq a \leq n-2$,
\item $1 \leq b \leq n-2$,
\item $a \neq b$,
\end{enumerate}
\subsubsection{Formula for $\vert P_m(3m)\vert$}
\label{subsubsec:formula}
Let $m>2$ be an integer and $n=3m$ and $\varphi(m)$ denote the Euler $\varphi$-function:\\

The number $\vert P_m(3m)\vert$ of the equivalence-classes of the $m$-axial $3m$-polygons is:
\begin{center}
$\vert P_m(3m)\vert=m\cdot \varphi(m)-\dfrac{\varphi(3m)}{2}$.
\end{center}
\begin{table}[!htp]
\centering
\begin{tabular}{| c | c | c || c | c | c || c | c | c || c | c | c |}
\hline
n & m & $\vert P_m(n) \vert$ & n & m & $\vert P_m(n) \vert$ & n & m & $\vert P_m(n) \vert$ & n & m & $\vert P_m(n) \vert$\\ \hline
9 & 3 & 3 & 12 & 4 & 6 & 15 & 5 & 16 & 18 & 6 & 9 \\
21 & 7 & 36 & 24 & 8 & 28 & 27 & 9 & 45 & 30 & 10 & 36\\
33 & 11 & 100 & 36 & 12 & 42 & 39 & 13 & 144 & 42 & 14 & 78\\
45 & 15 & 108 & 48 & 16 & 120 & 51 & 17 & 256 & 54 & 18 & 99\\
57 & 19 & 324 & 60 & 20 & 152 & 63 & 21 & 234 & 66 & 22 & 210\\
69 & 23 & 484 & 72 & 24 & 180 & 75 & 25 & 480 & 78 & 26 & 300\\
81 & 27 & 459 & 84 & 28 & 324 & 87 & 29 & 784 & 90 & 30 & 228\\
\hline
\end{tabular}
\caption{Number of equivalence classes of $3m$-polygons with $m$ axes}
\label{tab:number_of_equivalence_classes_of_n-polygons_with_m_axes}
\end{table}
\subsubsection{Special case: $\vert P_p(3p) \vert$ with  a prime number $p>3$}
\label{subsubsec:special case}
Let $p$ be a prime number with $p>3$.\\

The number $\vert P_p(3p)\vert$ of the equivalence-classes of the $p$-axial $3p$-polygons is:
\begin{center}
$\vert P_p(3p)\vert=\left(p-1\right)^2$.
\end{center}

\newpage
\subsubsection{Examples of $m$-axial $3m$-polygons}
\label{examples of m-axial 3m-polygons}
\paragraph[{$n=9$; $m=3$; $\vert P_3(9) \vert=3$}]{$n=9$; $m=3$; $\vert P_3(9) \vert=3$}.
\begin{figure}[!htp]
\centering
\begin{tabular}{c | c | c }
\includegraphics[width=0.3\textwidth]{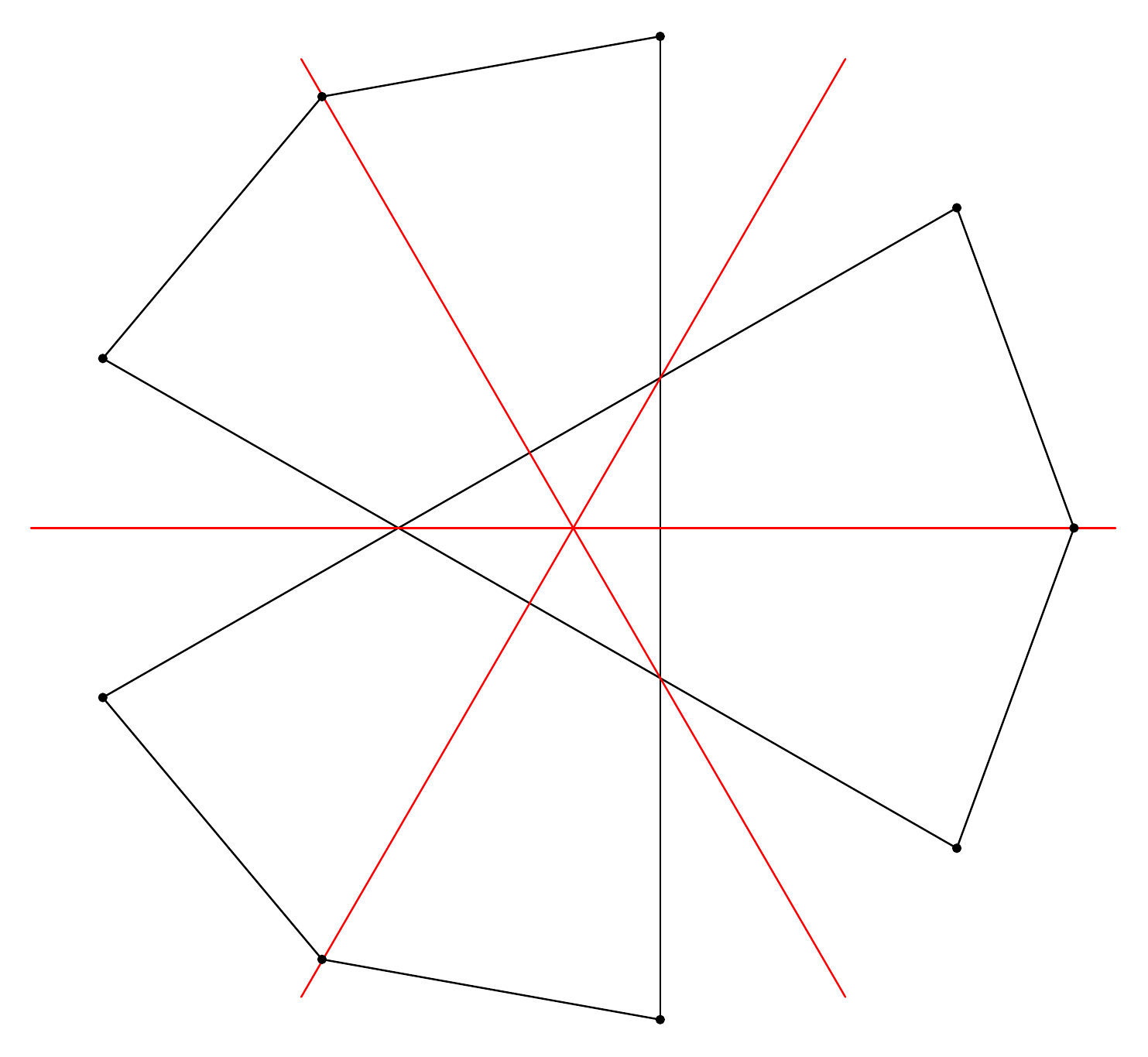} & \includegraphics[width=0.3\textwidth]{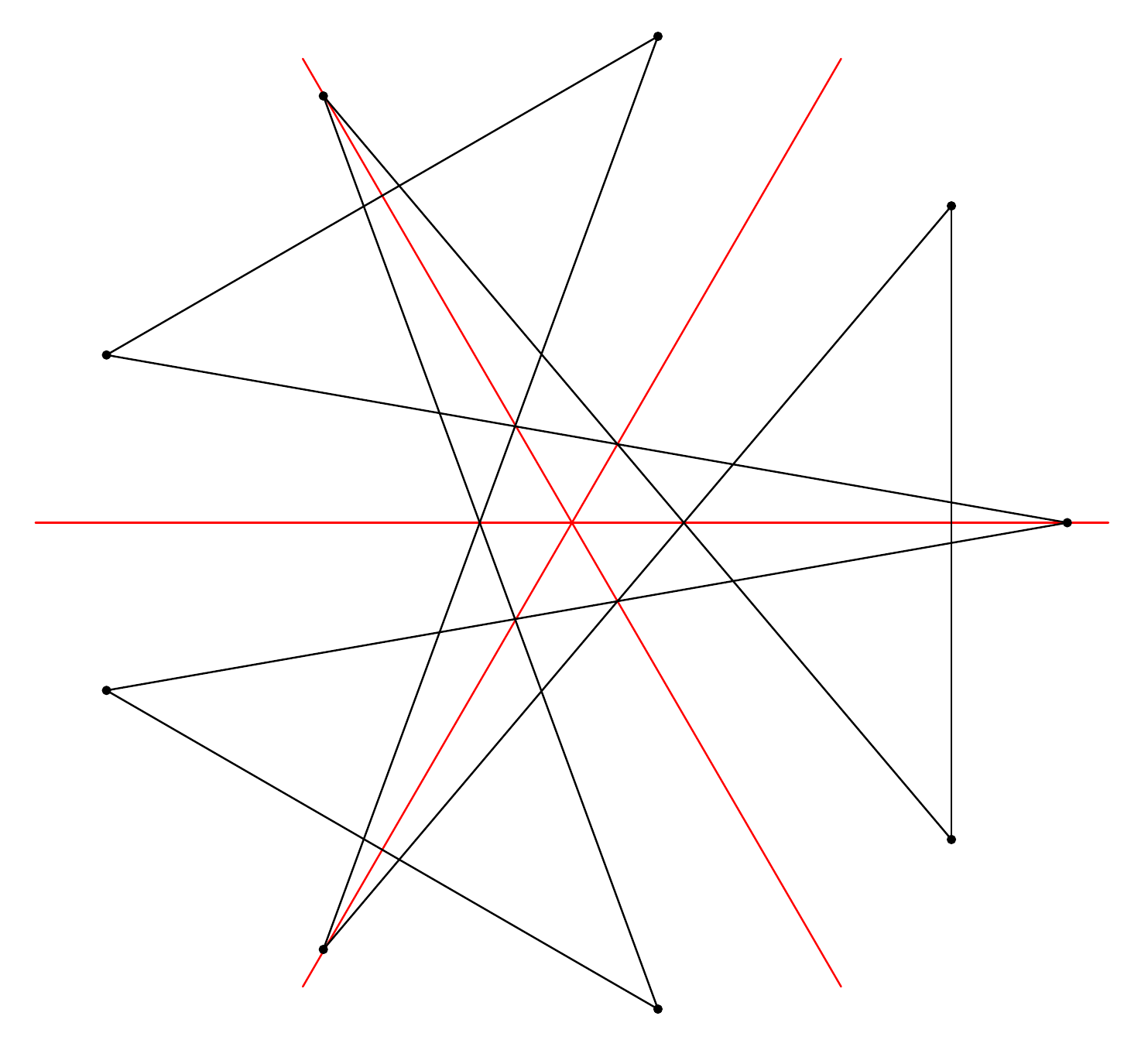} & \includegraphics[width=0.3\textwidth]{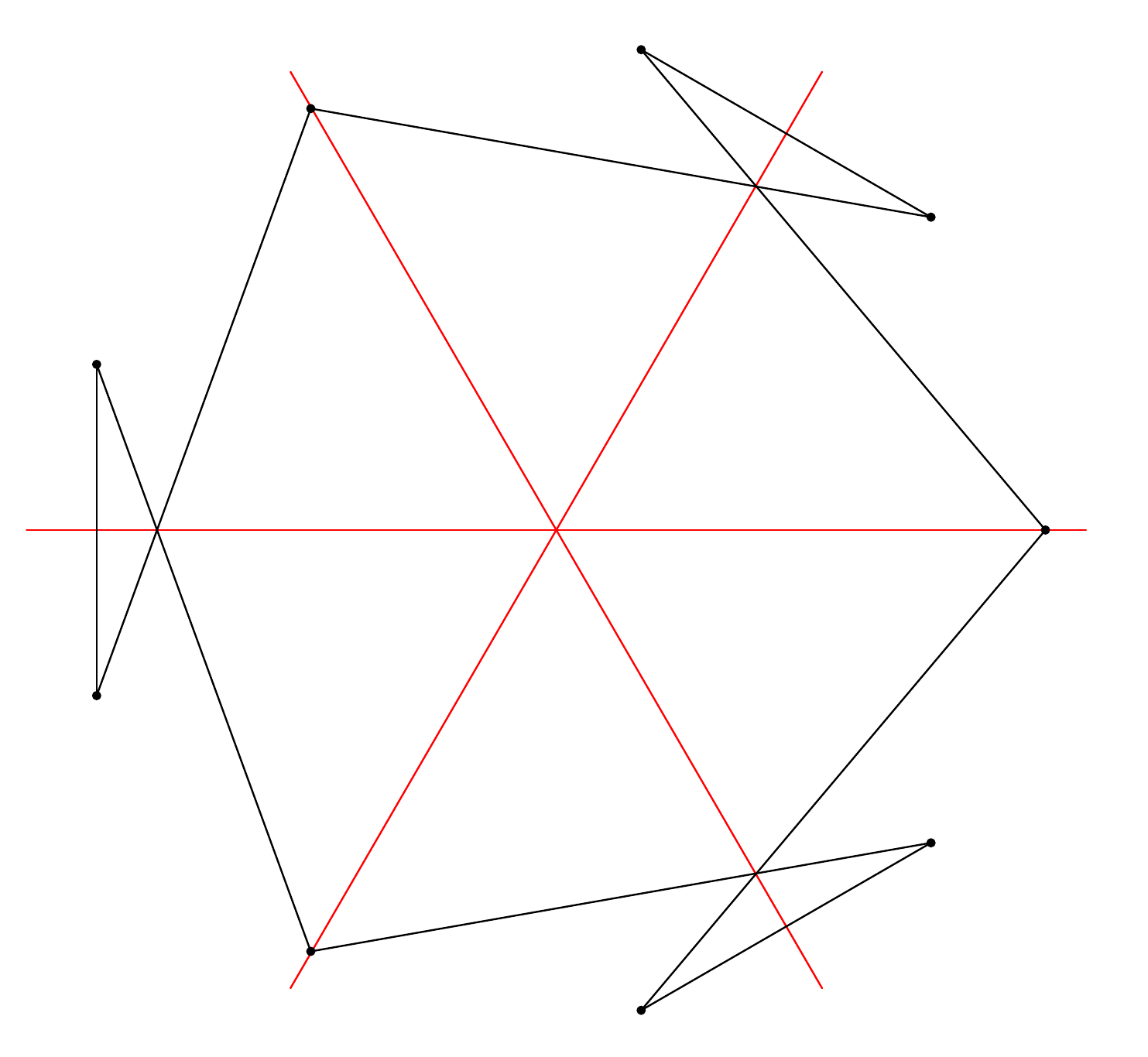}\\
$a=1;b=4$ & $a=4;b=7$ & $a=7;b=1$
\end{tabular}
\caption{$n=9$: A set of representatives of the different equivalence-classes of $P_3(9)$}
\end{figure}
\paragraph[$n=12$; $m=4$; $\vert P_4(12) \vert=6$]{$n=12$; $m=4$; $\vert P_4(12) \vert=6$}.
\begin{figure}[!htp]
\centering
\begin{tabular}{c | c | c }
\includegraphics[width=0.3\textwidth]{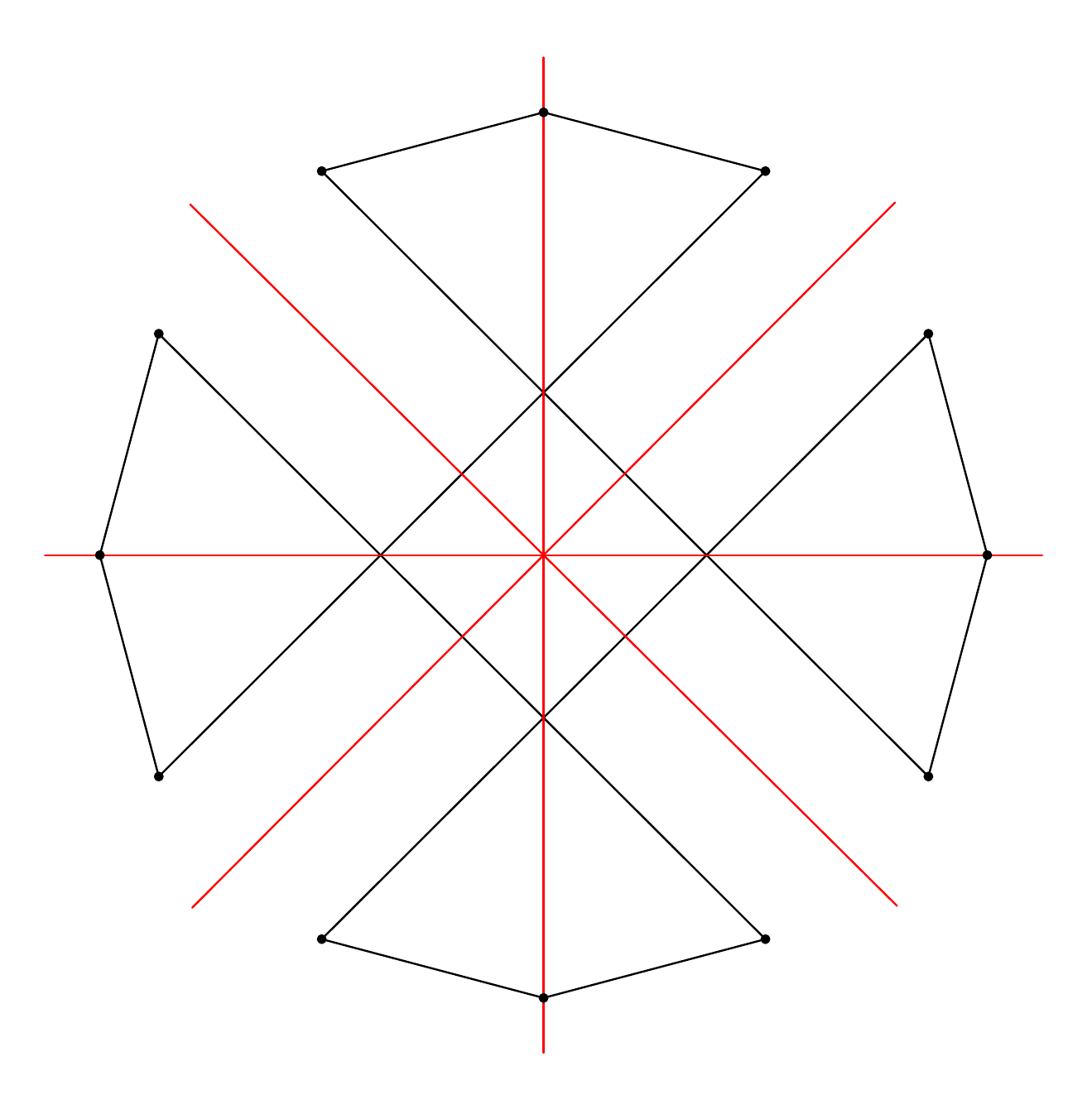} & \includegraphics[width=0.3\textwidth]{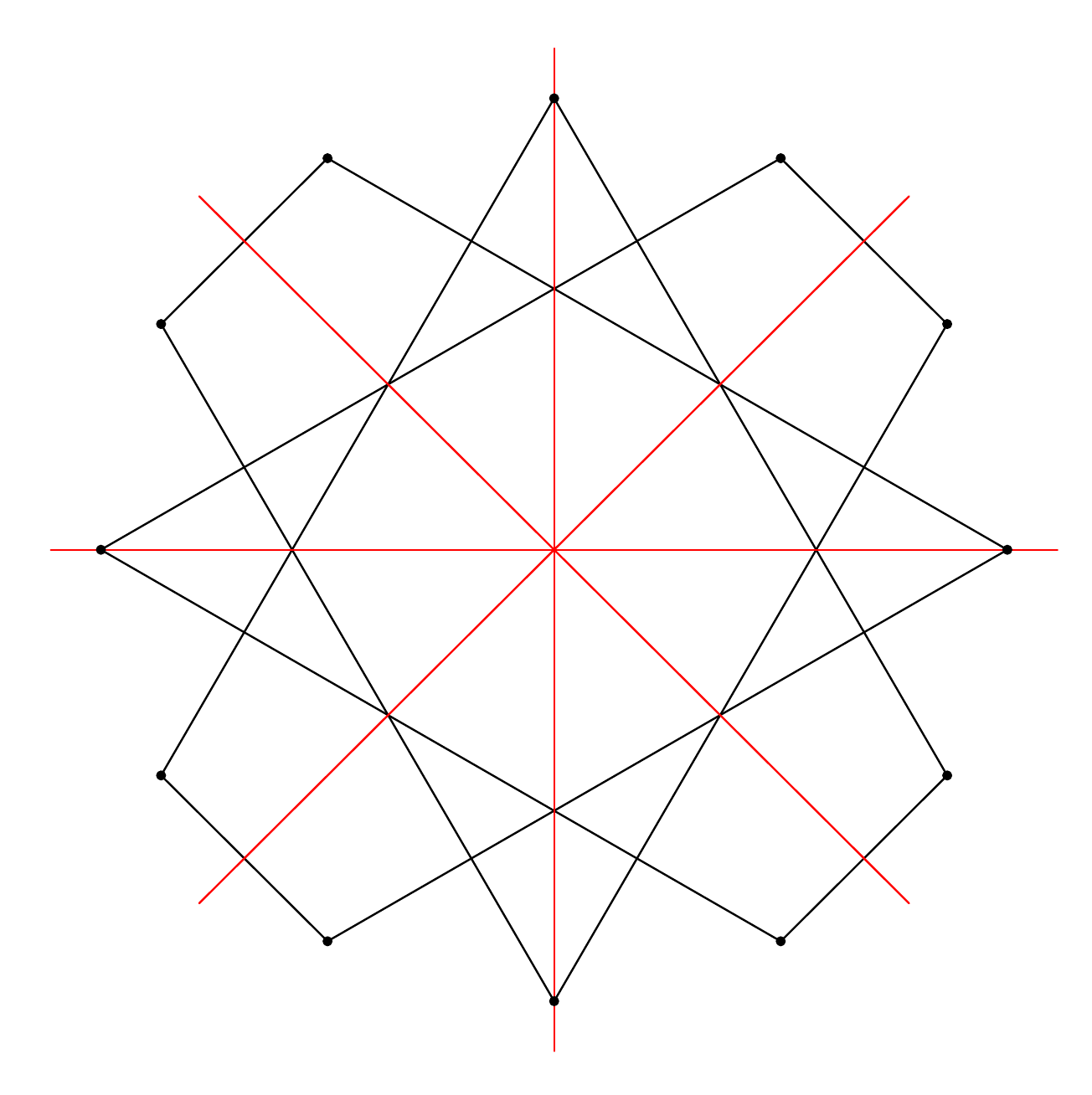} & \includegraphics[width=0.3\textwidth]{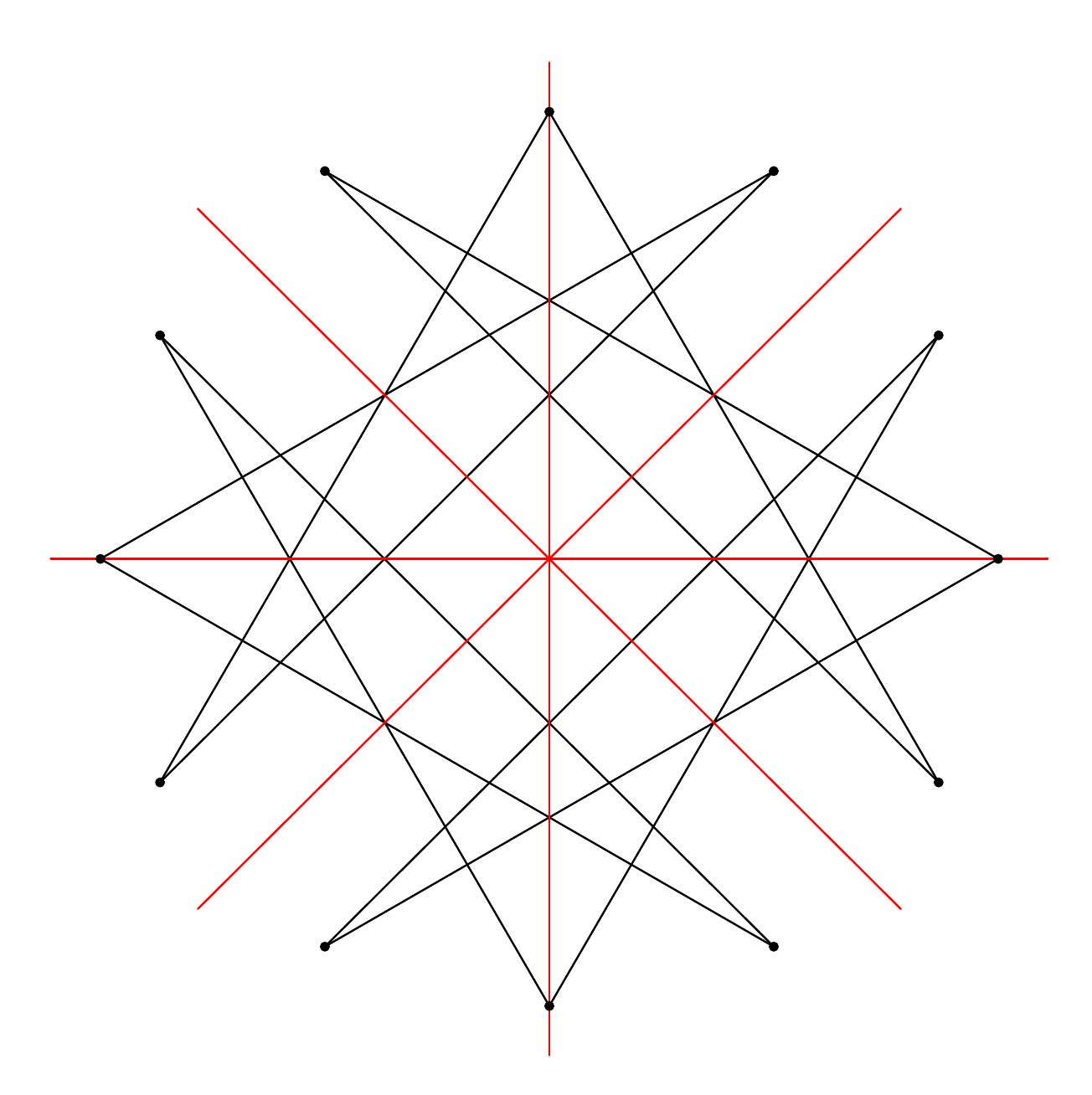}\\
$a=1;b=7$ & $a=4;b=1$ & $a=4;b=7$ \\ \hline
\includegraphics[width=0.3\textwidth]{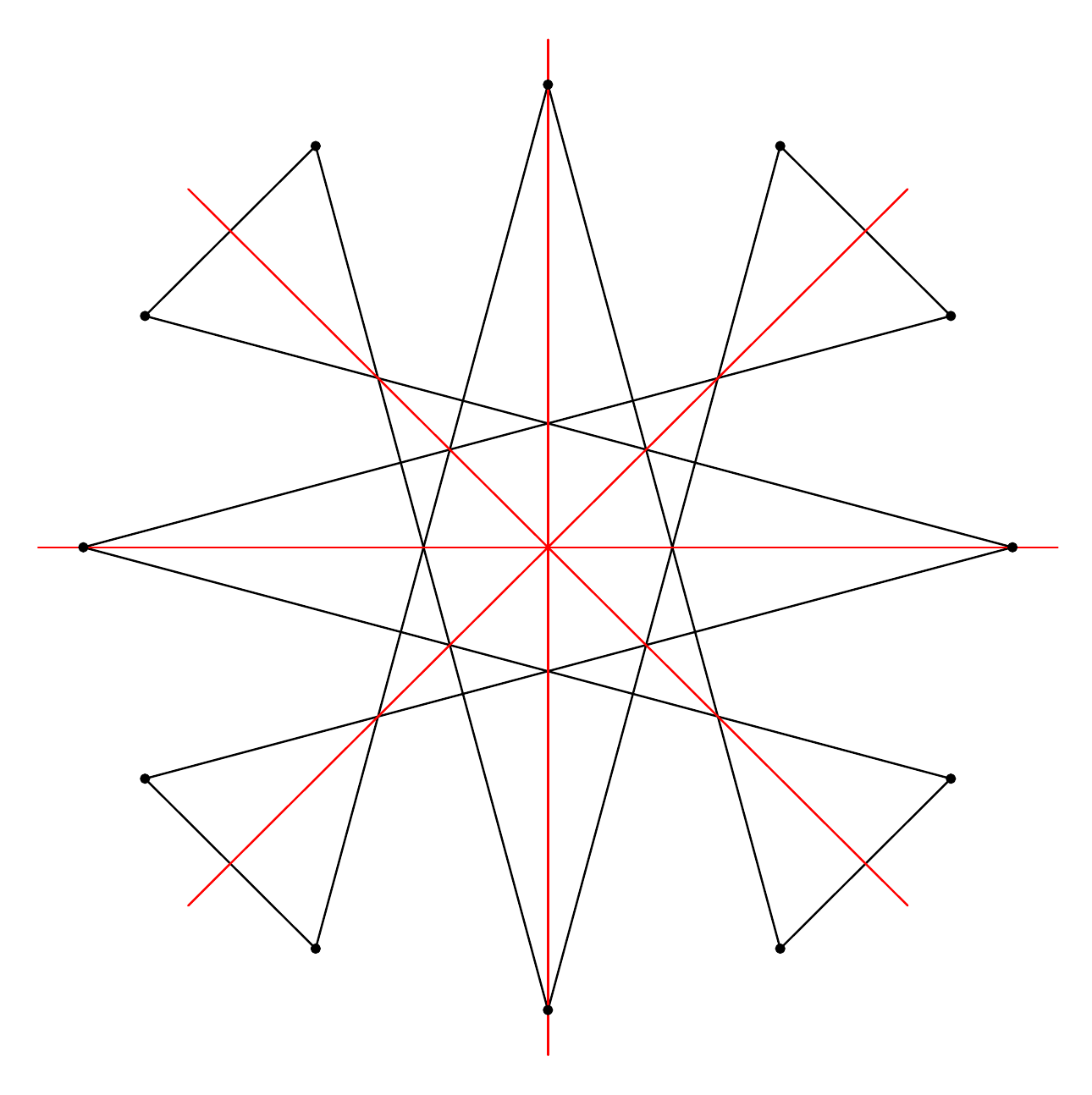} & \includegraphics[width=0.3\textwidth]{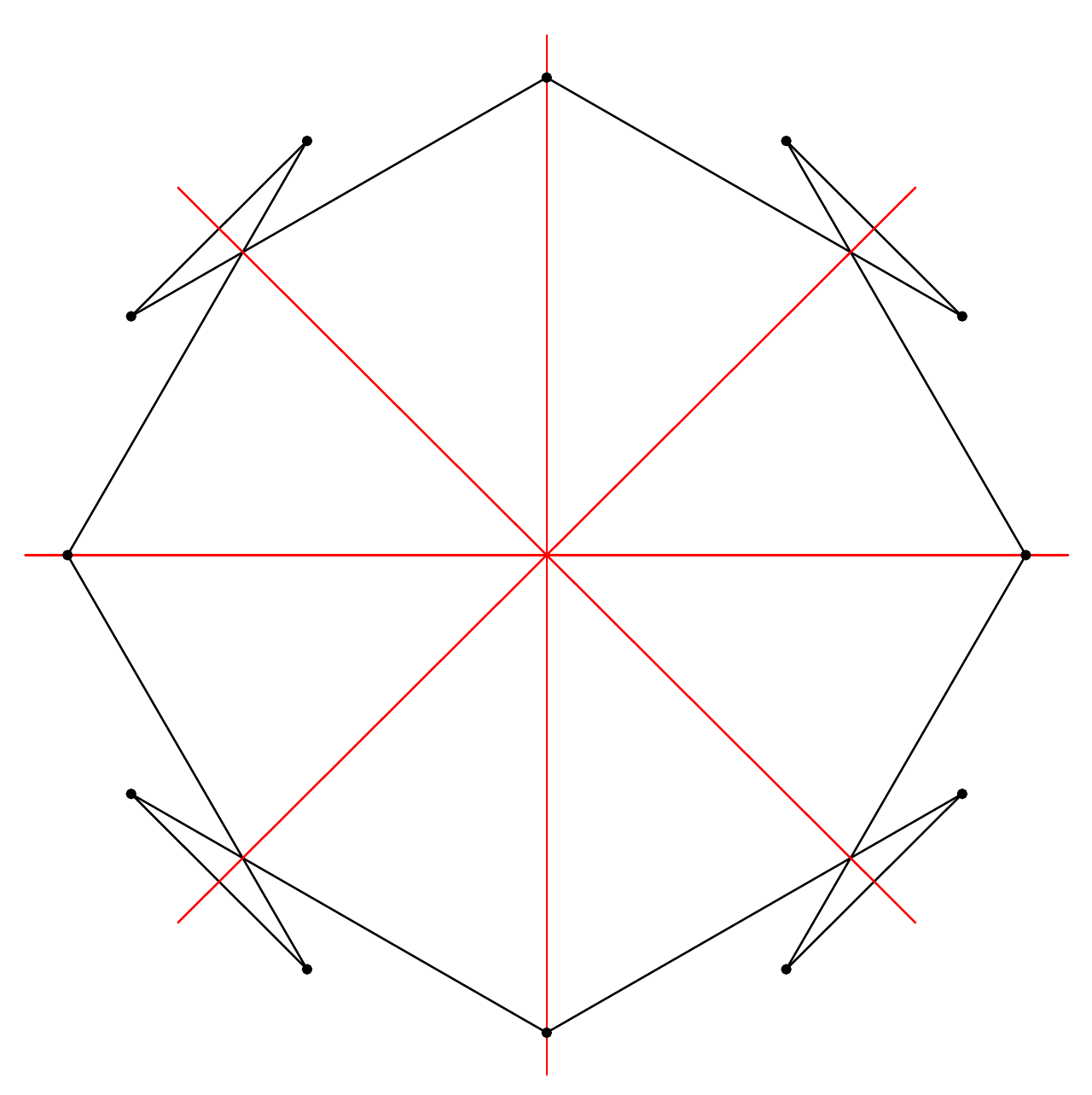} & \includegraphics[width=0.3\textwidth]{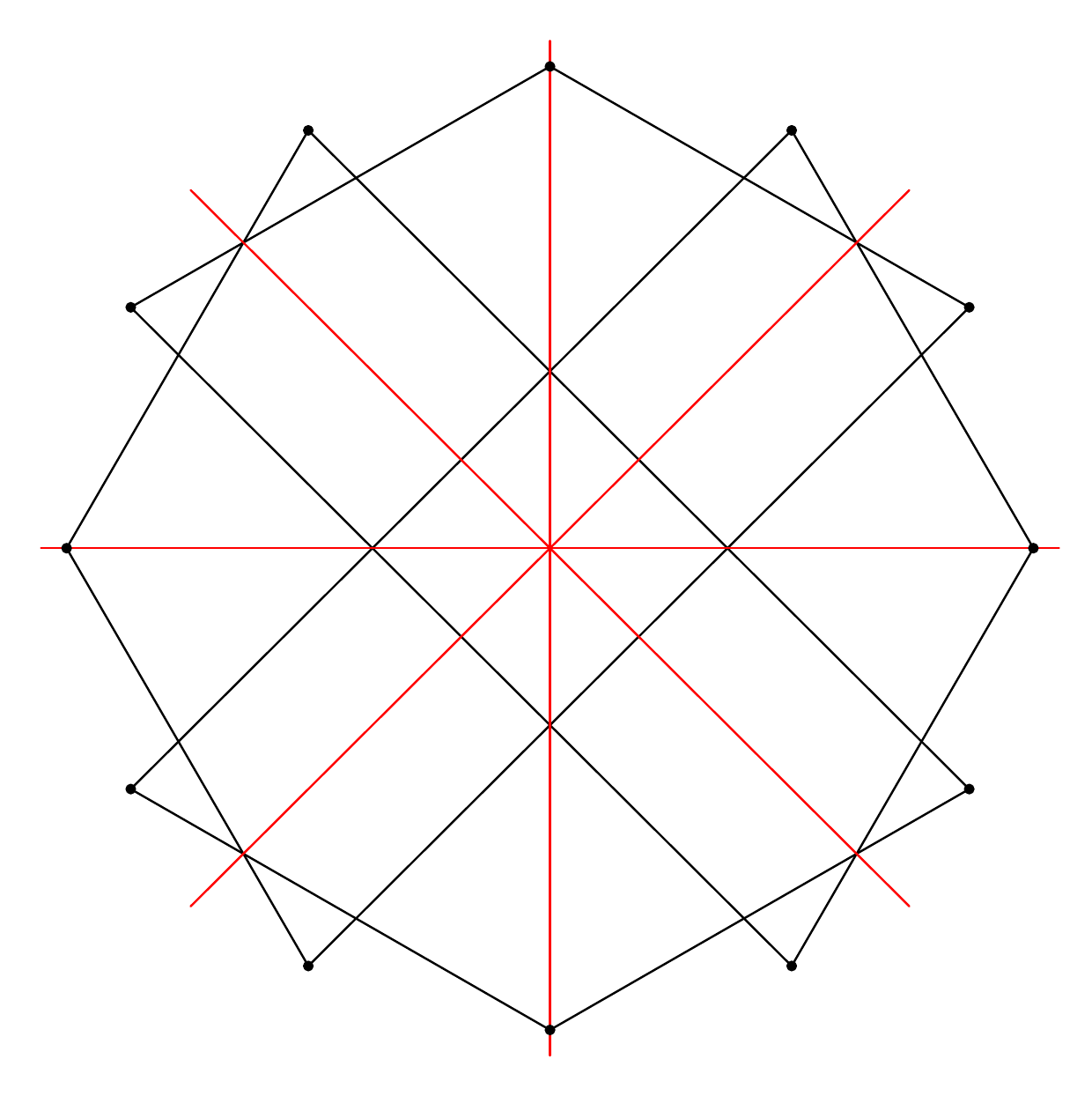}\\
$a=7;b=1$ & $a=10;b=1$ & $a=10;b=7$ \\
\end{tabular}
\caption{$n=12$: A set of representatives of the different equivalence-classes of $P_4(12)$}
\end{figure}
\newpage
\paragraph[$n=15$; $m=5$; $\vert P_5(15) \vert=16$]{$n=15$; $m=5$; $\vert P_5(15) \vert=16$}.
\begin{figure}[!htp]
\centering
\begin{tabular}{c | c | c | c }
\includegraphics[width=0.2\textwidth]{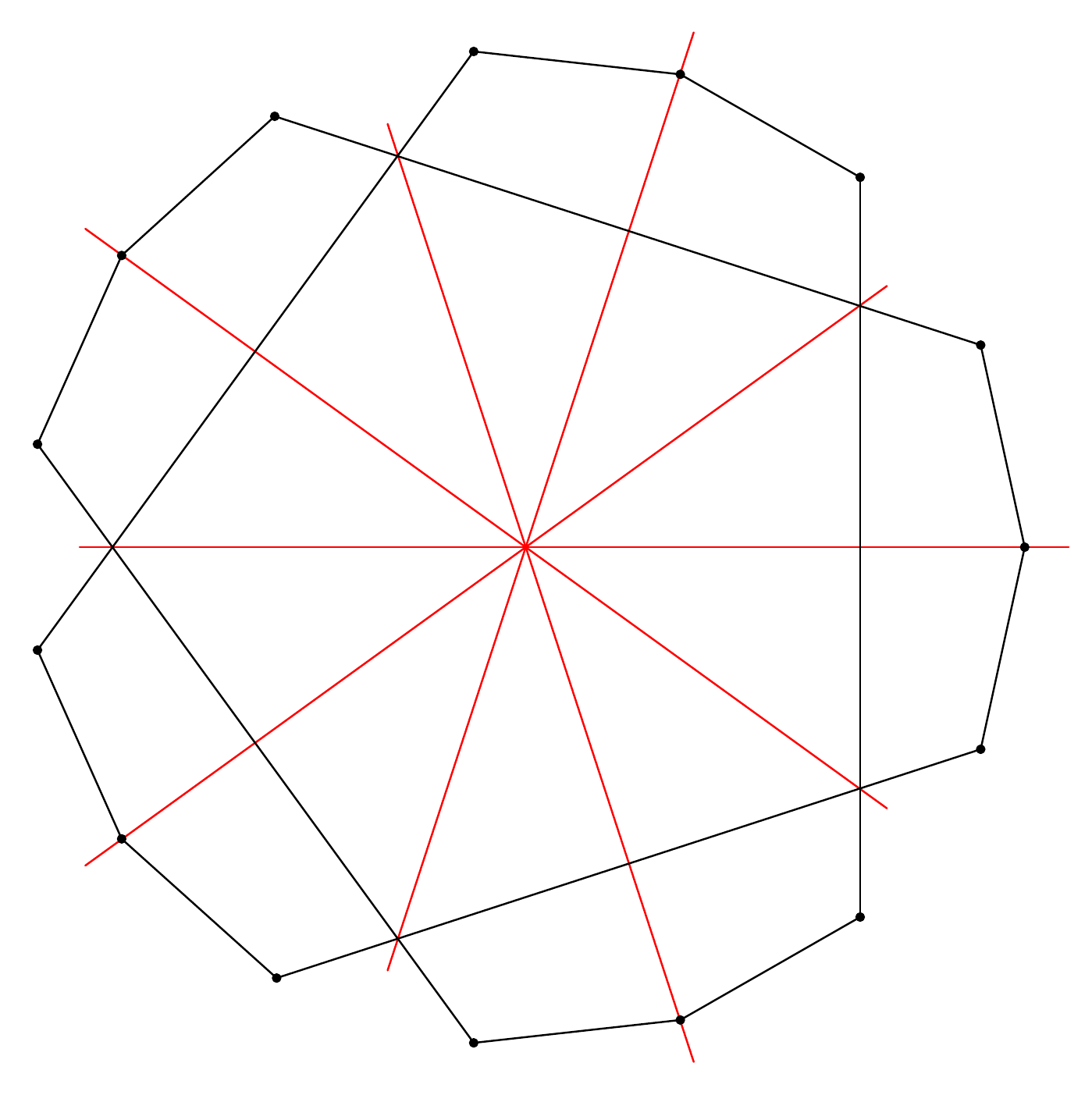} & \includegraphics[width=0.2\textwidth]{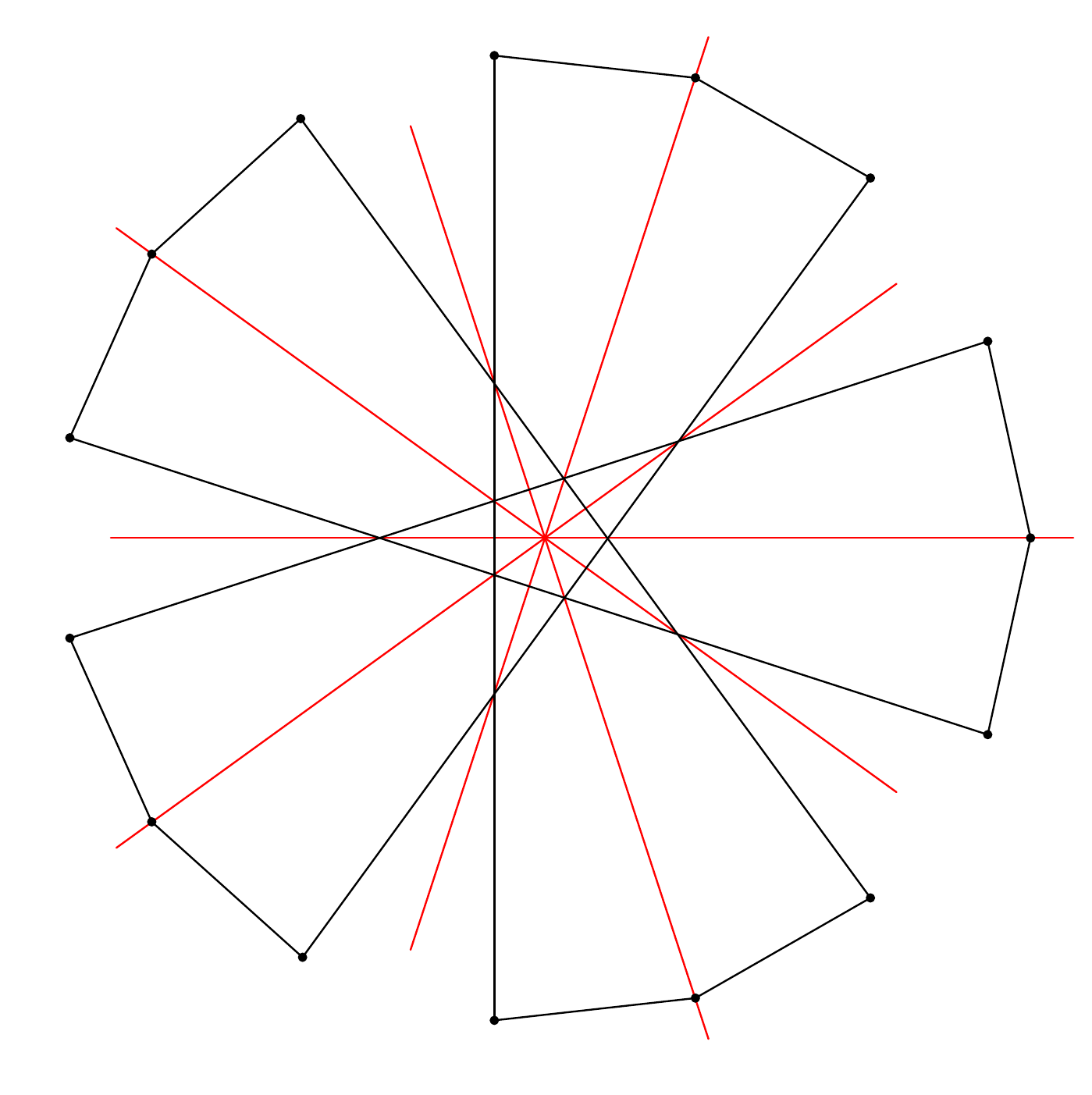} & \includegraphics[width=0.2\textwidth]{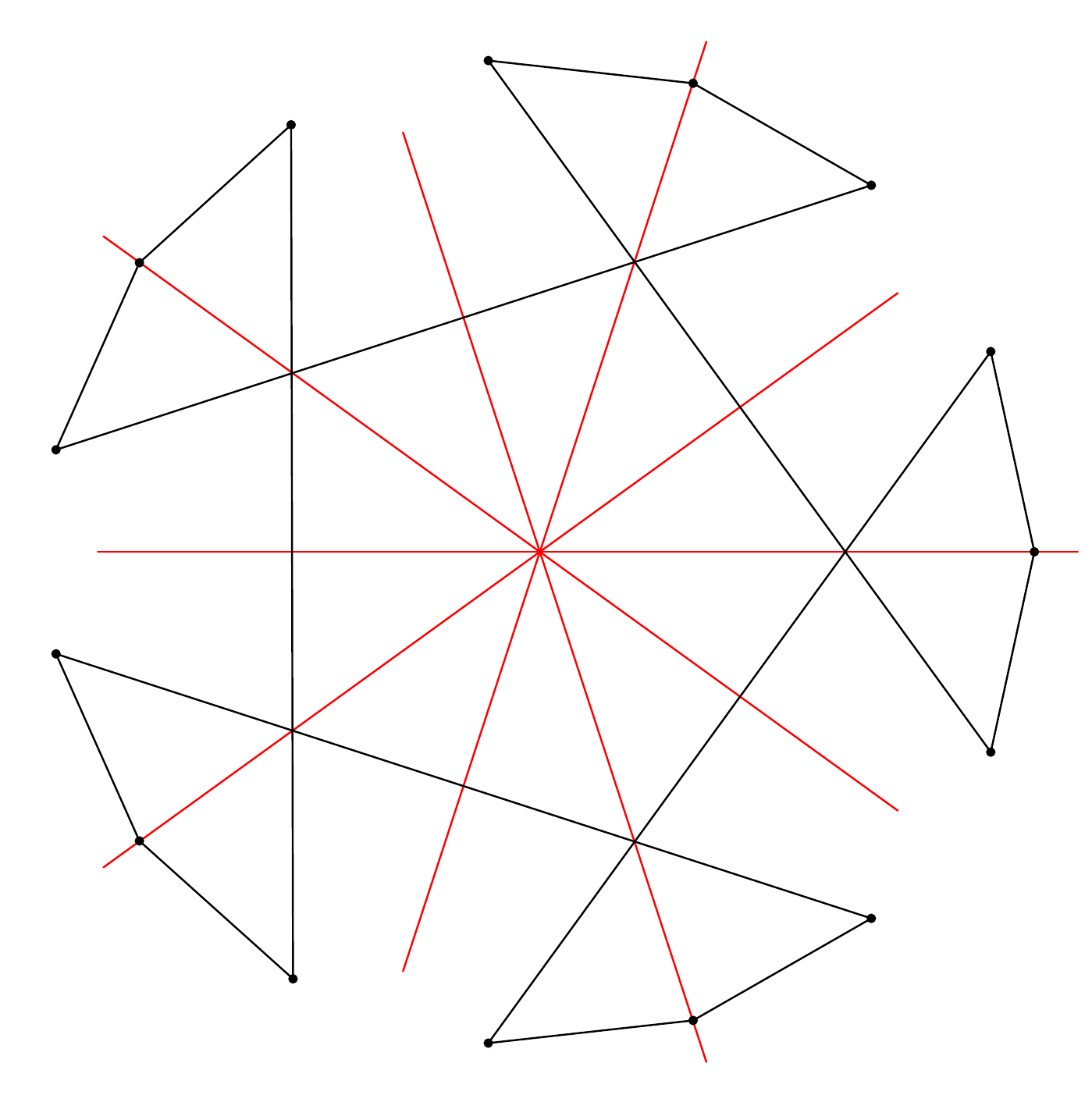} & \includegraphics[width=0.2\textwidth]{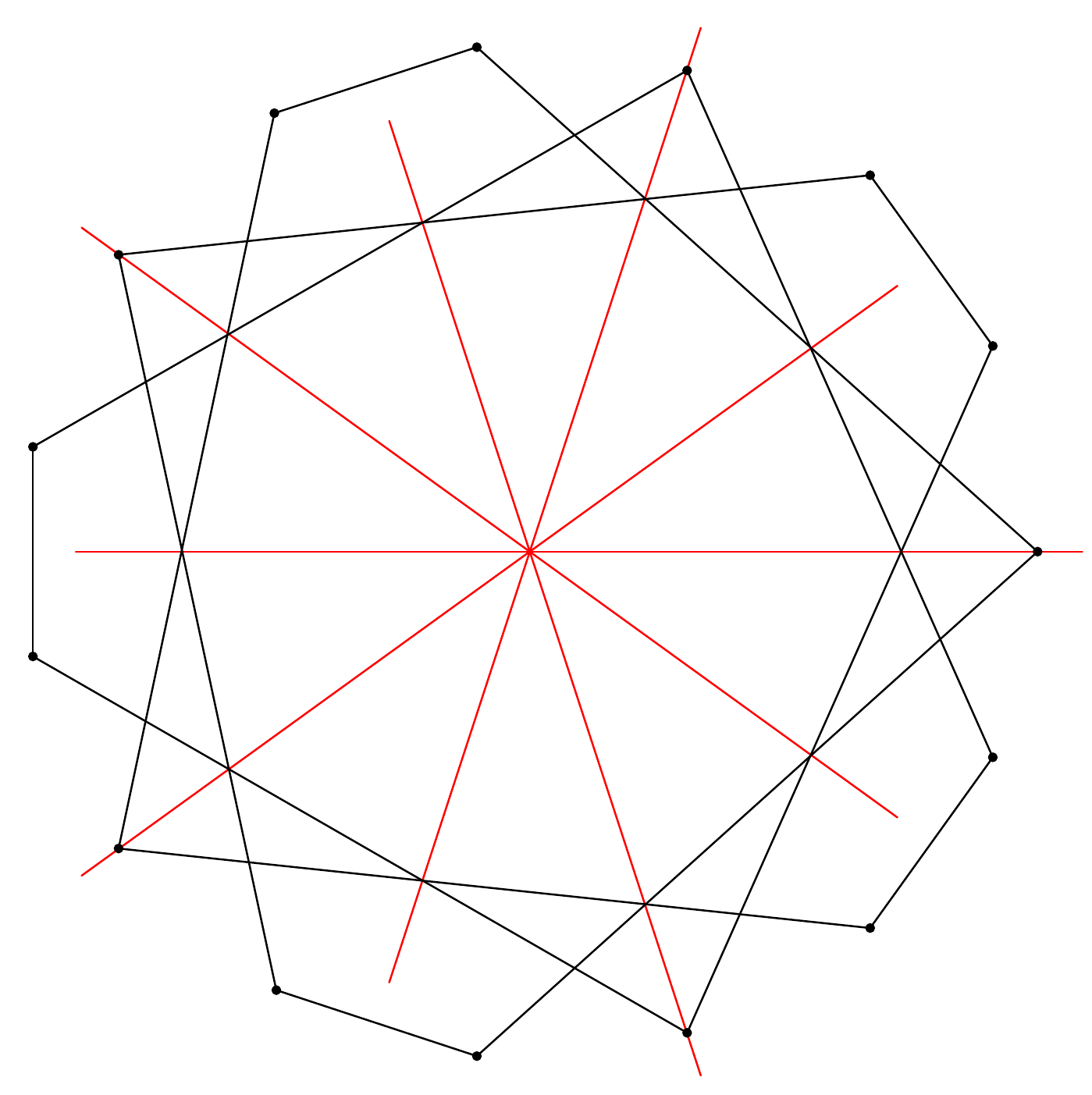}\\
$a=1;b=4$ & $a=1;b=7$ & $a=1;b=10$ & $a=4;b=1$\\ \hline
\includegraphics[width=0.2\textwidth]{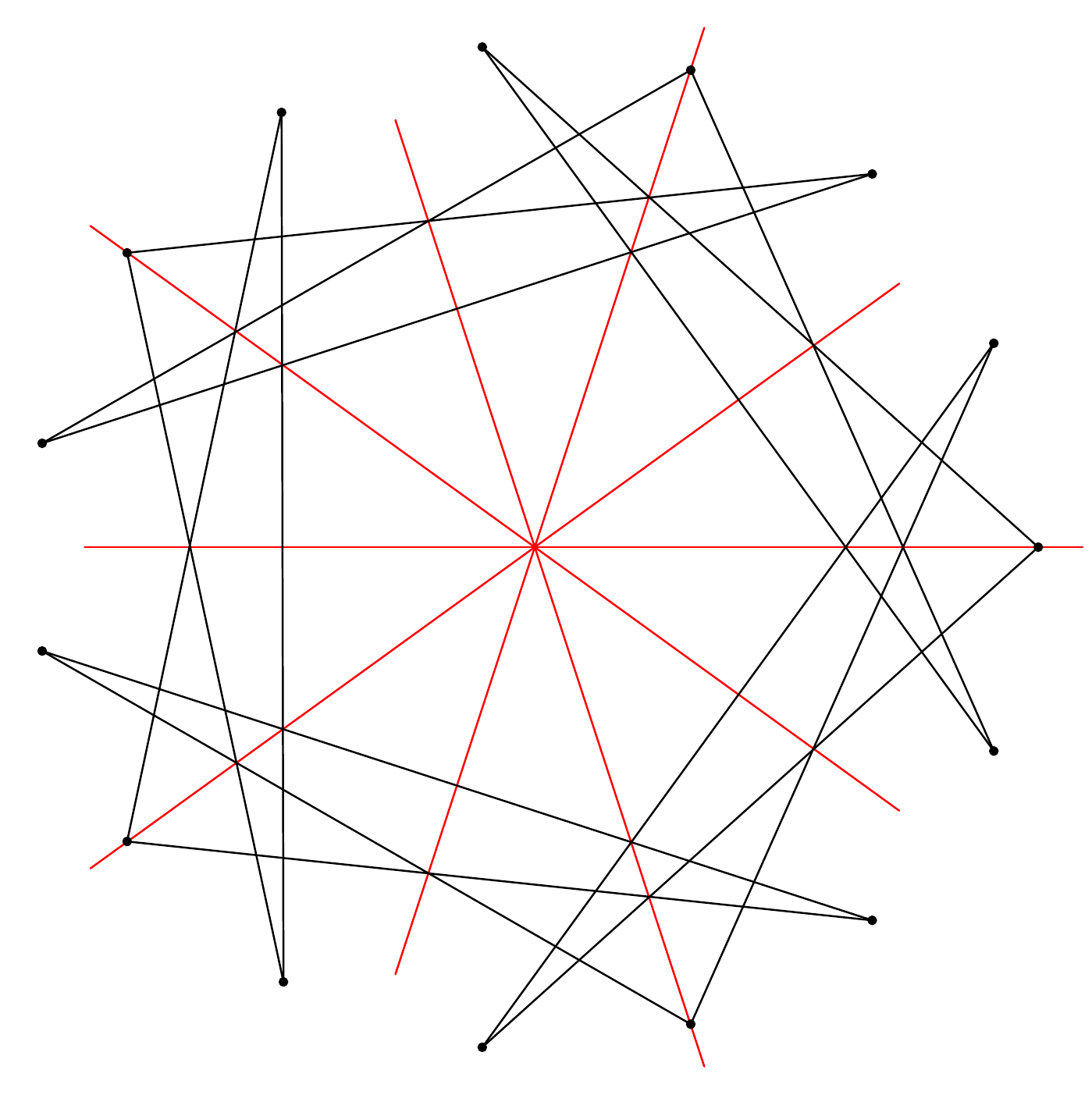} & \includegraphics[width=0.2\textwidth]{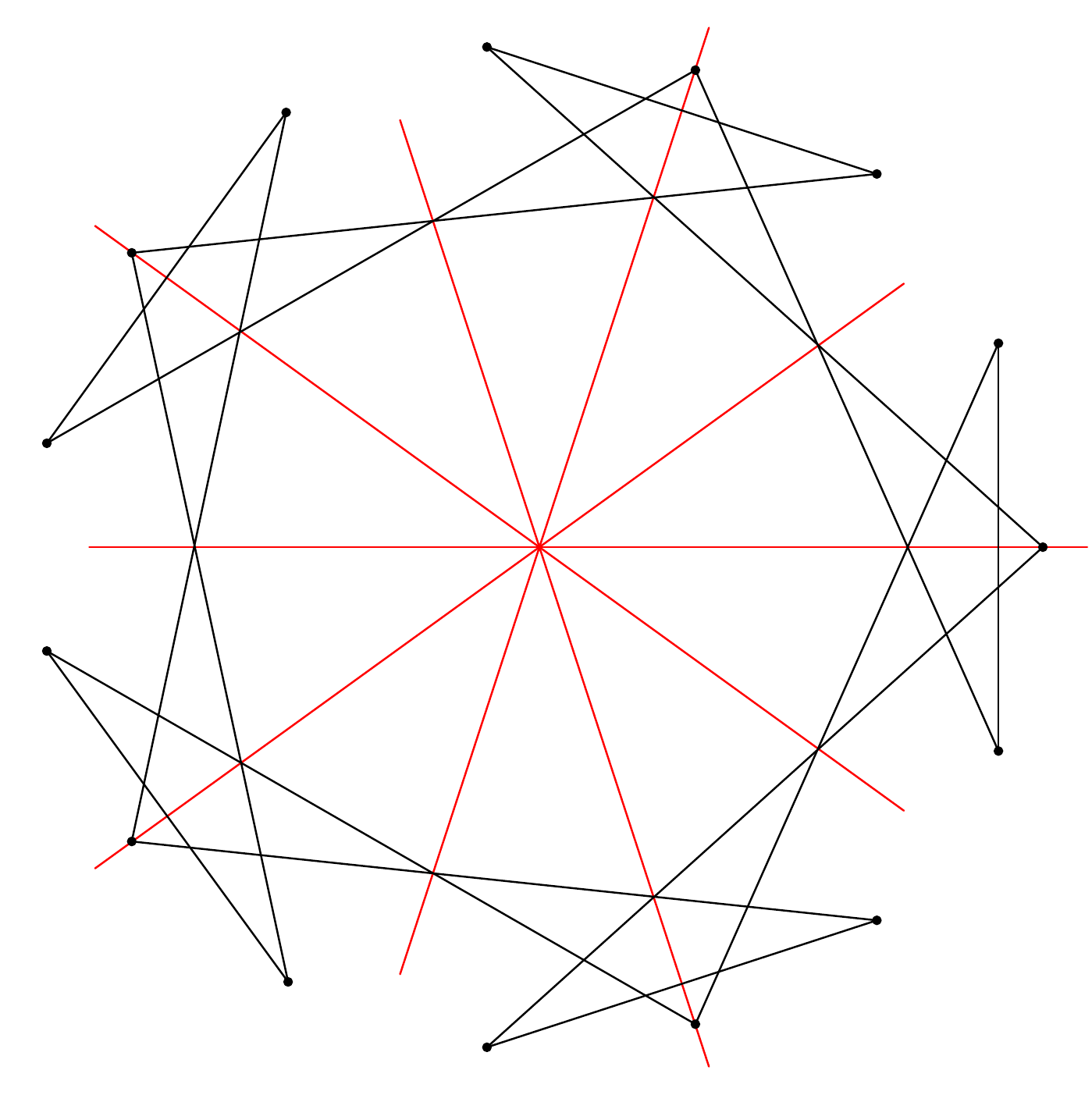} & \includegraphics[width=0.2\textwidth]{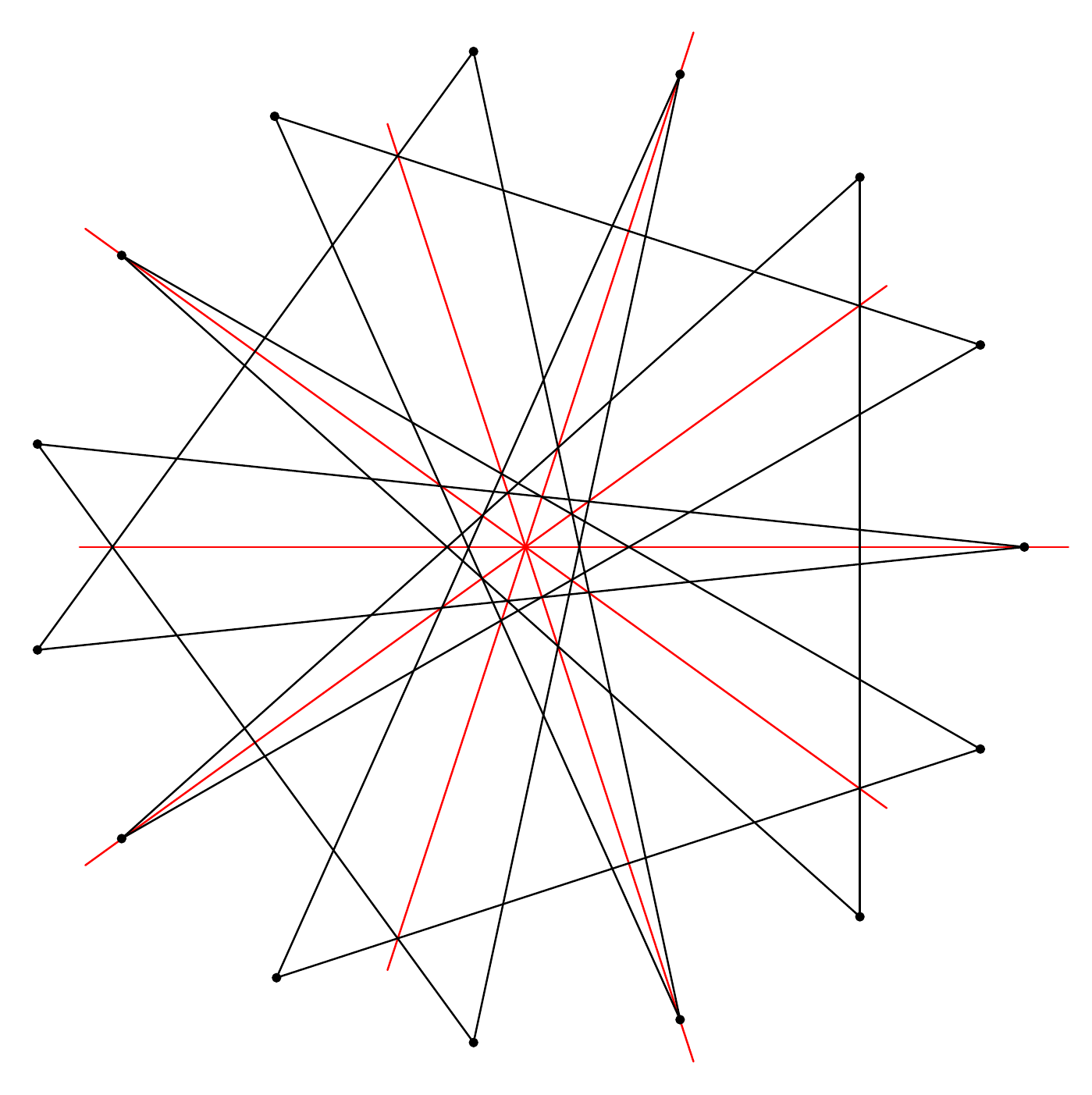} & \includegraphics[width=0.2\textwidth]{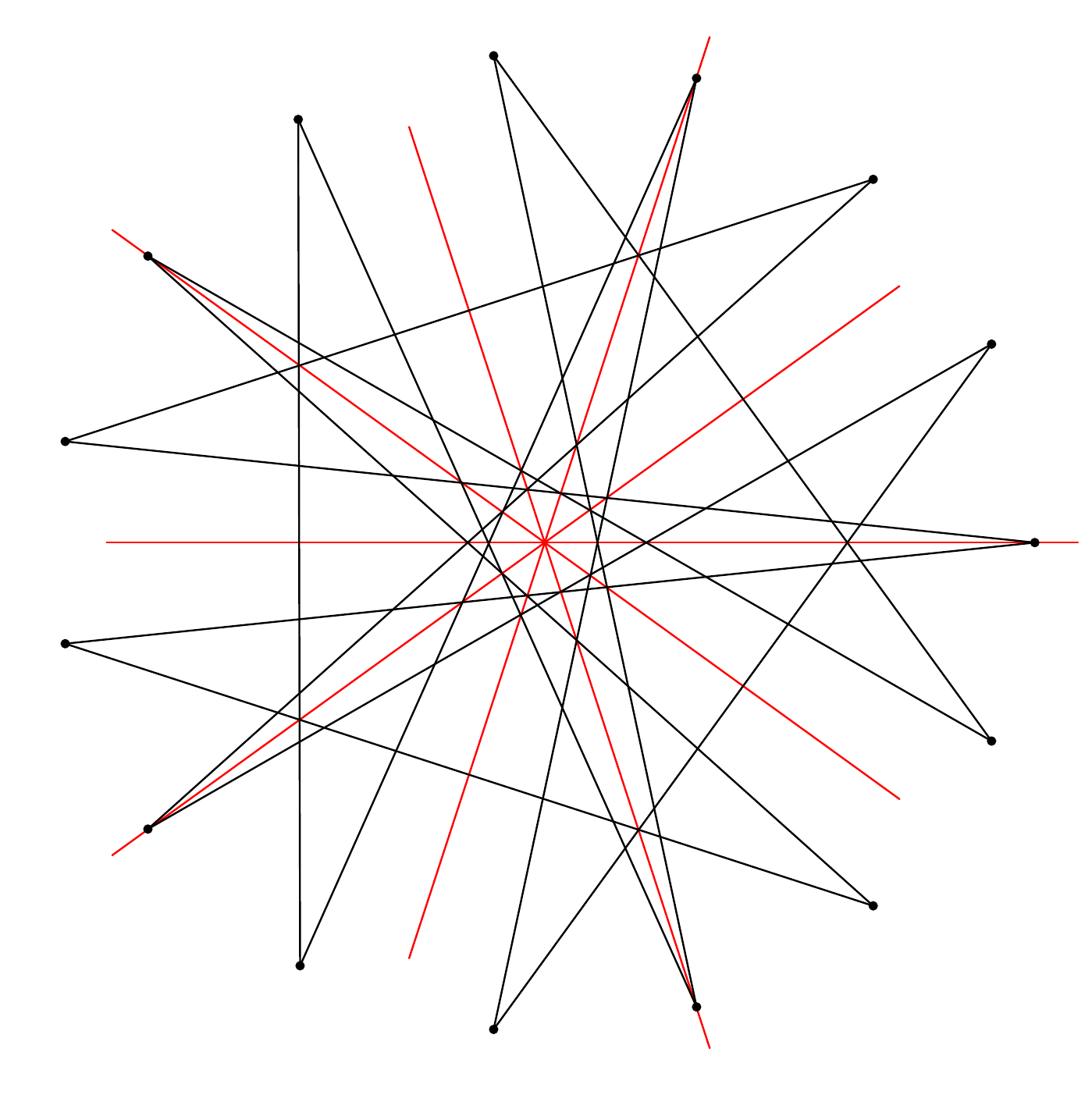}\\
$a=4;b=10$ & $a=4;b=13$ & $a=7;b=4$ & $a=7;b=10$\\ \hline
\includegraphics[width=0.2\textwidth]{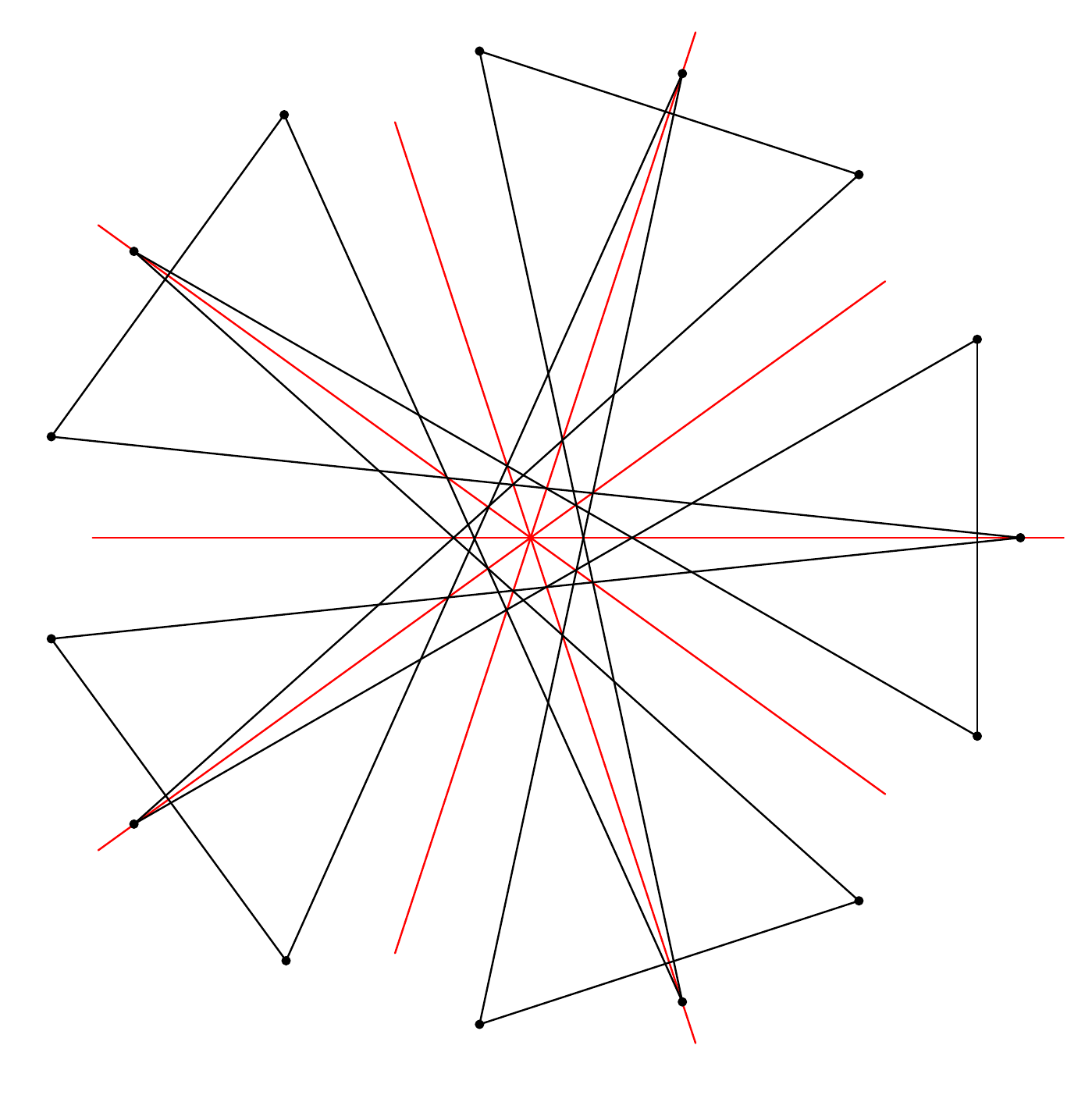} & \includegraphics[width=0.2\textwidth]{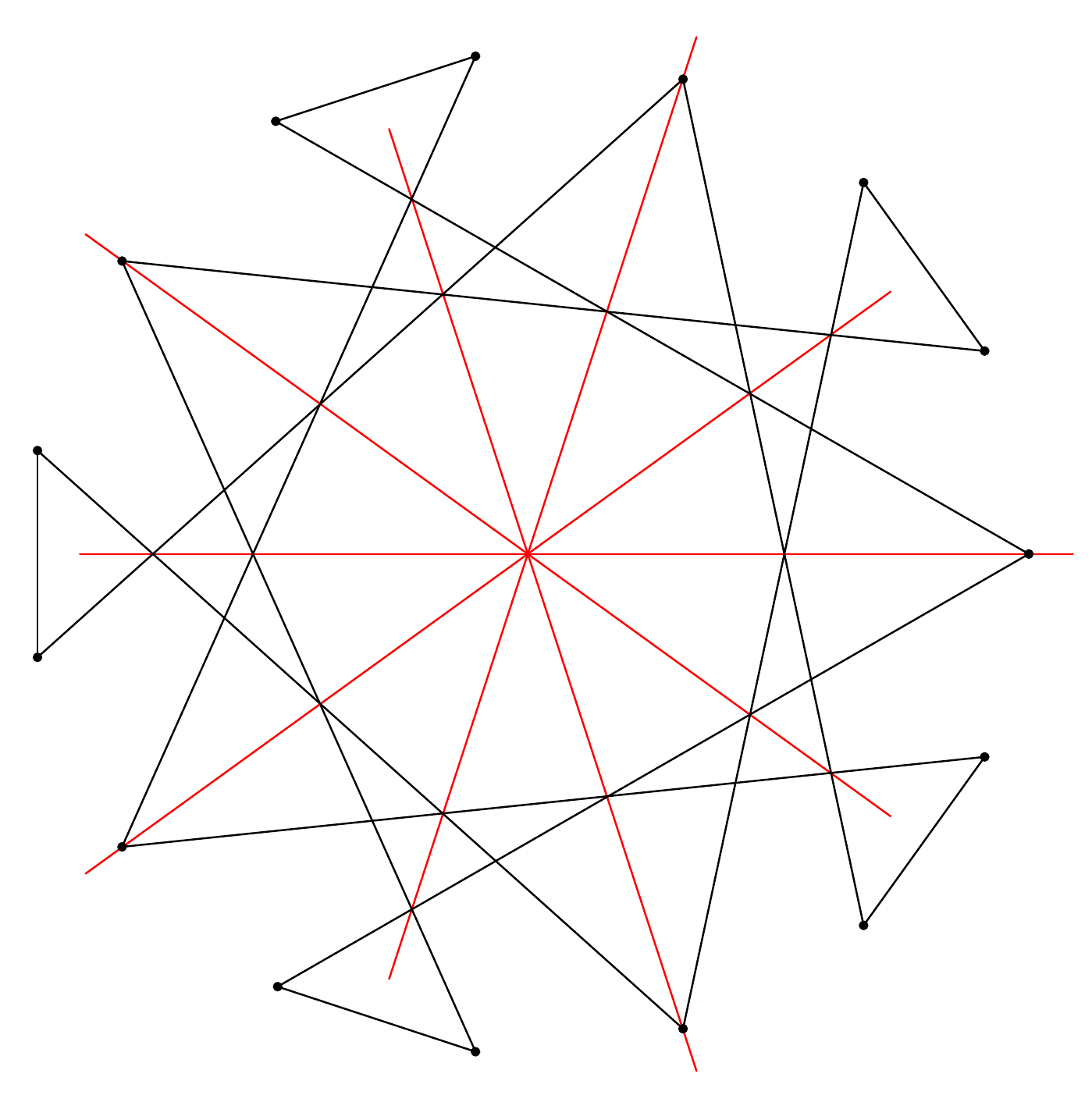} & \includegraphics[width=0.2\textwidth]{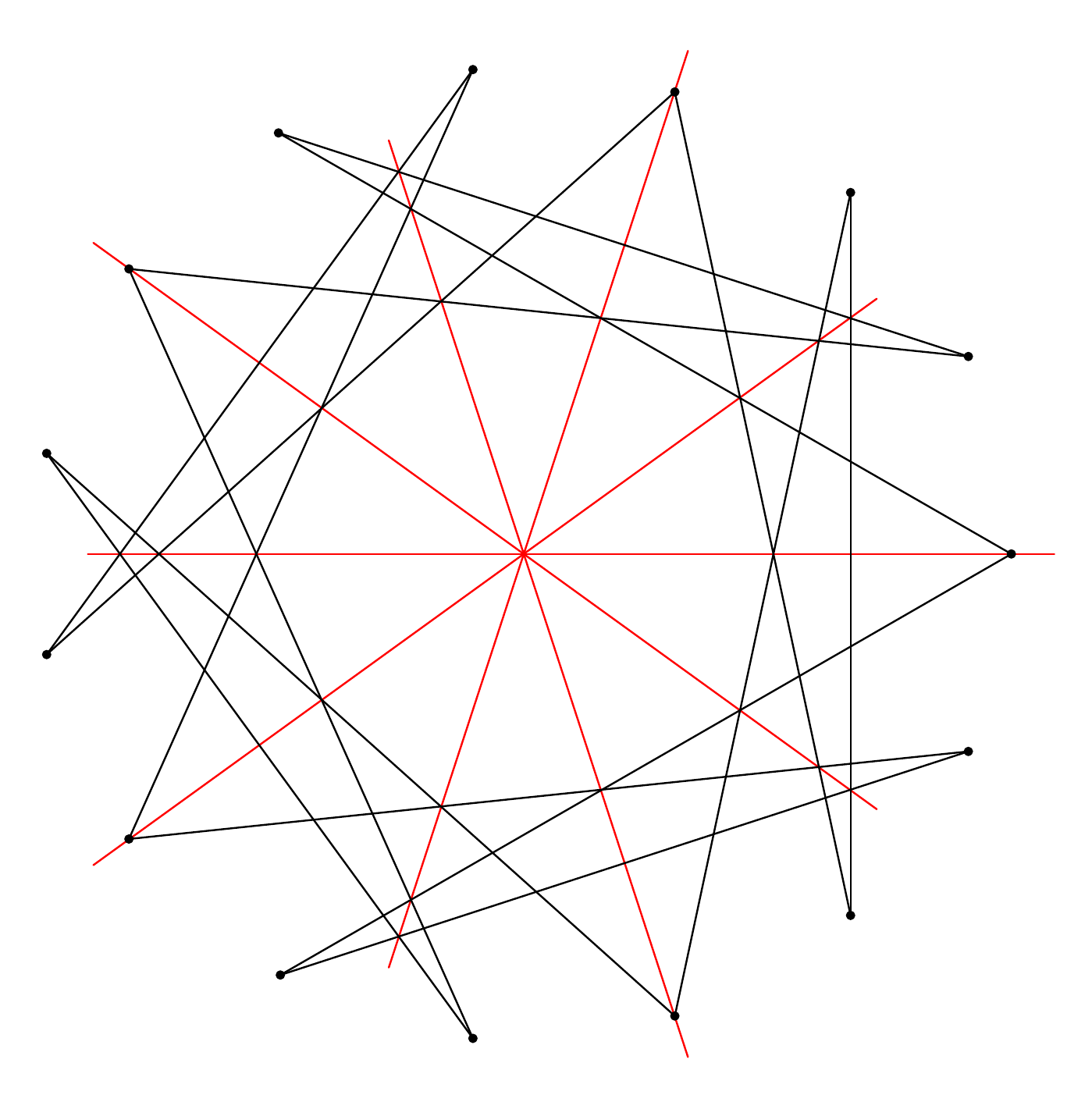} & \includegraphics[width=0.2\textwidth]{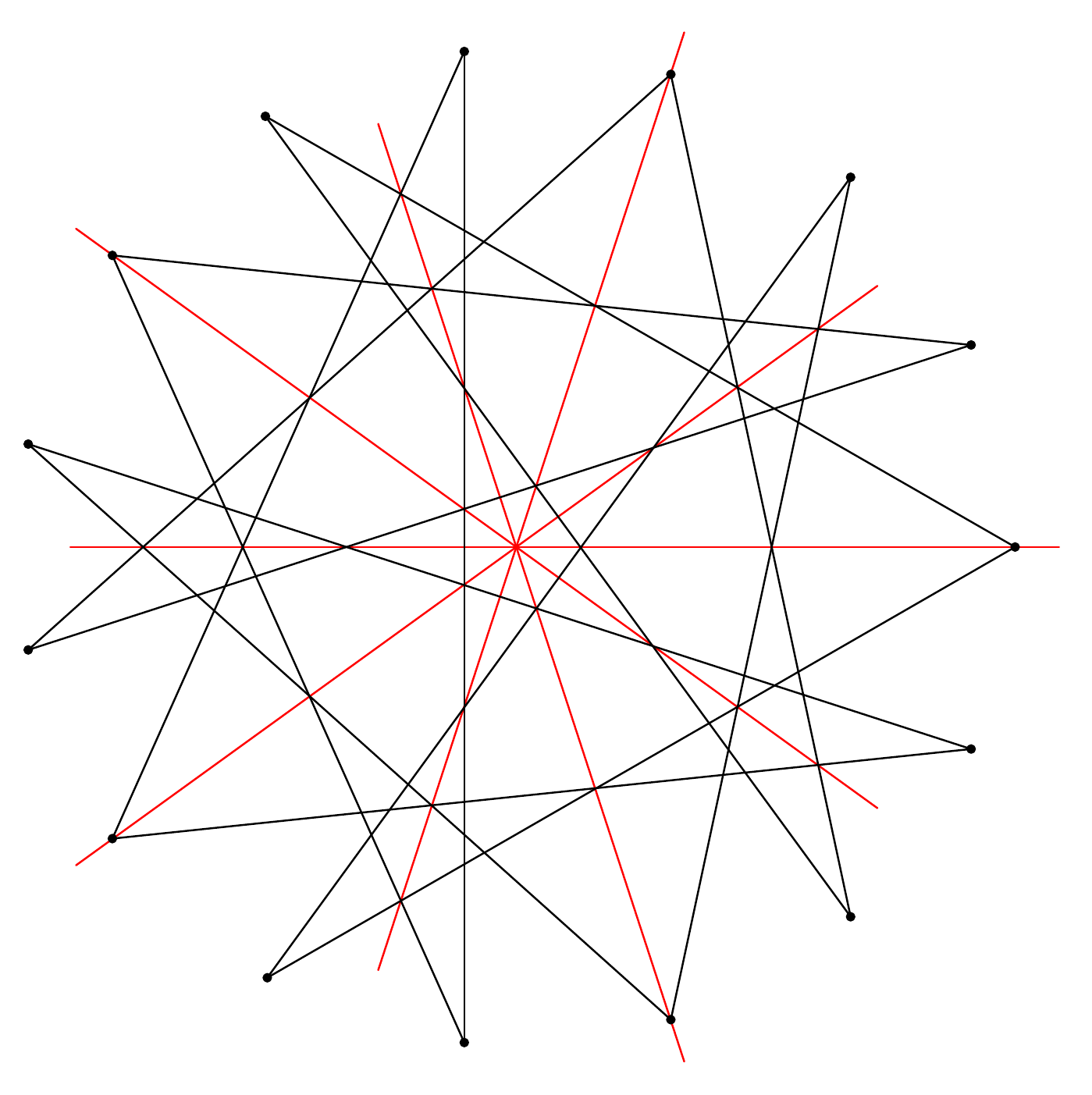}\\
$a=7;b=13$ & $a=10;b=1$ & $a=10;b=4$ & $a=10;b=7$\\ \hline
\includegraphics[width=0.2\textwidth]{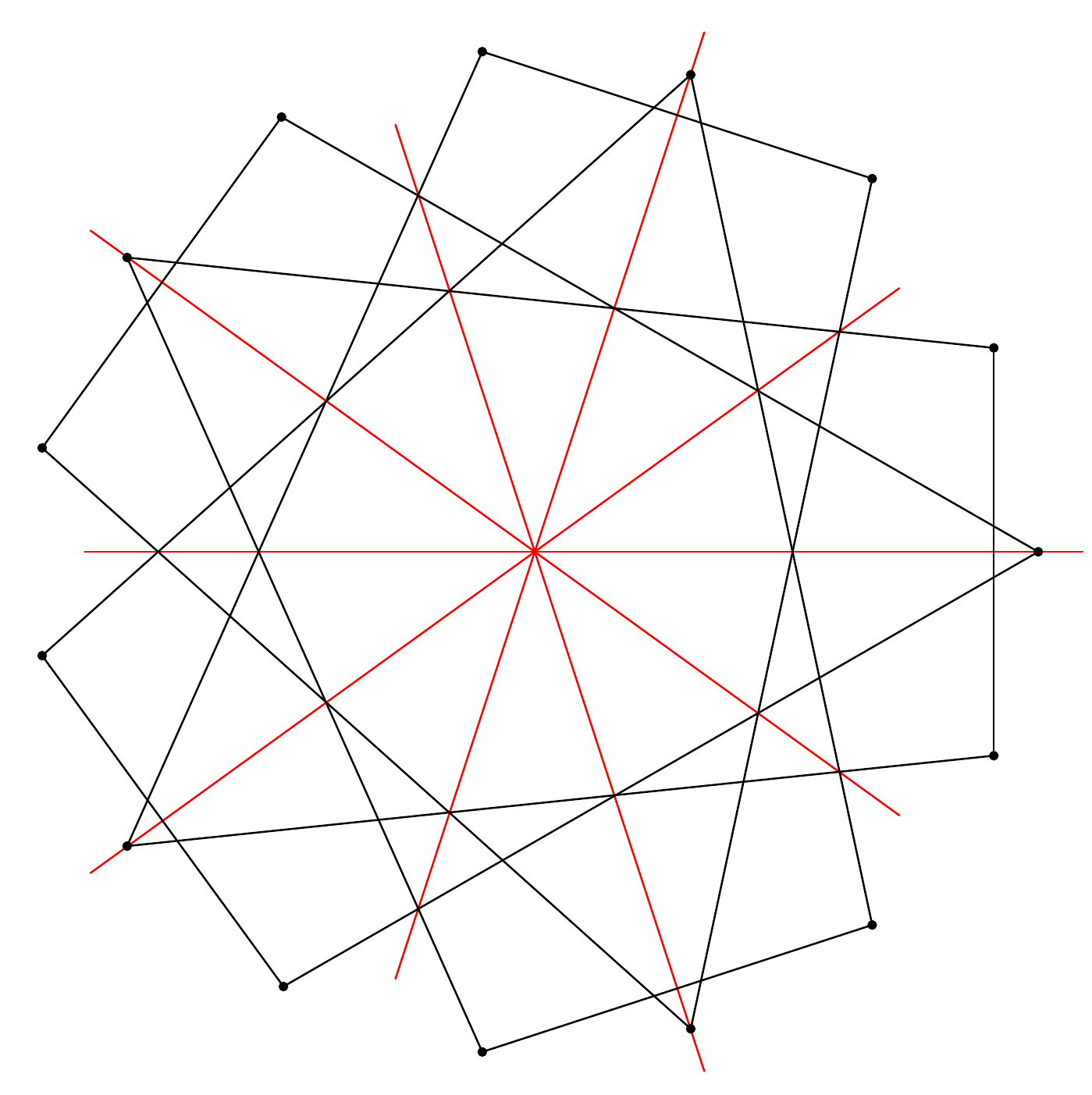} & \includegraphics[width=0.2\textwidth]{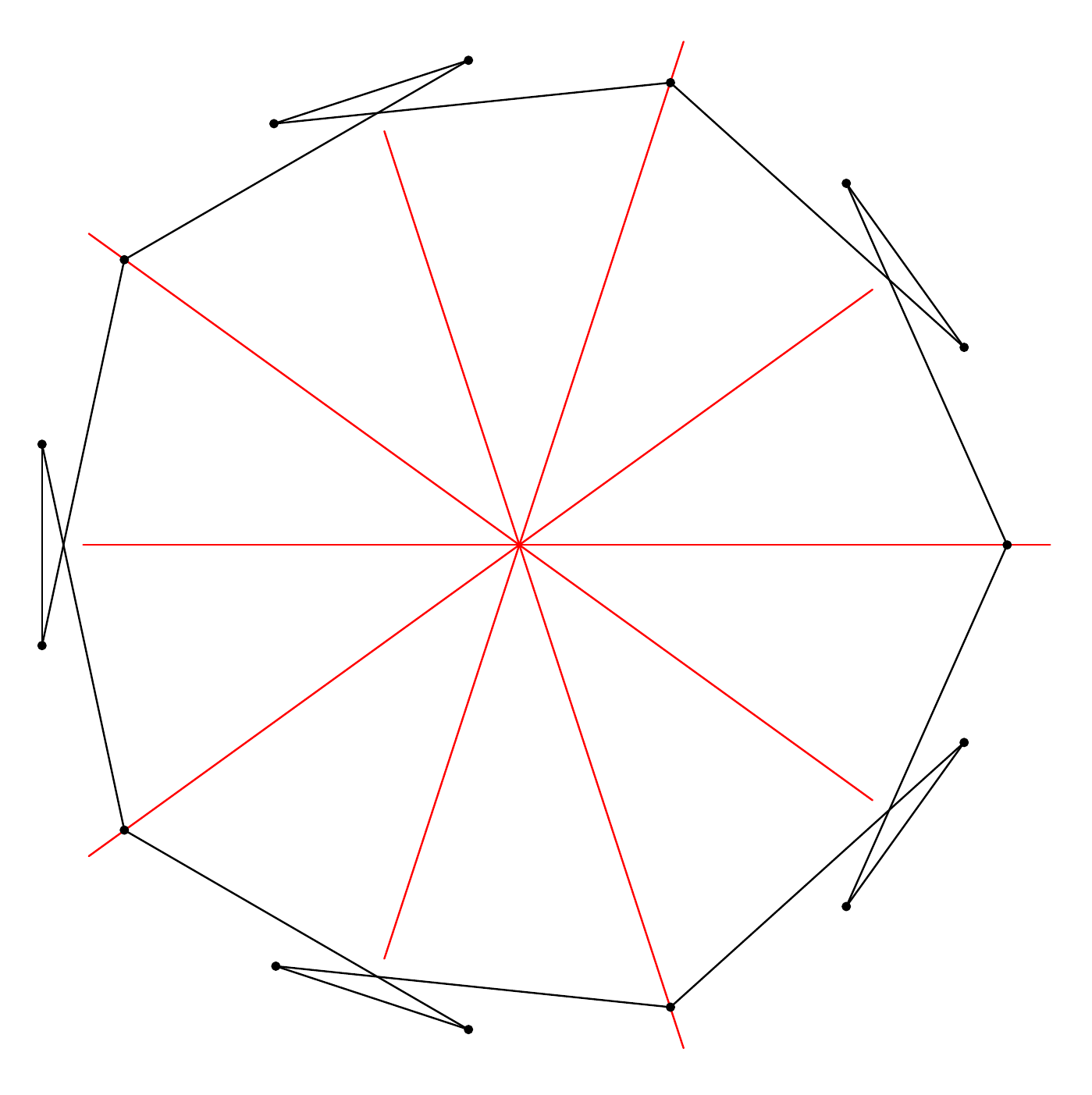} & \includegraphics[width=0.2\textwidth]{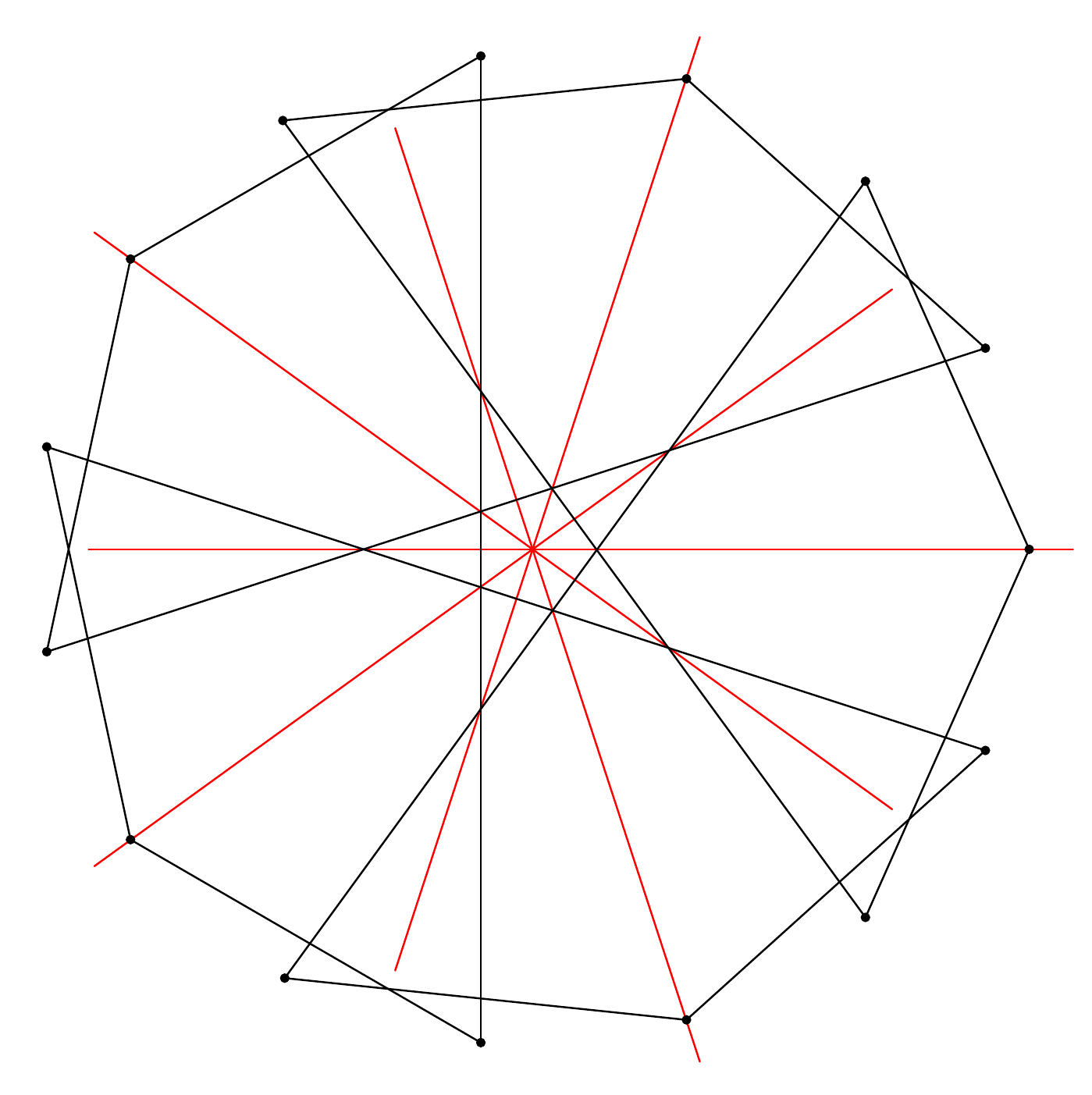} & \includegraphics[width=0.2\textwidth]{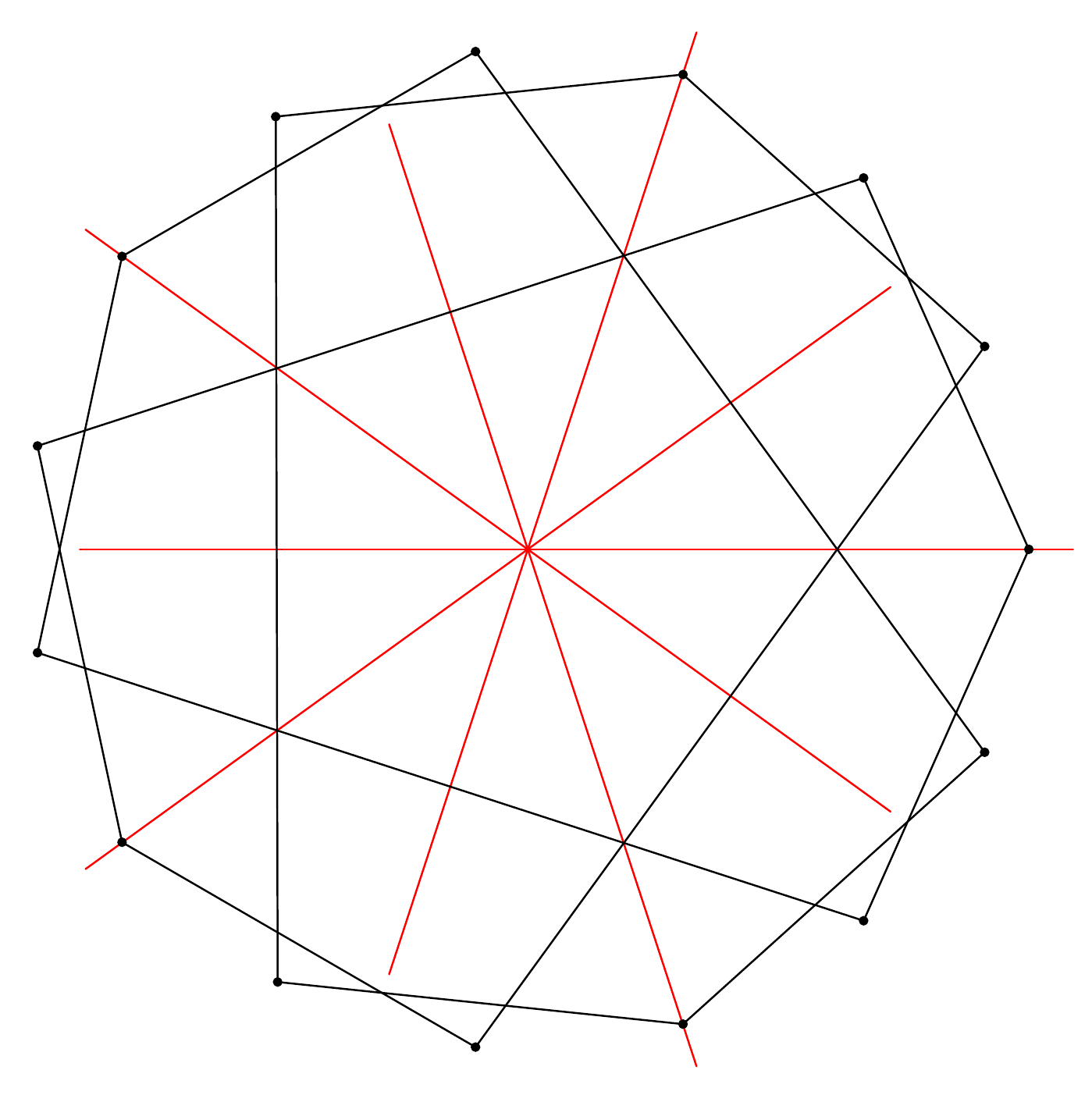}\\
$a=10;b=13$ & $a=13;b=1$ & $a=13;b=7$ & $a=13;b=10$\\
\end{tabular}
\caption{$n=15$: A set of representatives of the different equivalence-classes of $P_5(15)$}
\end{figure}
\newpage
\subsection{$m$-circular $3m$-polygons}
\label{subsec:m-circular_3m-polygons}
\subsubsection{Main-theorem 2}
\label{subsubsec:main_theorem_2}
Let $m>2$ be an integer and $n=3m$.\newline
The different equivalence classes of $m$-circular $3m$-polygons are represented by the $n$-tuples $(a, b, c, a, b, c, \ldots, a, b, c)$ of their sides, if $a$, $b$ and $c$ have the following ten properties:
\begin{enumerate}
\item $a\in \mathbb{N}$ with $a\equiv $1 mod 3,
\item $b\in \mathbb{N}$ with $b\equiv $1 mod 3,
\item $c\in \mathbb{N}$ with $c\equiv $1 mod 3,
\item $gcd\left(u,m \right)=1$,
\item $1 \leq a\leq n-2$,
\item $1 \leq b\leq n-2$,
\item $1 \leq c\leq n-2$,
\item $a \neq b$,
\item $b \neq c$,
\item $c \neq a$.
\end{enumerate}
\subsubsection{Formula for $\vert Q_m(3m)\vert$}
\label{subsubsec:formula}
Let $m>2$ be an integer and $n=3m$ and $\varphi(m)$ denote the Euler $\varphi$-function:\\

The number $\vert Q_m(3m)\vert$ of the equivalence-classes of the $m$-circular $3m$-polygons is:
\begin{center}
$\vert Q_m(3m)\vert=\dfrac{\varphi(m)\cdot m \cdot (m-3)+\varphi(3m)}{3}$.
\end{center}
\begin{table}[!htp]
\centering
\begin{tabular}{| c | c | c || c | c | c || c | c | c || c | c | c |}
\hline
n & m & $\vert Q_m(n) \vert$ & n & m & $\vert Q_m(n) \vert$ & n & m & $\vert Q_m(n) \vert$ & n & m & $\vert Q_m(n) \vert$\\ \hline
9 & 3 & 2 & 12 & 4 & 4 & 15 & 5 & 16 & 18 & 6 & 14 \\
21 & 7 & 60 & 24 & 8 & 56 & 27 & 9 & 114 & 30 & 10 & 96\\
33 & 11 & 300 & 36 & 12 & 148 & 39 & 13 & 528 & 42 & 14 & 312\\
45 & 15 & 488 & 48 & 16 & 560 & 51 & 17 & 1280 & 54 & 18 & 546\\
57 & 19 & 1836 & 60 & 20 & 912 & 63 & 21 & 1524 & 66 & 22 & 1400\\
69 & 23 & 3388 & 72 & 24 & 1352 & 75 & 25 & 3680 & 78 & 26 & 2400\\
81 & 27 & 3906 & 84 & 28 & 2808 & 87 & 29 & 7056 & 90 & 30 & 2168\\
\hline
\end{tabular}
\caption{Number of equivalence classes of $m$-circular $3m$-polygons}
\label{tab:number_of_equivalence_classes_of_m-circular 3m-polygons}
\end{table}
\newpage
\subsubsection{Examples of $m$-circular $3m$-polygons}
\label{examples of m-circular 3m-polygons}
\paragraph[{$n=9$; $m=3$; $\vert Q_3(9) \vert=2$}]{$n=9$; $m=3$; $\vert Q_3(9) \vert=2$}.
\begin{figure}[!htp]
\centering
\begin{tabular}{c | c }
\includegraphics[width=0.4\textwidth]{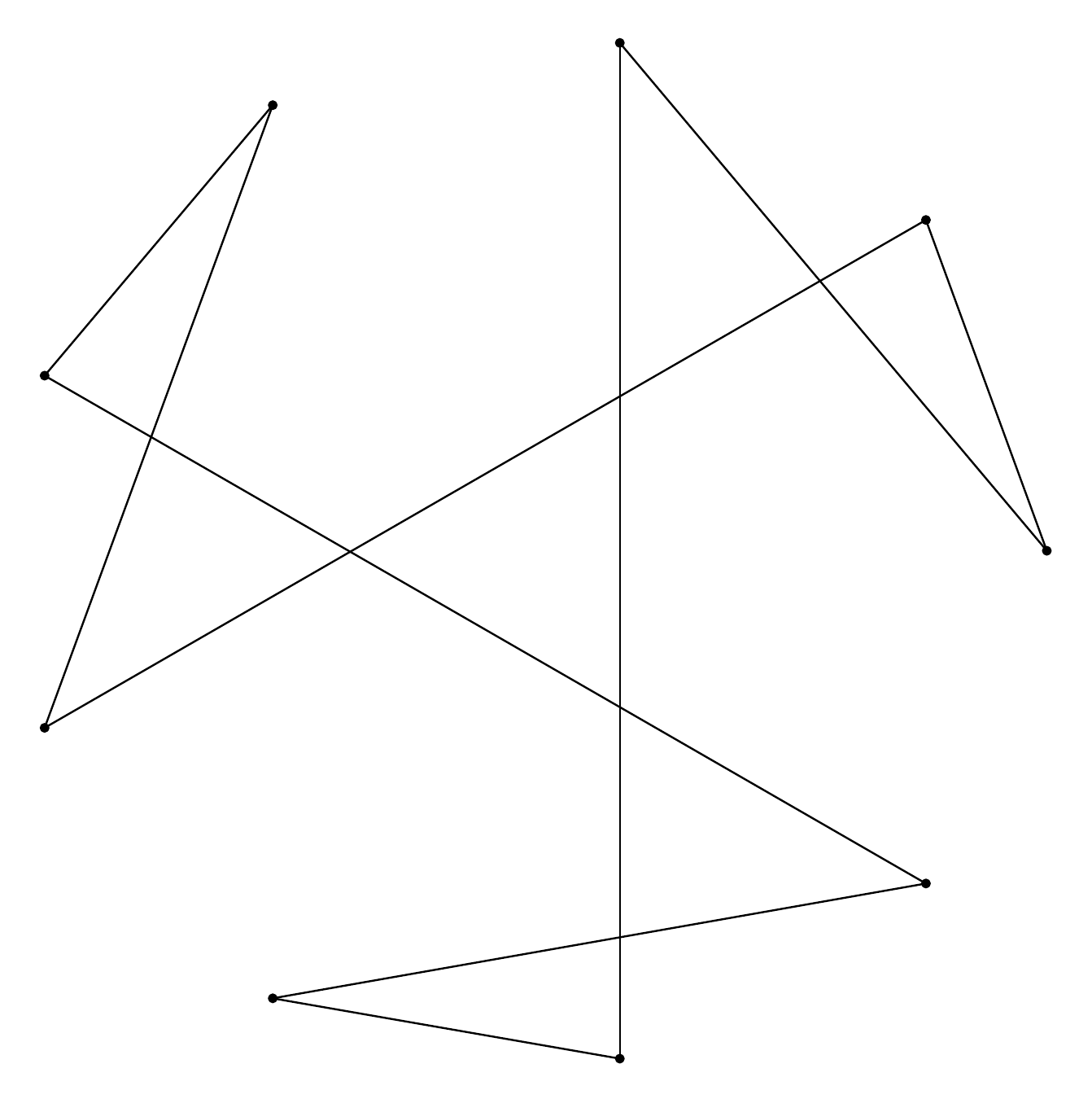} & \includegraphics[width=0.4\textwidth]{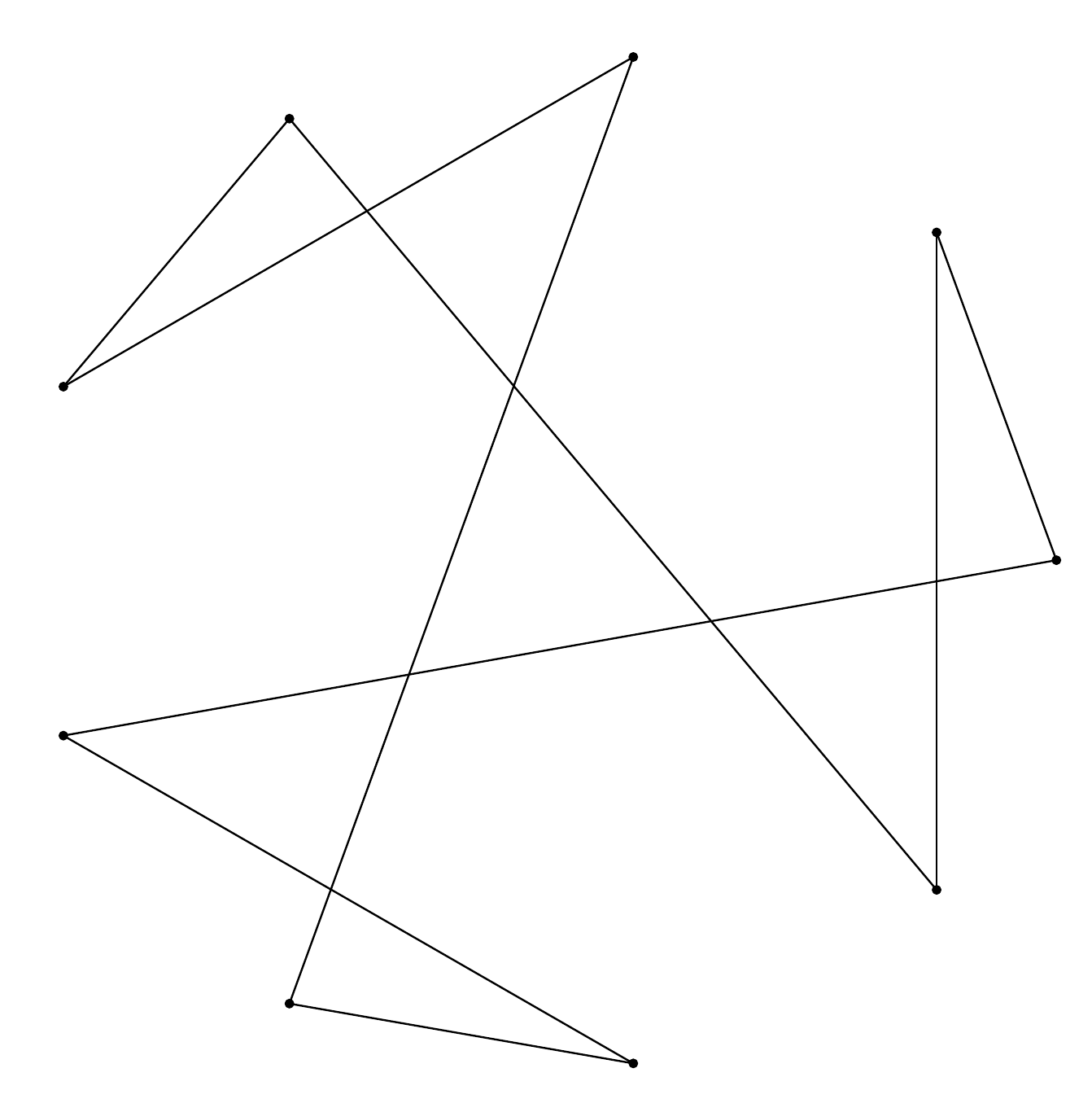}\\
$a=1;b=4;c=7$ & $a=1;b=7;c=4$ 
\end{tabular}
\caption{$n=9$: A set of representatives of the two different equivalence-classes of $Q_3(9)$}
\end{figure}
\paragraph[$n=12$; $m=4$; $\vert Q_4(12) \vert=4$]{$n=12$; $m=4$; $\vert Q_4(12) \vert=4$}.
\begin{figure}[!htp]
\centering
\begin{tabular}{c | c }
\includegraphics[width=0.4\textwidth]{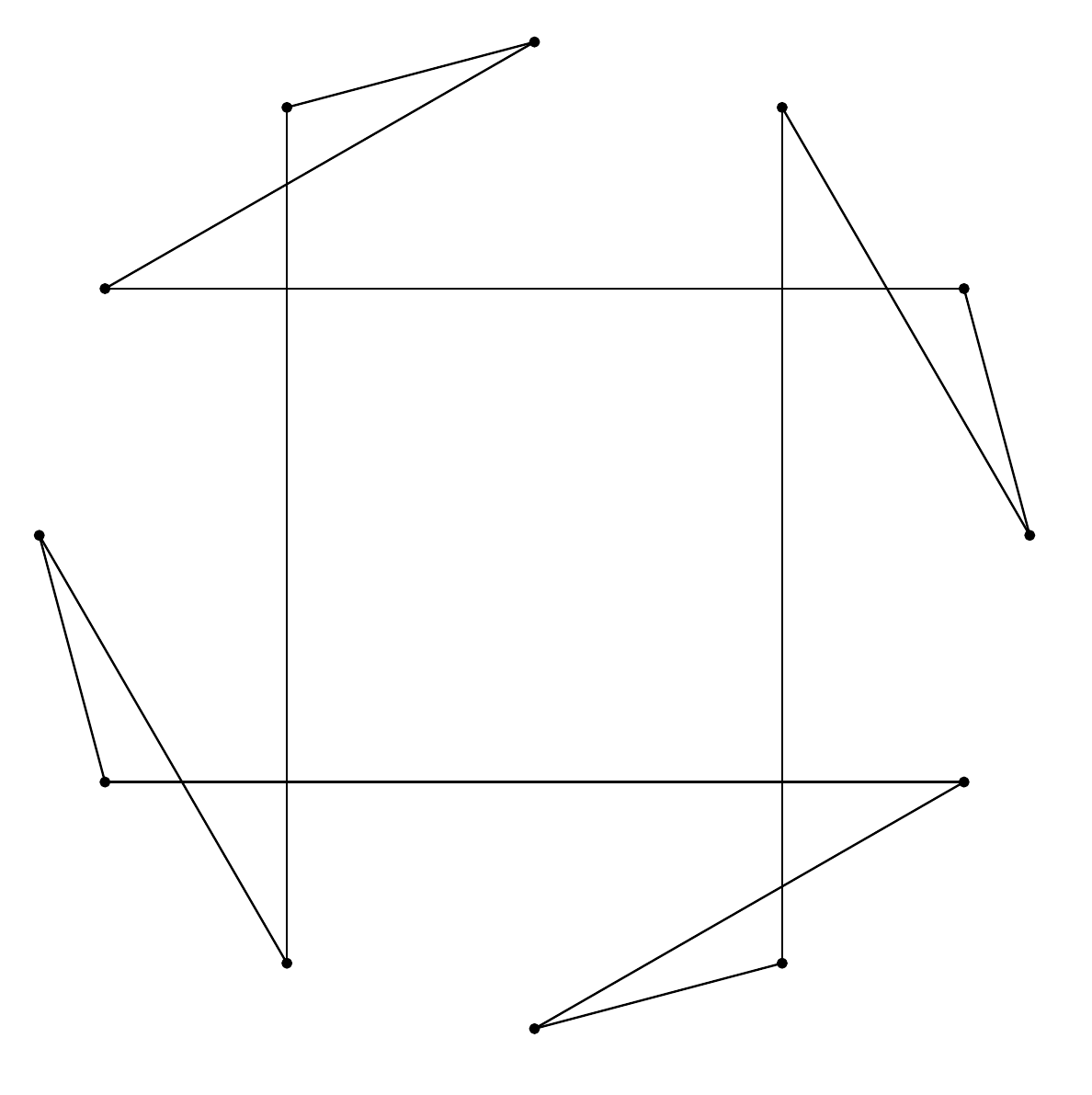} & \includegraphics[width=0.4\textwidth]{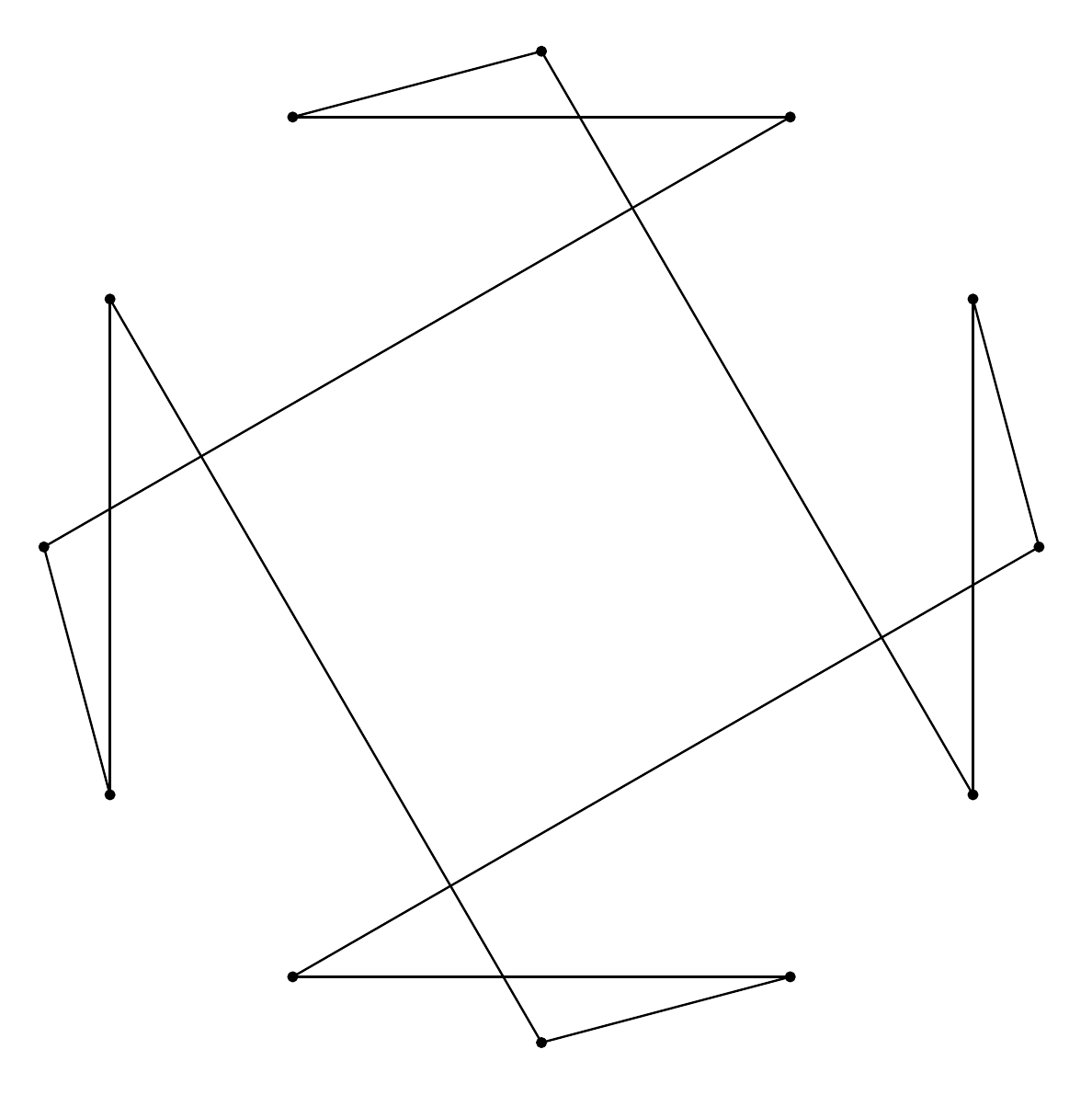}\\
$a=1;b=4;c=10$ & $a=1;b=10;c=4$\\ \hline
\includegraphics[width=0.4\textwidth]{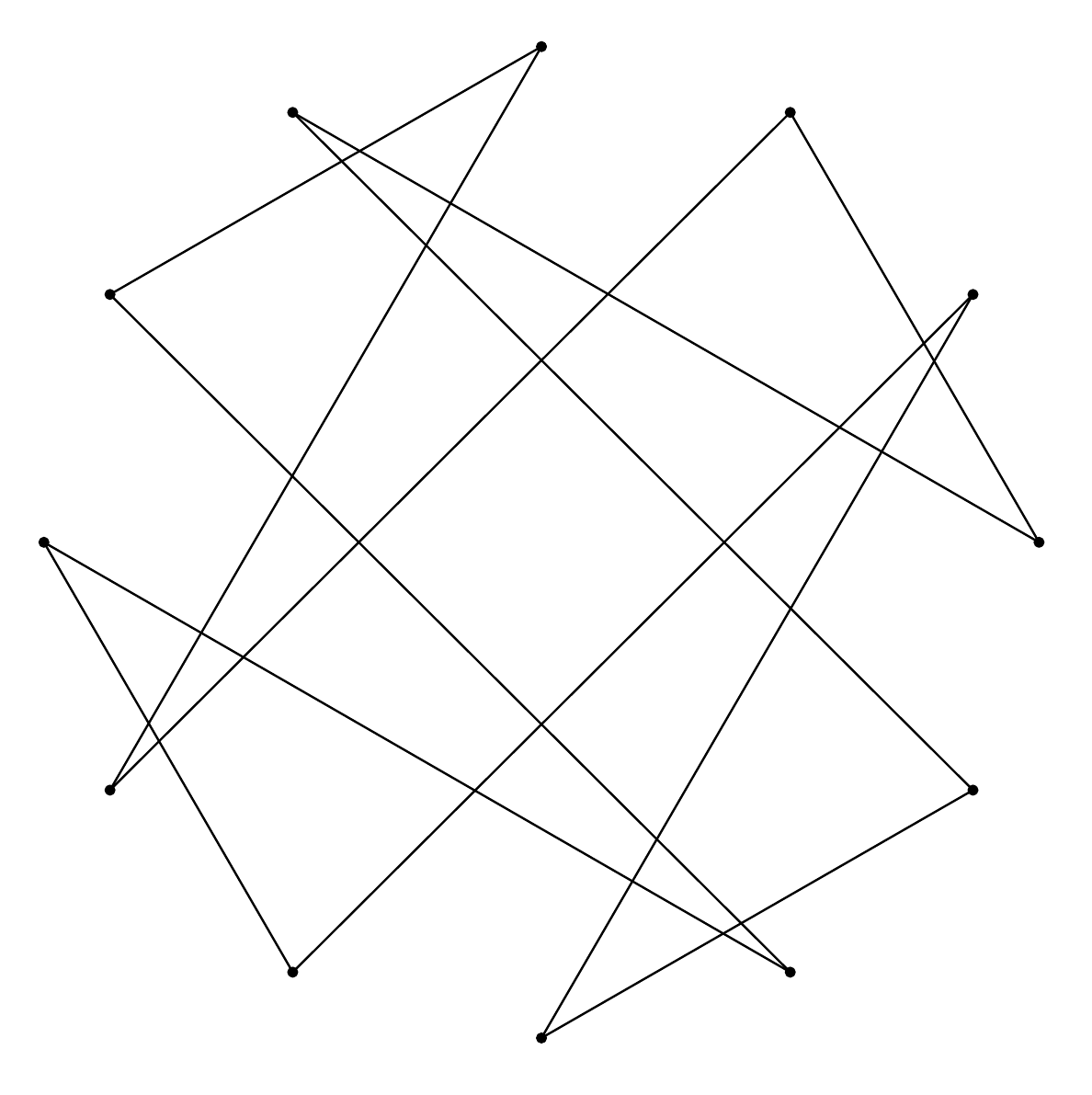} & \includegraphics[width=0.4\textwidth]{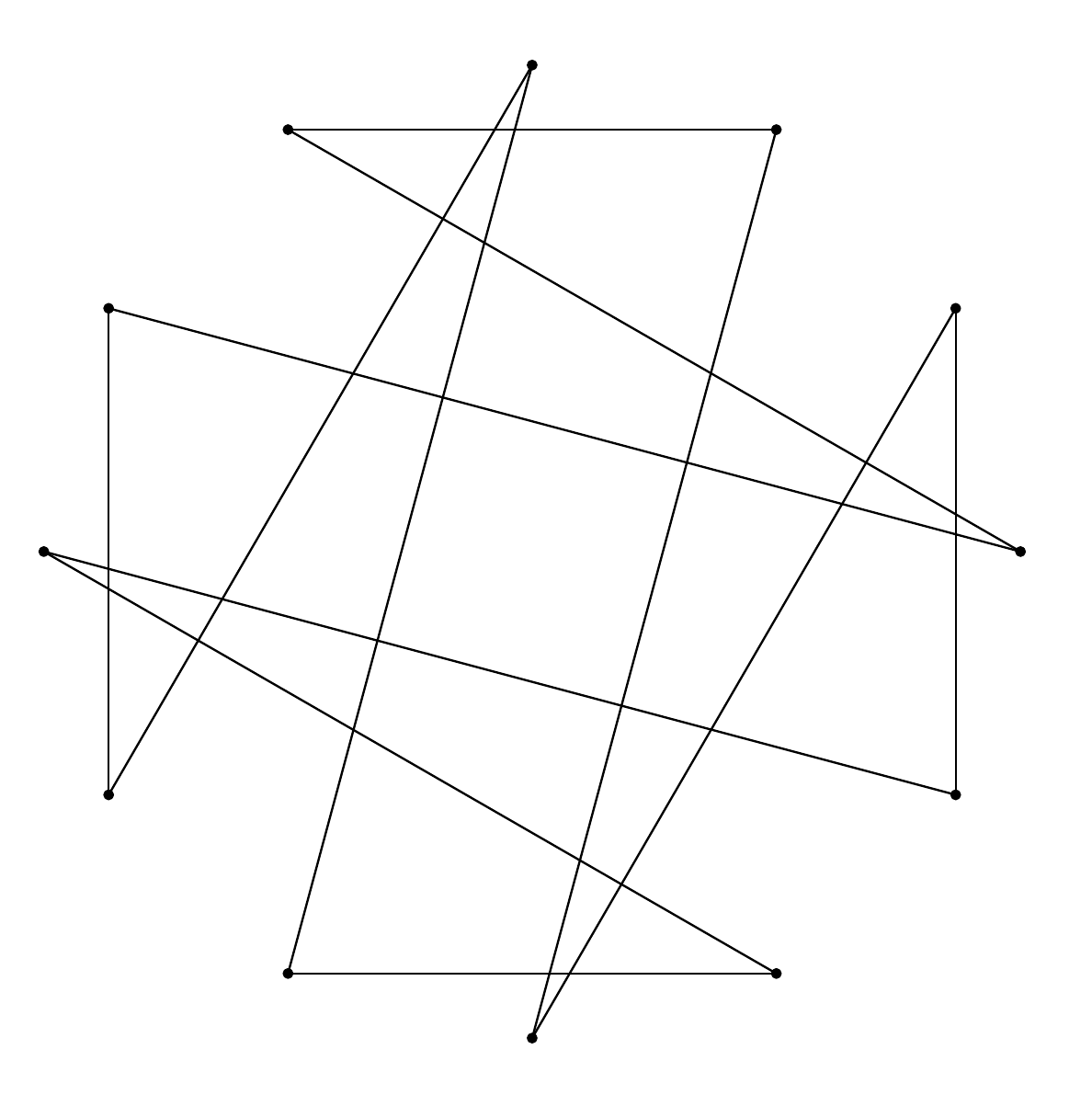}\\
$a=4;b=7;c=10$ & $a=4;b=10;c=7$\\
\end{tabular}
\caption{$n=12$: A set of representatives of the four different equivalence-classes of $Q_4(12)$}
\end{figure}
\newpage
\paragraph[$n=15$; $m=5$; $\vert Q_5(15) \vert=16$]{$n=15$; $m=5$; $\vert Q_5(15) \vert=16$}.
\begin{figure}[!htp]
\centering
\begin{tabular}{c | c | c | c }
\includegraphics[width=0.25\textwidth]{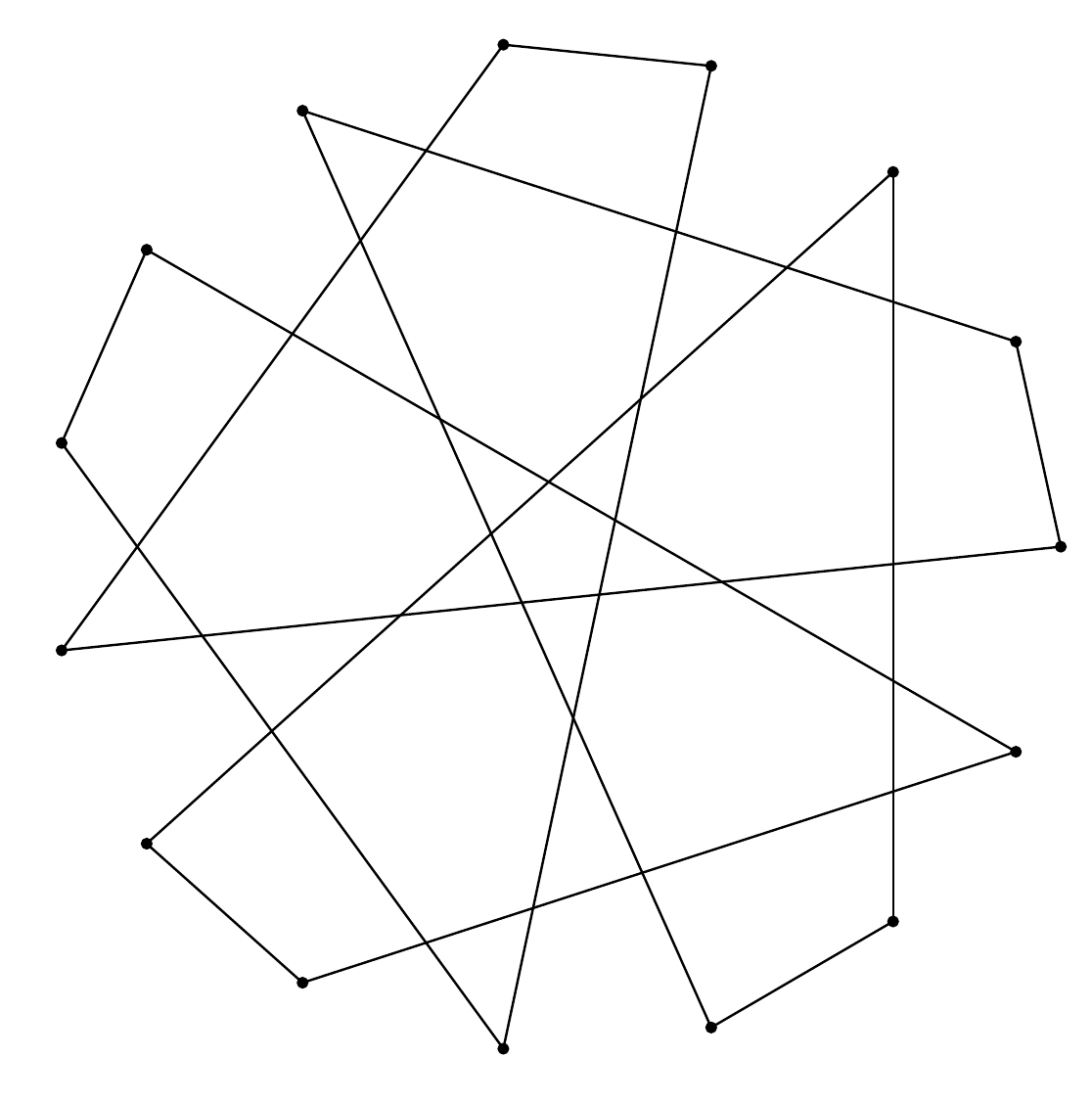} & \includegraphics[width=0.25\textwidth]{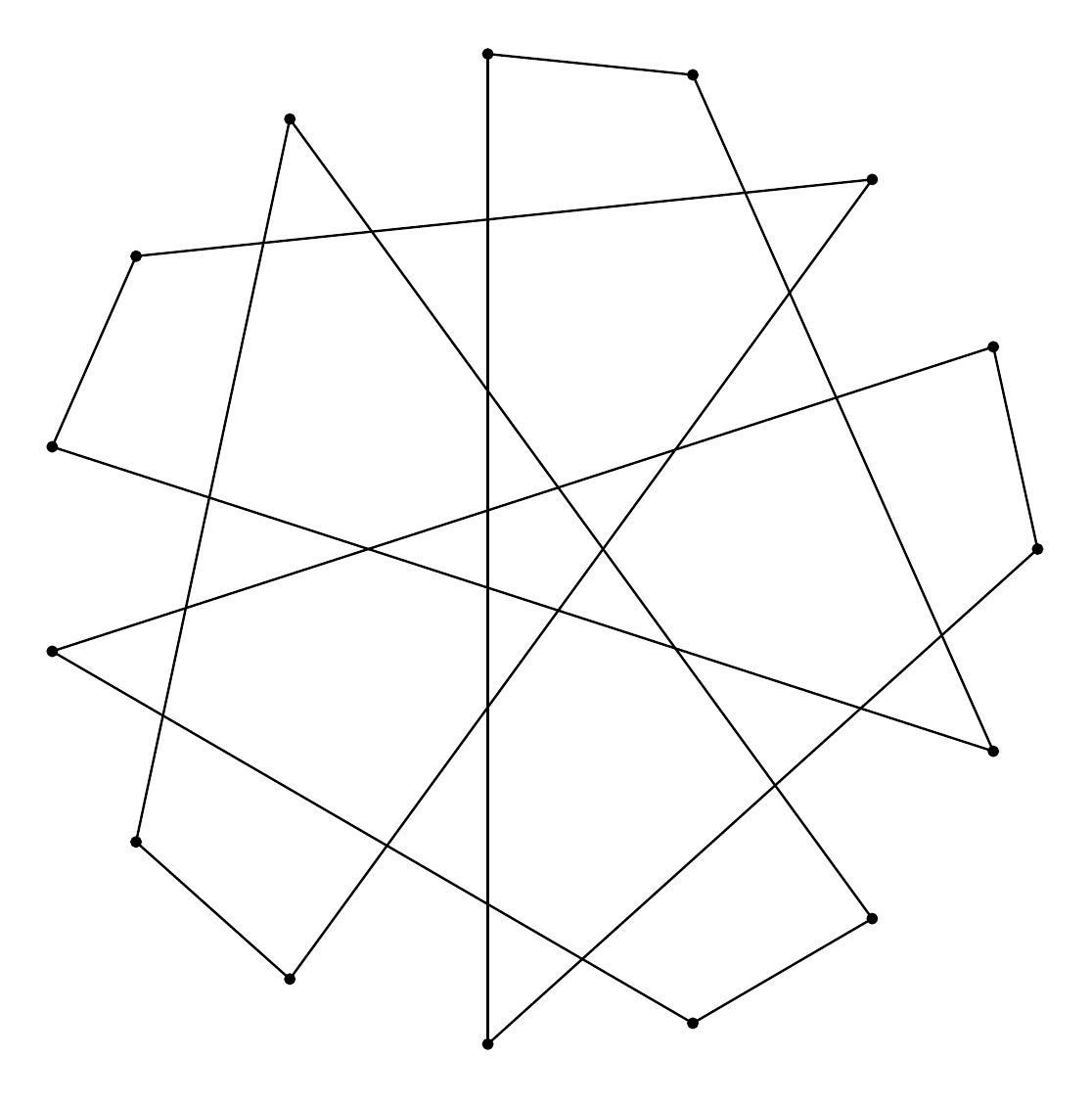} & \includegraphics[width=0.25\textwidth]{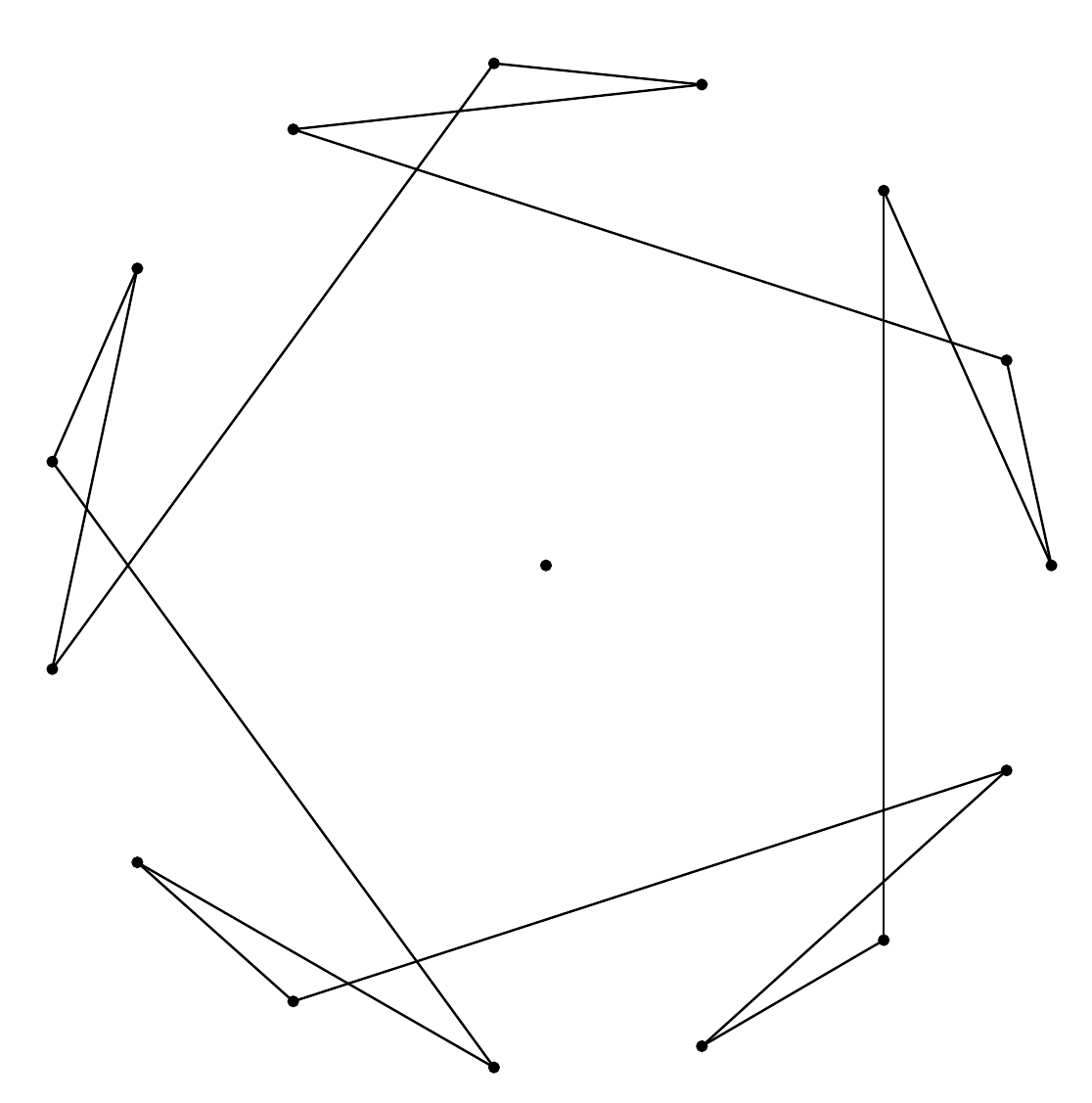} & \includegraphics[width=0.25\textwidth]{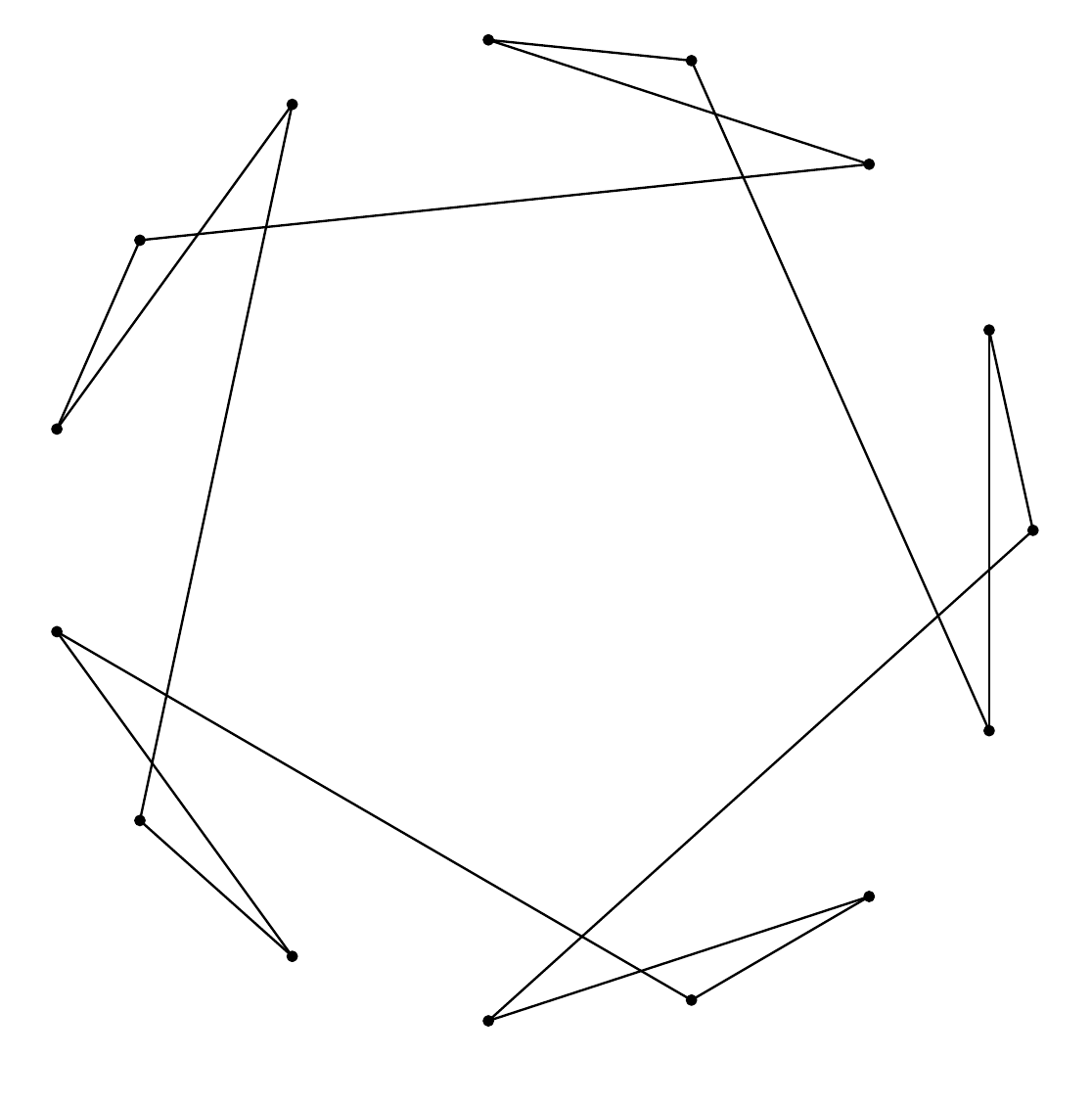}\\
$a=1;b=4;c=7$ & $a=1;b=7;c=4$ & $a=1;b=4;c=13$ & $a=1;b=13;c=4$\\ \hline
\includegraphics[width=0.25\textwidth]{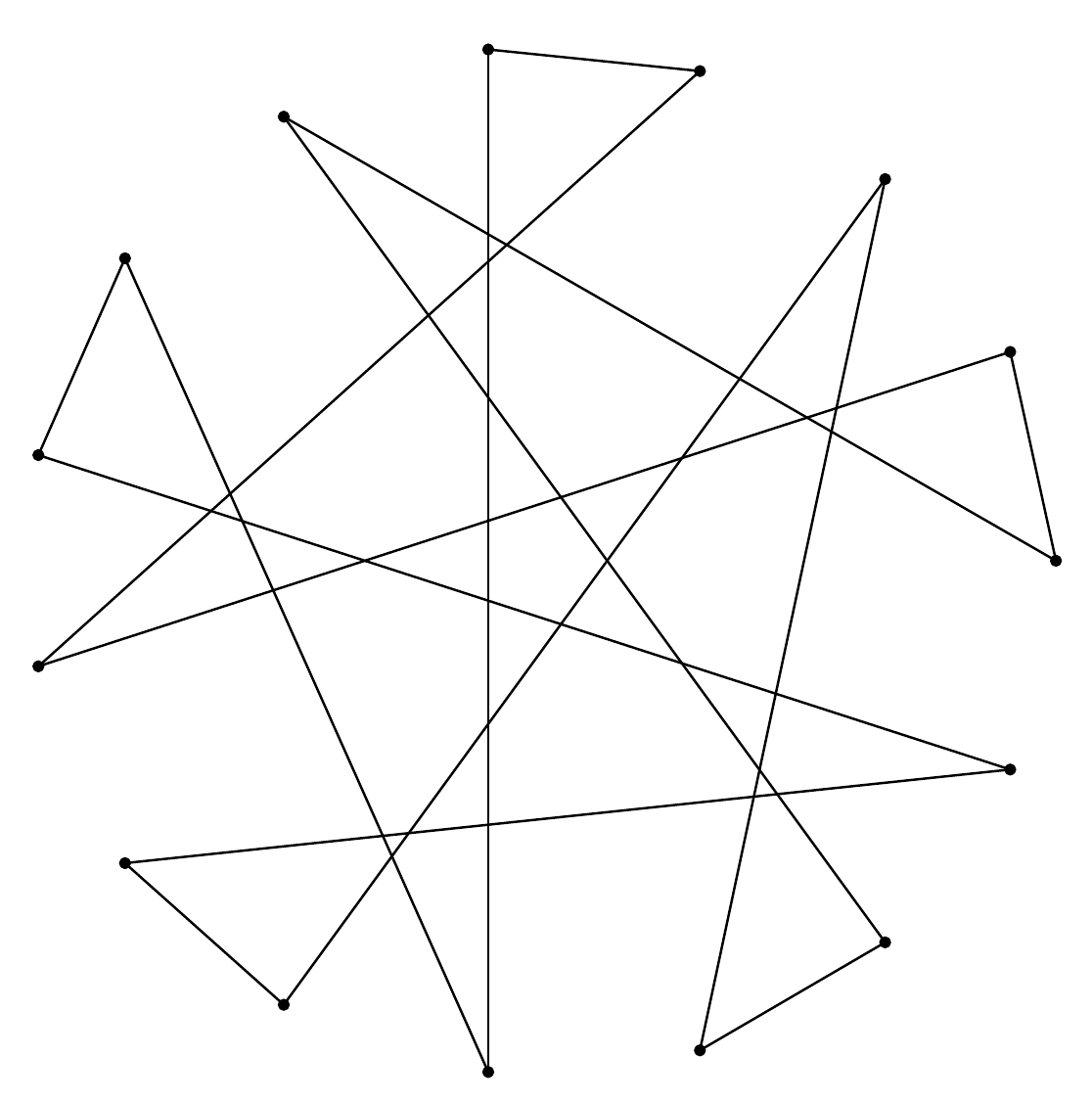} & \includegraphics[width=0.25\textwidth]{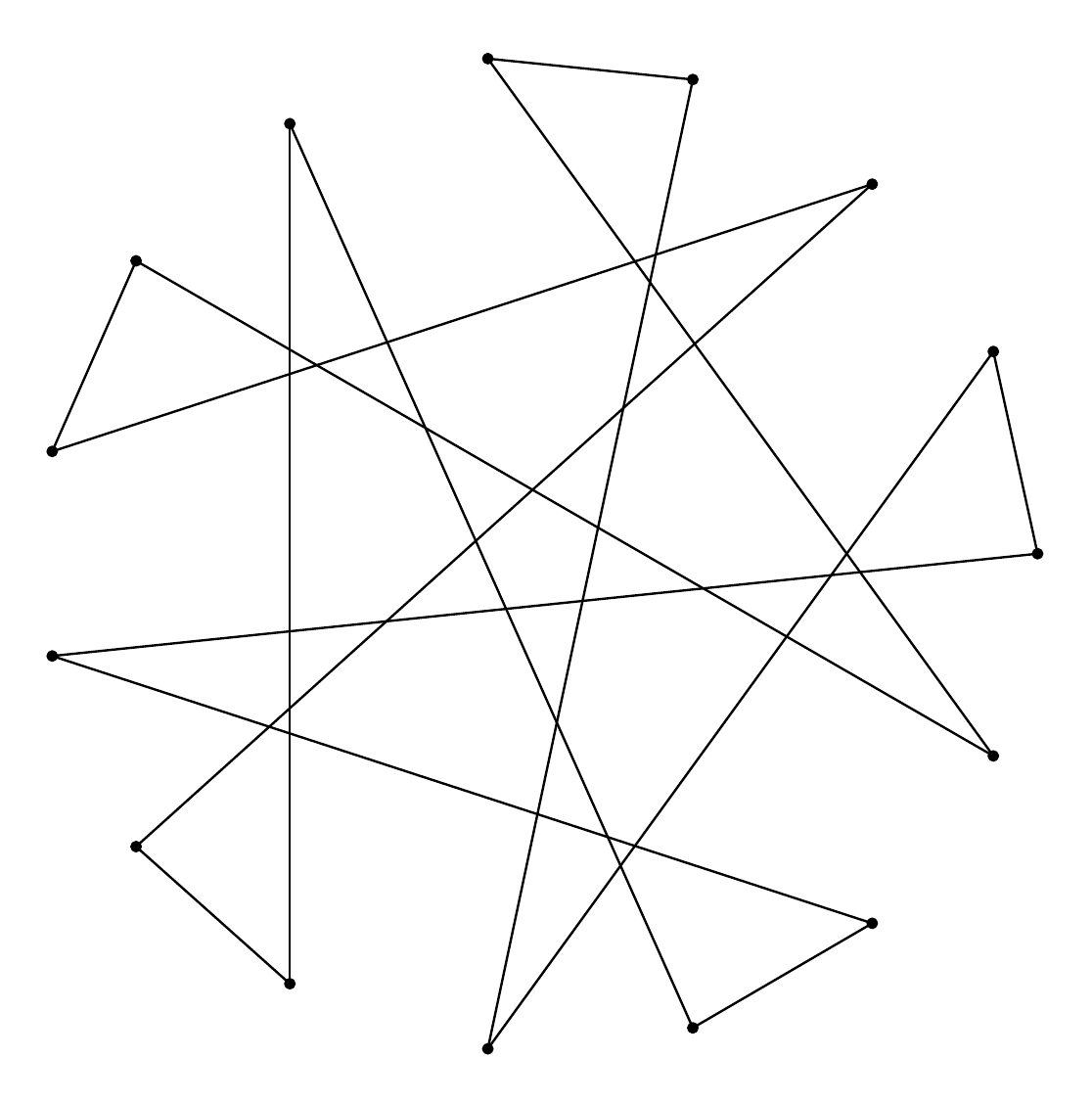} & \includegraphics[width=0.25\textwidth]{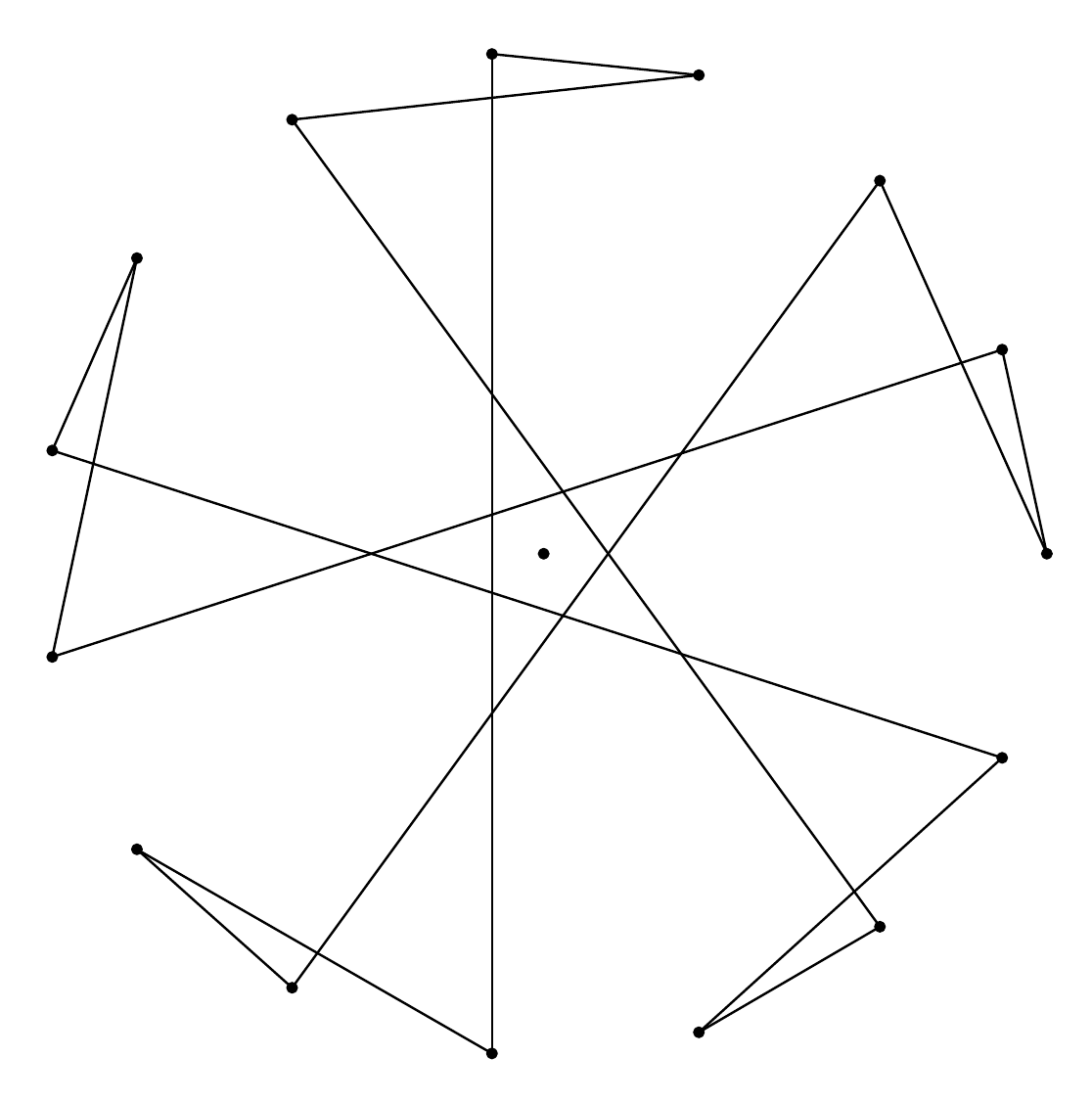} & \includegraphics[width=0.25\textwidth]{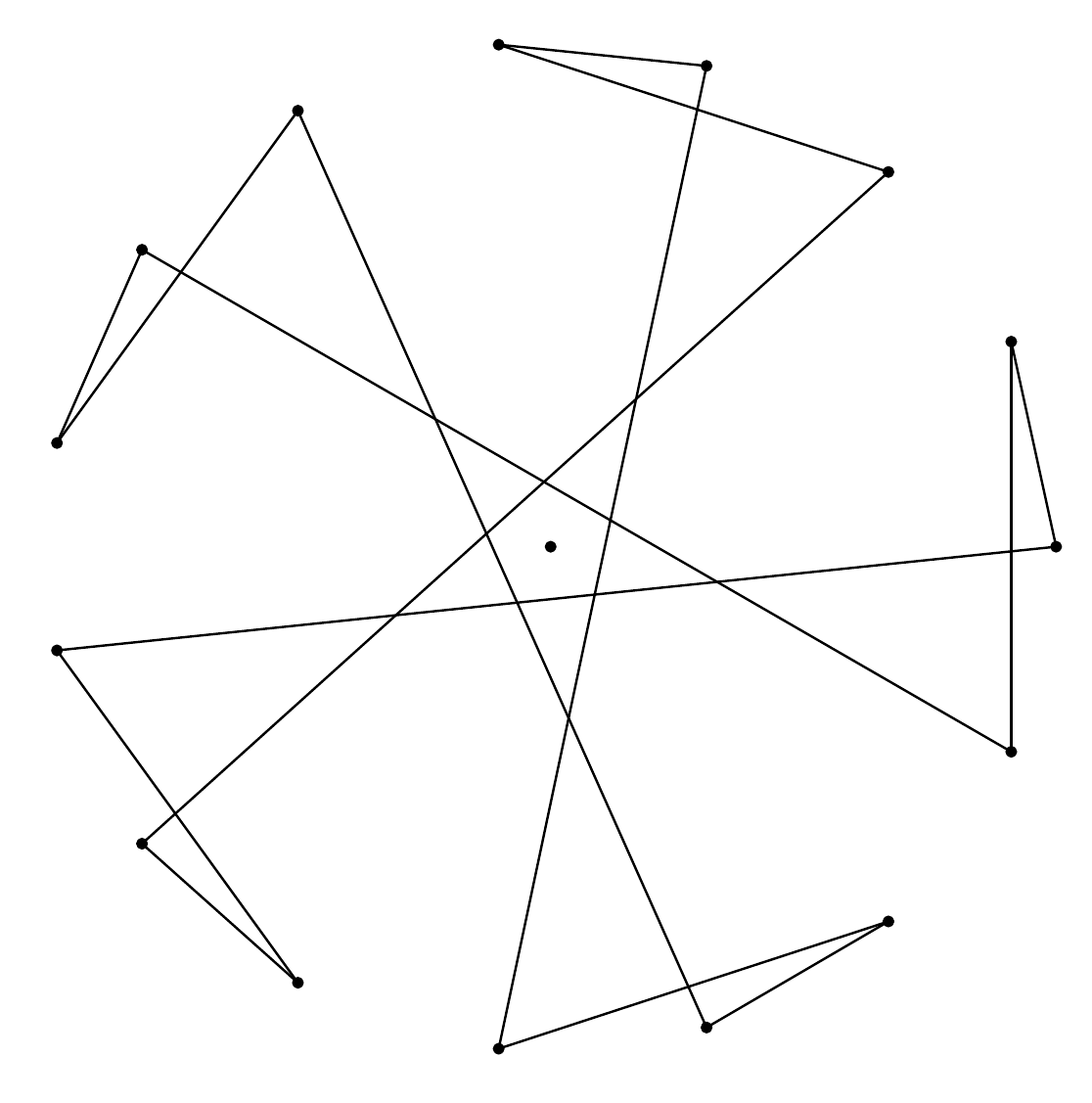}\\
$a=1;b=7;c=10$ & $a=1;b=10;c=7$ & $a=1;b=7;c=13$ & $a=1;b=13;c=7$\\ \hline
\includegraphics[width=0.25\textwidth]{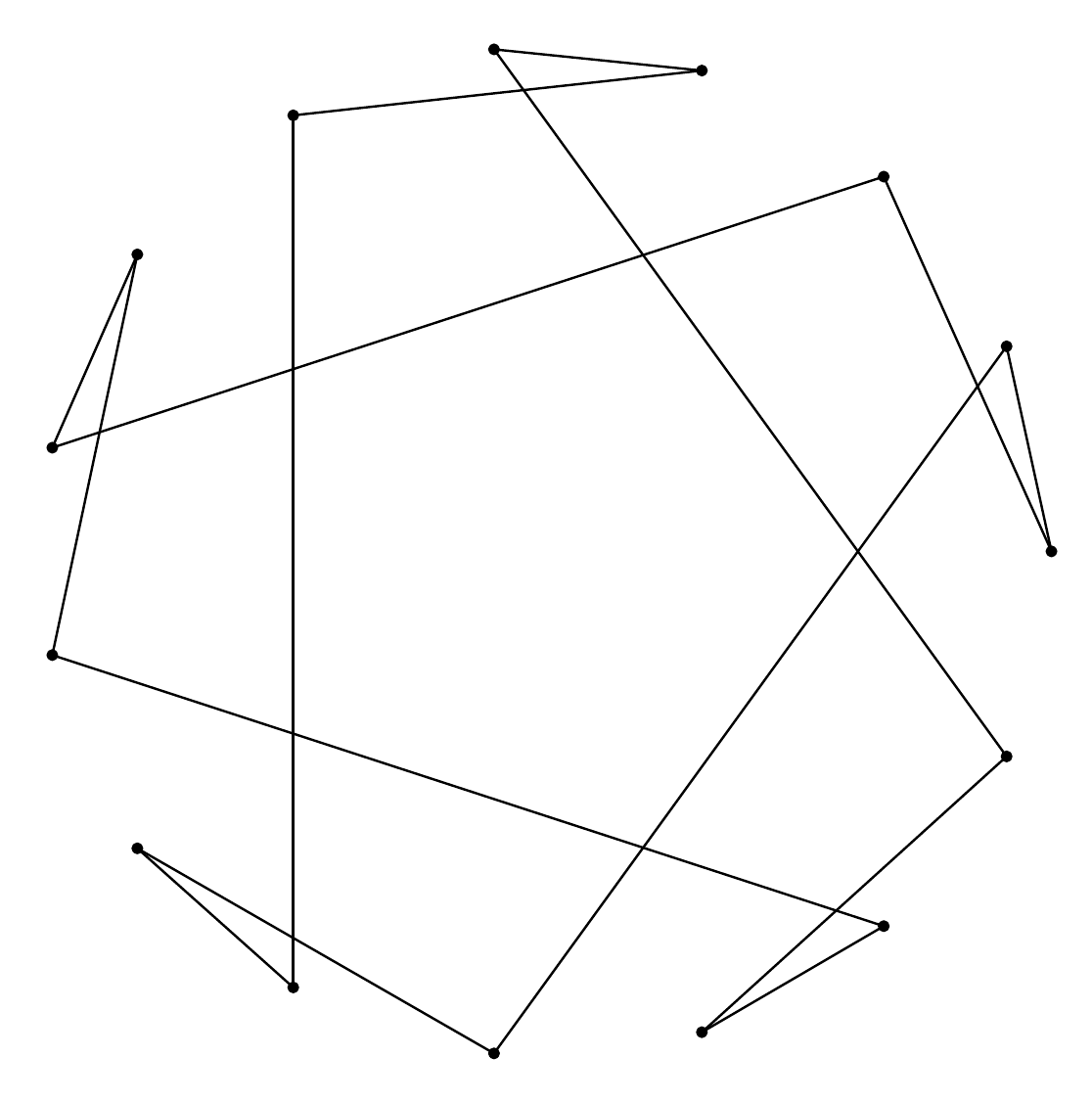} & \includegraphics[width=0.25\textwidth]{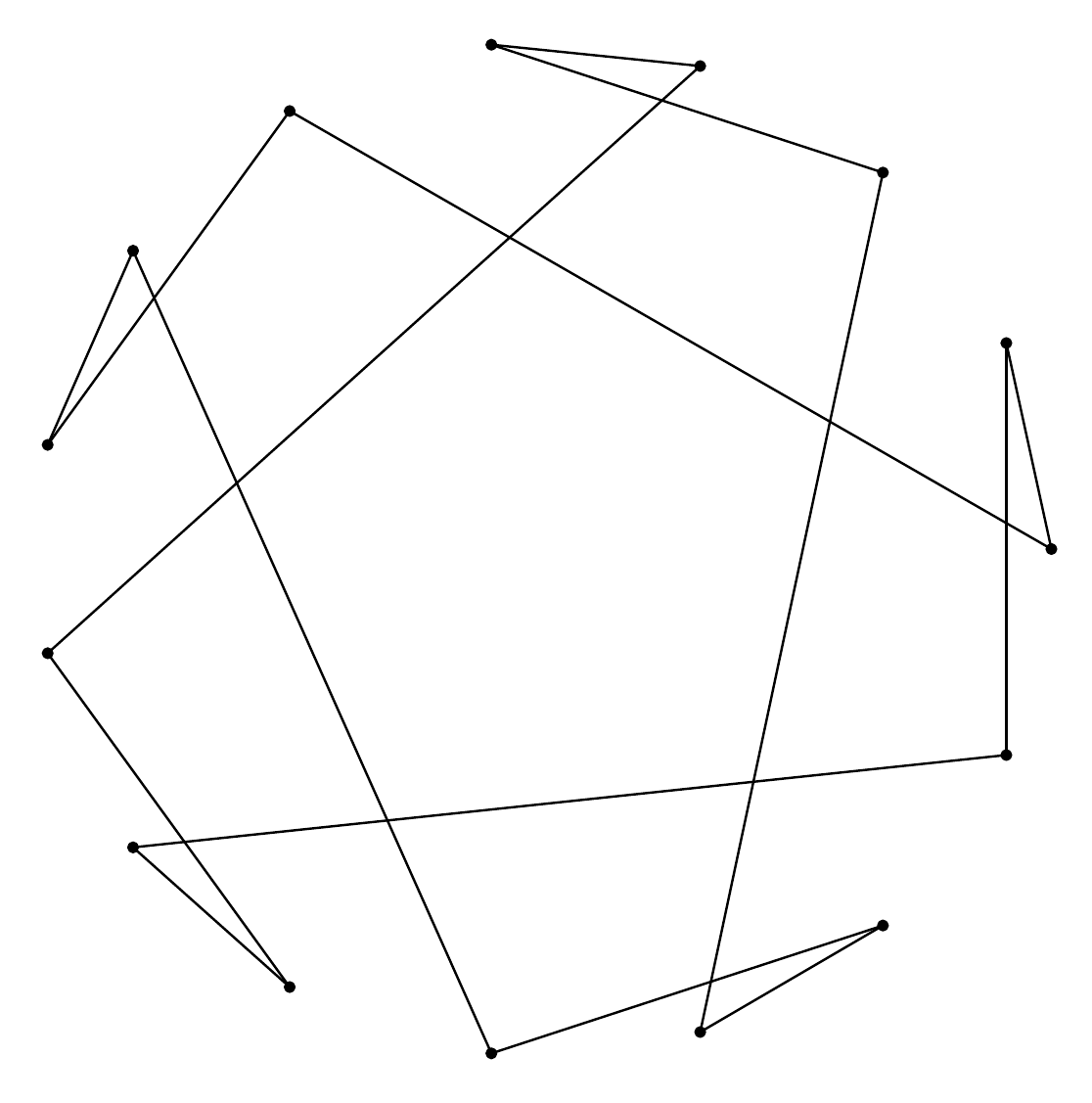} & \includegraphics[width=0.25\textwidth]{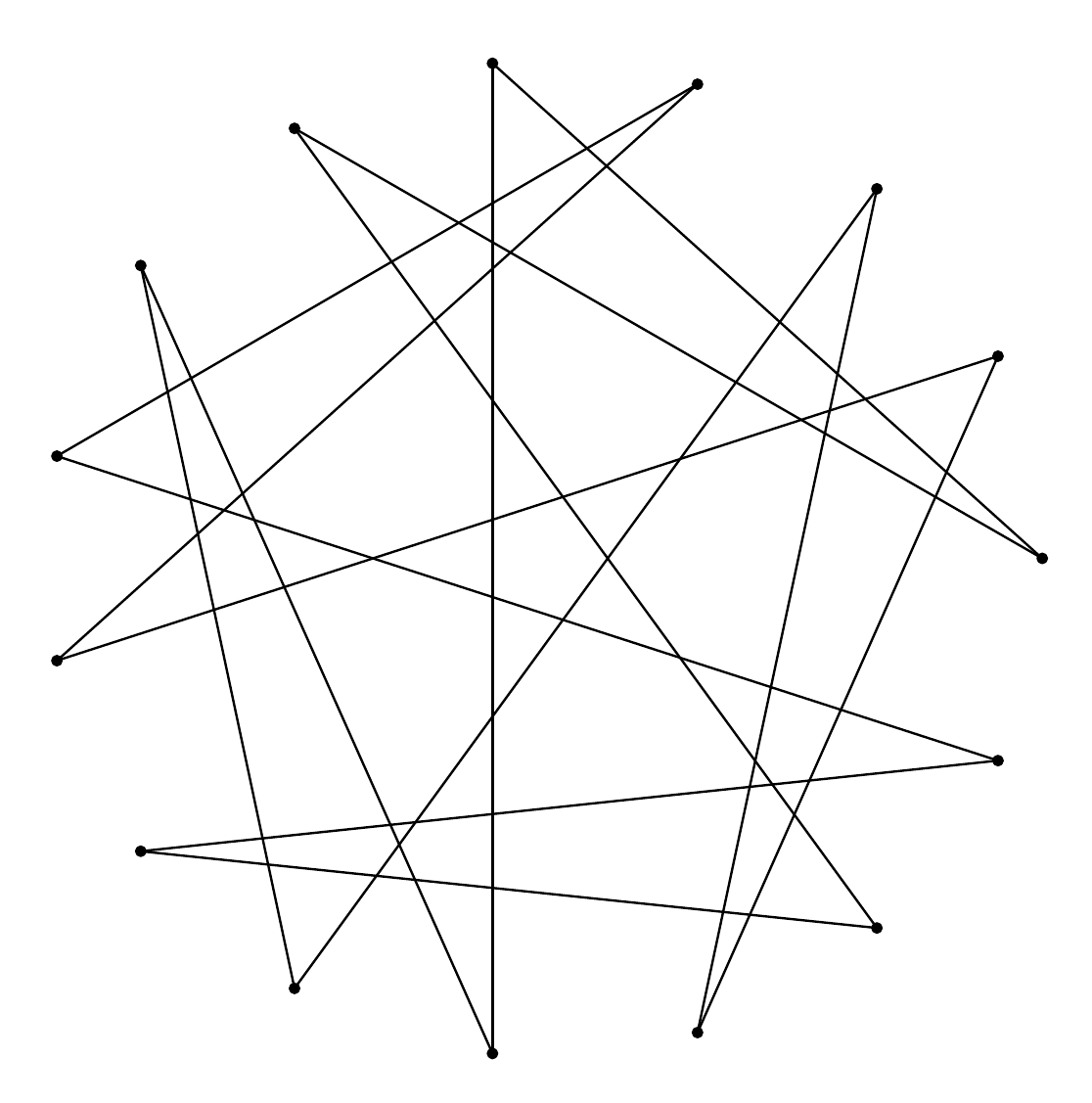} & \includegraphics[width=0.25\textwidth]{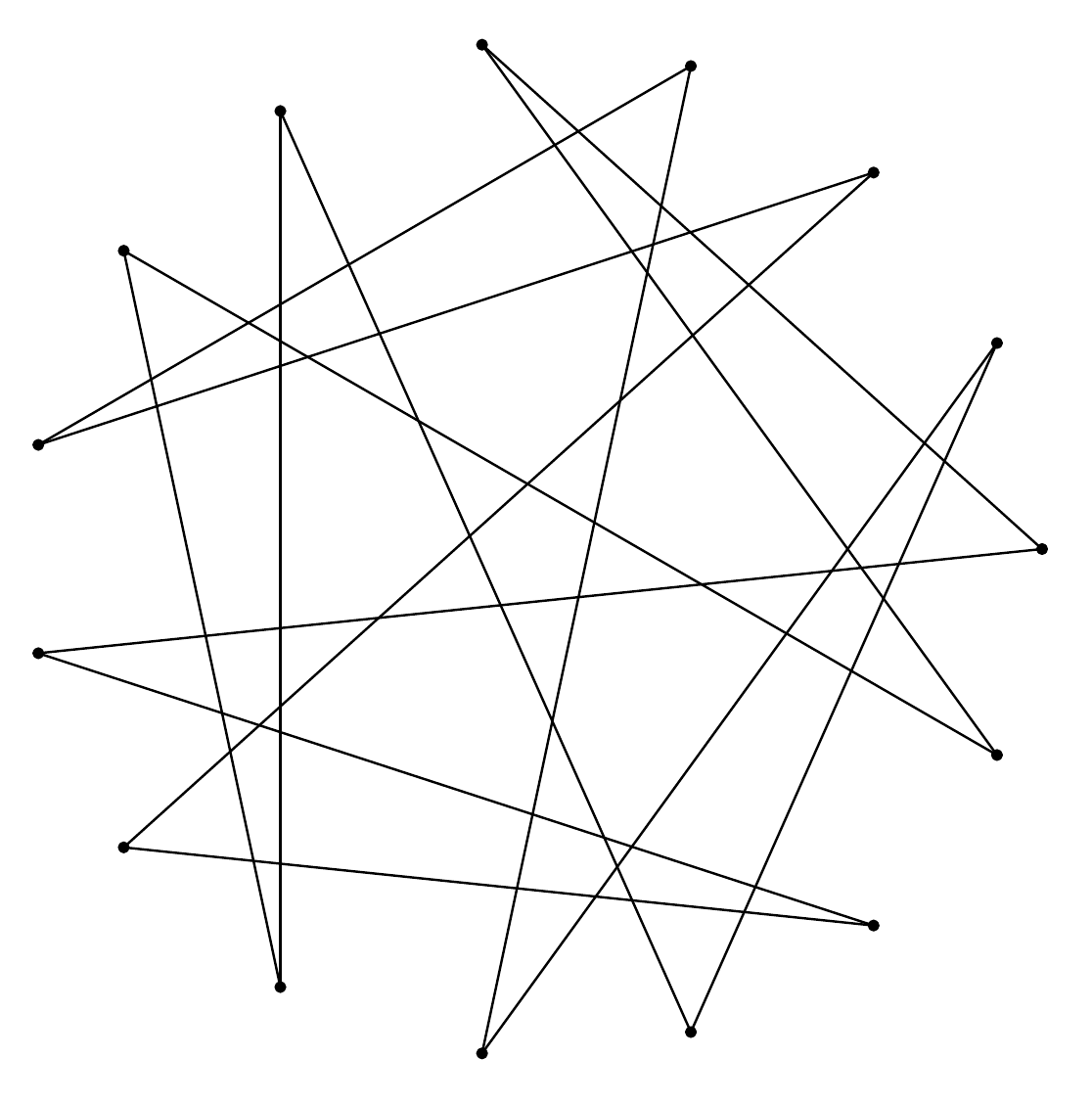}\\
$a=1;b=10;c=13$ & $a=1;b=13;c=10$ & $a=4;b=7;c=10$ & $a=4;b=10;c=7$\\ \hline
\includegraphics[width=0.25\textwidth]{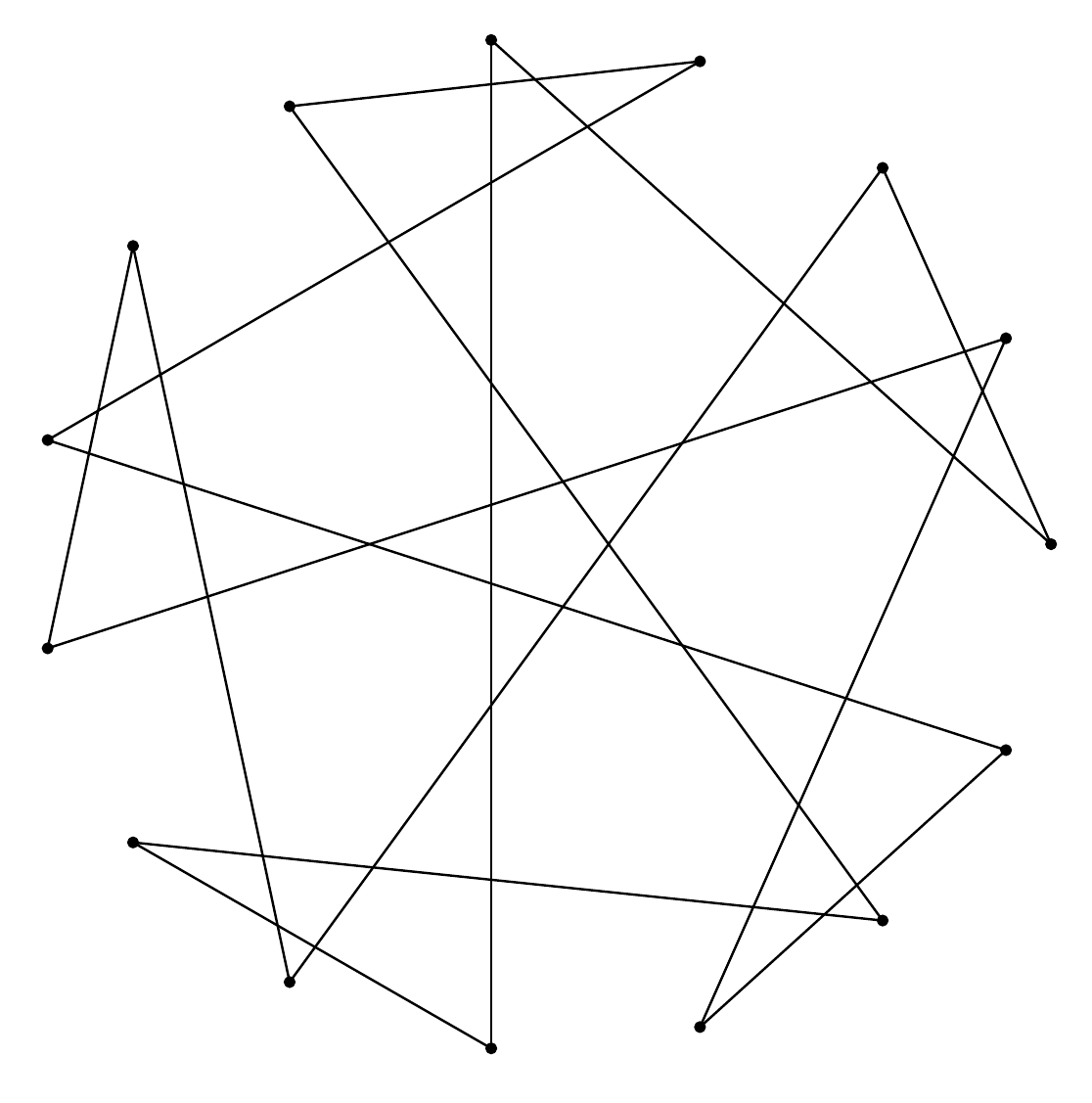} & \includegraphics[width=0.25\textwidth]{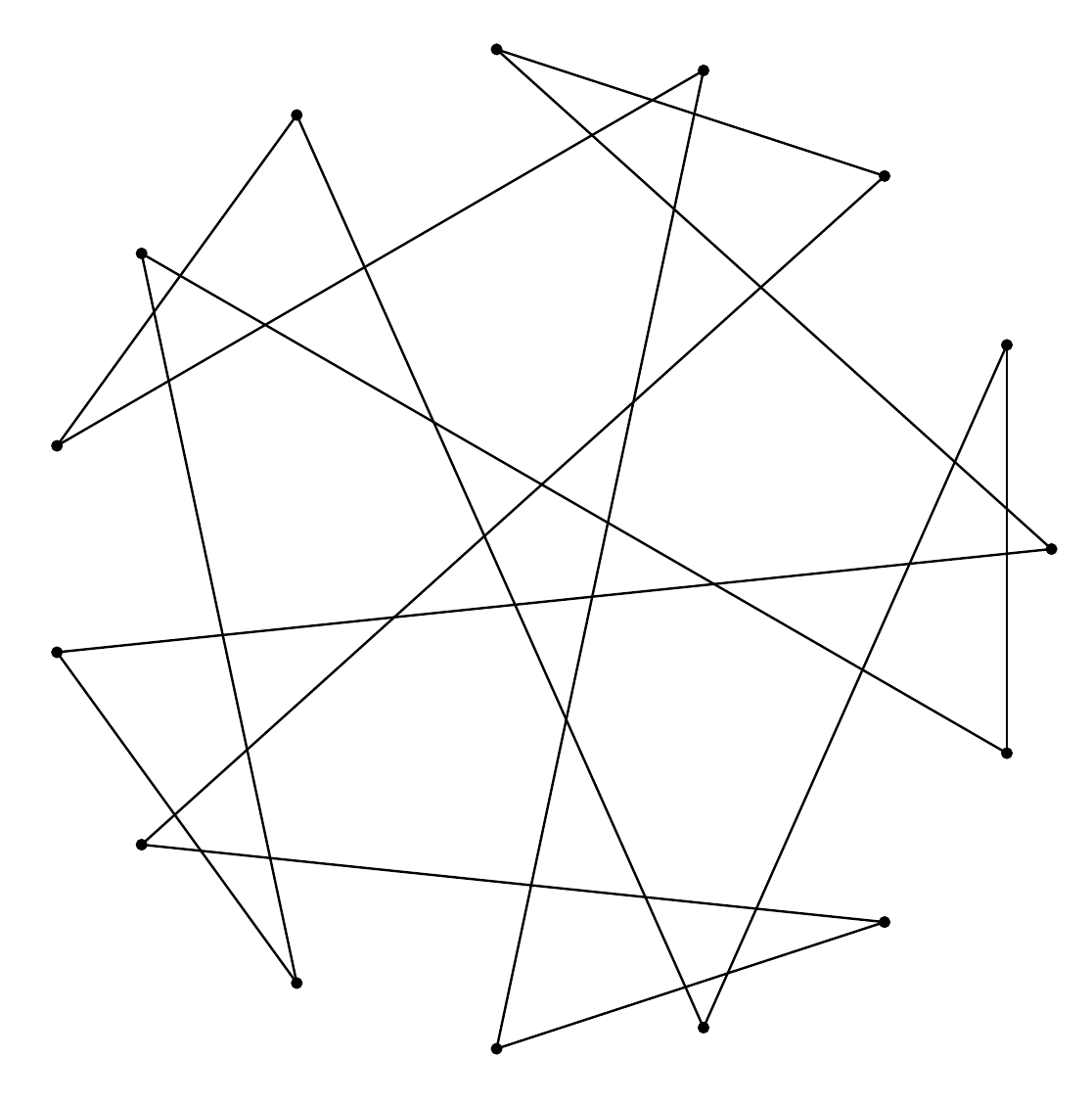} & \includegraphics[width=0.25\textwidth]{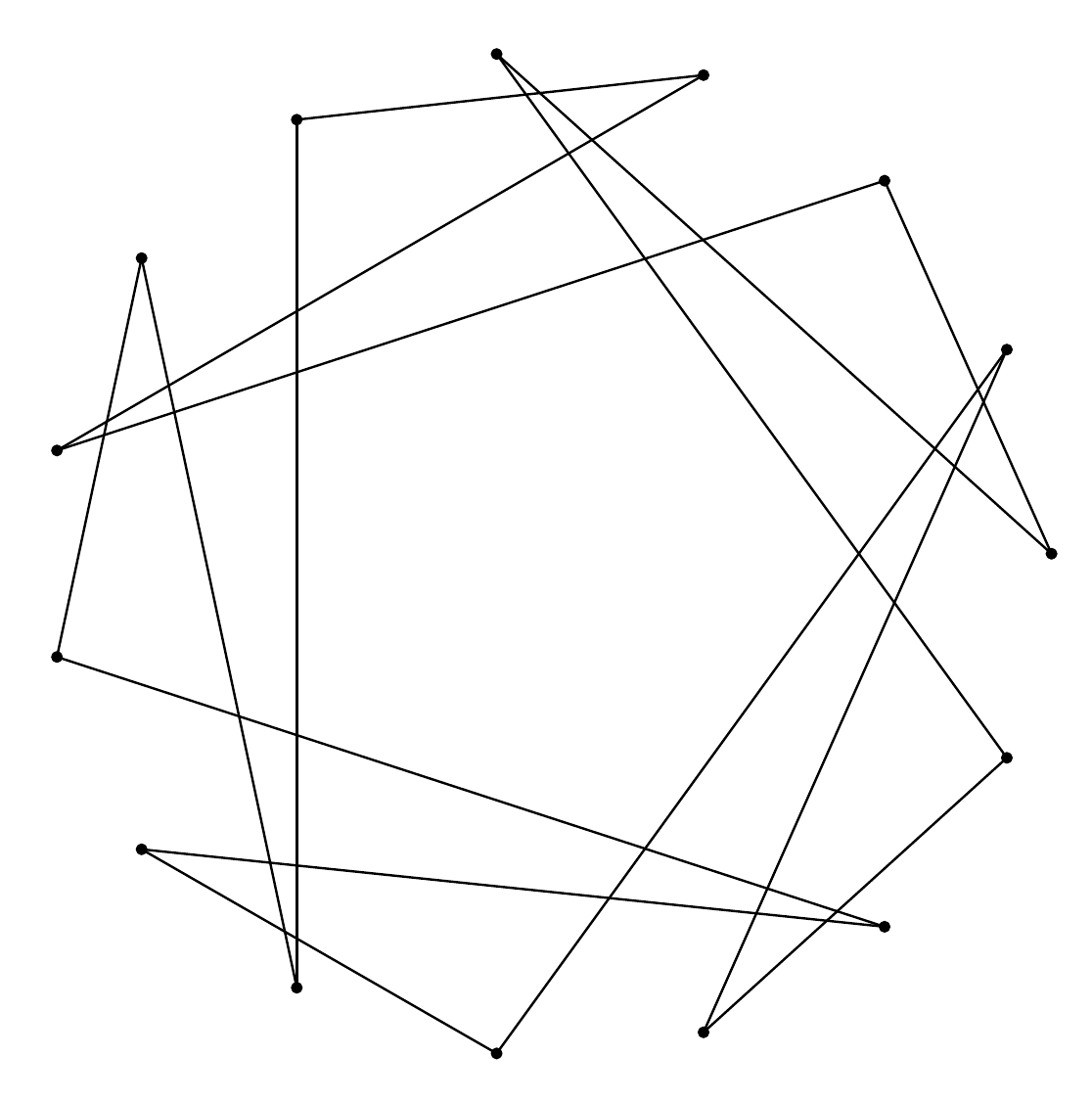} & \includegraphics[width=0.25\textwidth]{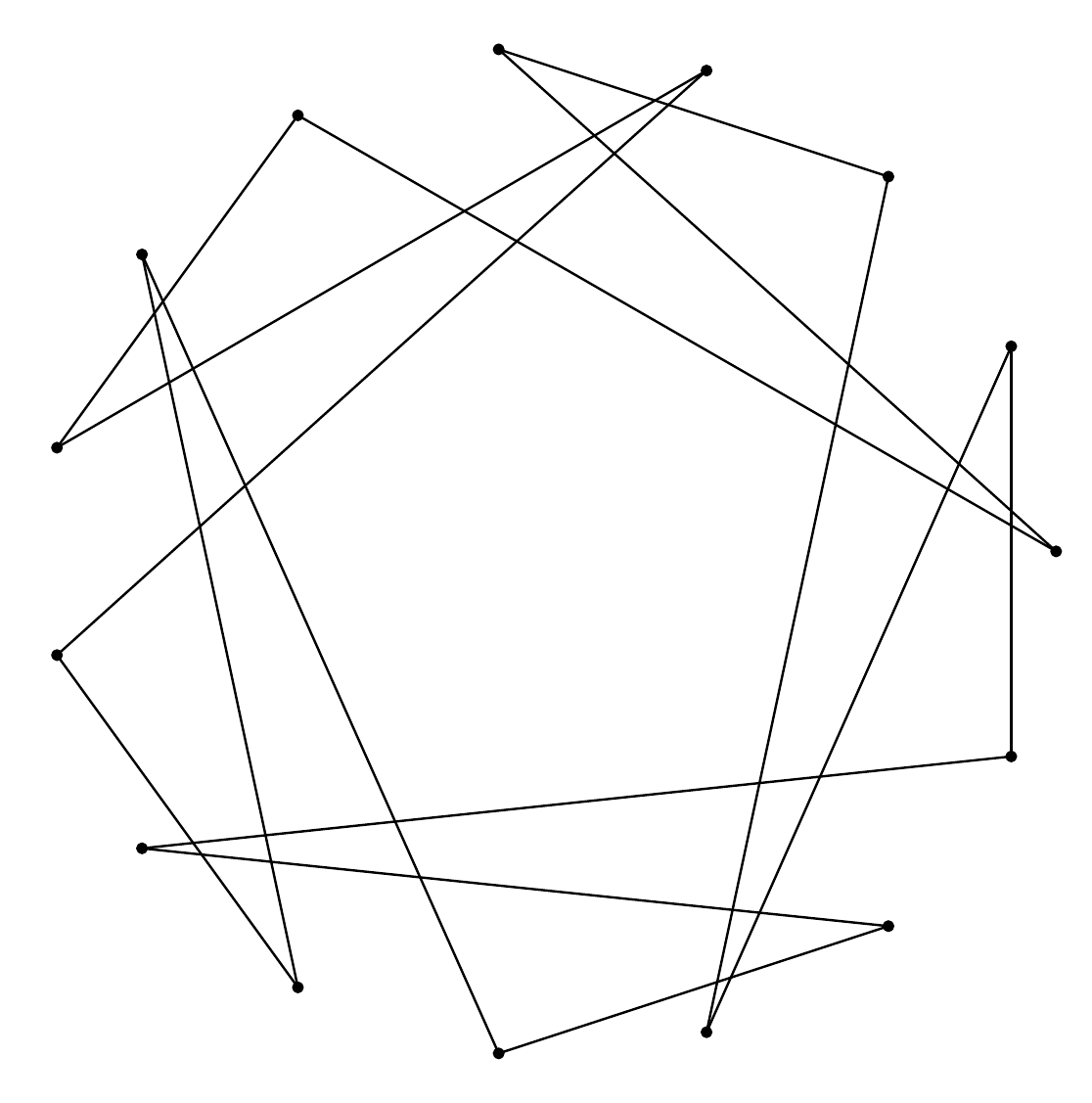}\\
$a=4;b=7;c=13$ & $a=4;b=13;c=7$ & $a=4;b=10;c=13$ & $a=4;b=13;c=10$\\
\end{tabular}
\caption{$n=15$: A set of representatives of the different equivalence-classes of $Q_5(15)$}
\end{figure}

\section{Proofs}
\label{sec:proofs}

\subsection{Preparations}
\label{subsec:preparations}

\subsubsection{Lemma 1}
\label{subsubsec:lemma_1}
Let $m>2$ be an integer and $n=3m$ and $a, b, c$ three integers.
\begin{enumerate}
\item A $m$-circular $3m$-polygon $Q_m(3m)$ owns a sequence of the $n$ sides: $(a, b, c, a, b, c, \ldots, a, b, c)$, with $a \neq b$ and $a \neq c$.
\item A $m$-axial $3m$-polygon $P_m(3m)$ owns a sequence of the $n$ sides: $(a, b, c, a, b, c, \ldots, a, b, c)$, with $a \neq b$ but $a=c$.
\item A completely regular $3m$-polygon $R(3m)$ owns a sequence of the $n$ sides: $(a, b, c, a, b, c, \ldots, a, b, c)$, with $a=b$ and $a=c$.
\end{enumerate}

\textbf{Proof of Lemma 1}\newline
Let $m>2$ be an integer and $n=3m$ and $Q_m(3m)$ a $m$-circular $3m$-polygon.\newline
If $Q_m(3m)$ gets rotated three times by the angle $\dfrac{2\pi}{n}$, then it changes over into itself. A first rotation by the angle $\dfrac{2\pi}{n}$ changes $Q_m(3m)$ with the sequence of sides $(e_1, e_2, \ldots, e_i, \ldots, e_{n-2}, e_{n-1}, e_n)$  into $\overline{Q_m(3m)}$ with the sequence $(e_n, e_1, e_2, \ldots, e_i, \ldots, e_{n-2}, e_{n-1})$ A second rotation of $Q_m(3m)$ by the angle $\dfrac{2\pi}{n}$ changes the sequence $(e_1, e_2, \ldots, e_i, \ldots, e_{n-2}, e_{n-1}, e_n)$ into the sequence \newline $(e_{n-1}, e_n, e_1, e_2, \ldots, e_i, \ldots, e_{n-2})$, a third rotation of $Q_m(3m)$ by the angle $\dfrac{2\pi}{n}$ changes the sequence $(e_1, e_2, \ldots, e_i, \ldots, e_{n-2}, e_{n-1}, e_n)$ into the sequence $(e_{n-2}, e_{n-1}, e_n, e_1, e_2, \ldots, e_i, \ldots, e_{n-3})$. After every third rotation we must get back our  $m$-circular $3m$-polygon $Q_m(3m)$. Therefore we get the following chains of equations, which allows us to define the integers $a,b$ and $c$:
\begin{enumerate}
\item $e_1=e_4=e_7=\cdots:=a$
\item $e_2=e_5=e_8=\cdots:=b$
\item $e_3=e_6=e_9=\cdots:=c$
\end{enumerate}
Thus a sequence of the $3m$ sides $(a, b, c, a, b, c, \ldots, a, b, c)$ describes a $m$-circular $3m$-polygon $Q_m(3m)$.\\

Because a completely regular $3m$-polygon $R(3m)$ and a $3m$-polygon with $m$ axes $P_m(3m)$ can be considered as special cases of a $m$-circular $3m$-polygon $Q_m(3m)$, a sequence of the $3m$ sides $(a, b, c, a, b, c, \ldots, a, b, c)$ describes also $R(3m)$  and $P_m(3m)$. In case of $R(3m)$ we must set $a=b=c$. In case $P_m(3m)$ we get $a=c$ and $b \neq a$.\hfill $\square$\\

\newpage
\subsection{Proof of Main-theorem 1}
\label{subsec:proof_of_main_theorem_1}

We start with an investigation upon the sums of sides and the number $u$ of revolutions of the $3m$-polygon $P_m(3m)$ with $m$ axes:

\subsubsection{Lemma 2 and its proof}
\label{subsubsec:lemma_2_and_its_proof}

Let $m>2$ be an integer and $n=3m$ and let $P_m(3m)$ be a $m$-axial $3m$-polygon.\newline
$P_m(n) = (a, b, a, a, b, a, \ldots, a, b, a)$ has the following sums of sides:

\begin{table}[!htp]
\centering
\begin{tabular}{| c | c |}
\hline
$i$ & $ s _i$\\ \hline
1 & $a$ \\
2 & $a+b$\\
3 & $2a+b$\\
4 & $3a+b$\\
5 & $3a+2b$\\
6 & $4a+2b=2(2a+b)$\\
7 & $5a+2b$\\
8 & $5a+3b$\\
9 & $6a+3b=3(2a+b)$\\
$\cdots$ & $\cdots$\\
$i-2$ & $\dfrac{2i+1}{3}\cdot a+\dfrac{i-1}{3}\cdot b$\\
$i-1$ & $\dfrac{2i-1}{3}\cdot a+\dfrac{i+1}{3}\cdot b$\\
$i$ & $\dfrac{i}{3}\cdot 2a+\dfrac{i}{3}\cdot b=\dfrac{i}{3}\cdot(2a+b)$\\
$\cdots$ & $\cdots$\\
$n-2$ & $\dfrac{2n+1}{3}\cdot a+\dfrac{n-1}{3}\cdot b$\\
$n-1$ & $\dfrac{2n-1}{3}\cdot a+\dfrac{n+1}{3}\cdot b$\\
$n$ & $\dfrac{n}{3}\cdot 2a+\dfrac{n}{3}\cdot b=\dfrac{n}{3}\cdot(2a+b)$\\
\hline
\end{tabular}
\caption{Sums of the sides of a $m$-axial $3m$-polygon}
\label{tab:sums_of_the_sides_of_a_a_m-axial_3m_-polygon}
\end{table}

Thus, we have proved that for a $m$\textbf{-axial} $3m$-polygon, the sums $s_i$ of the sides and \newline the number $u$ of revolutions of the polygon are given as follows:
\begin{enumerate}
\item For $i\equiv $ -1 mod 3: $s_i=\dfrac{2i-1}{3}\cdot a+\dfrac{i+1}{3}\cdot b$
\item For $i\equiv $ 0 mod 3: $s_i=\dfrac{i}{3}\cdot 2a+\dfrac{i}{3}\cdot b=\dfrac{i}{3}\cdot(2a+b)$
\item For $i\equiv $ 1 mod 3: $s_i=\dfrac{2i+1}{3}\cdot a+\dfrac{i-1}{3}\cdot b$
\item $u=\dfrac{2a+b}{3}$
\end{enumerate}\hfill $\square$

\newpage

We now examine the conditions that must be placed on the two integers $a$ and $b$ so that the $n$-tuple $(a,b,a, a,b,a, \ldots, a,b,a)$ represents an $m$-\textbf{axial} $3m$-polygon $P_m(3m)$.\\

We have already shown that $1 \leq a \leq n-1$, $1 \leq b \leq n-1$ and $a \neq b$ are necessary conditions.\\

\subsubsection{Lemma 3 and its proof}
\label{subsubsec:lemma_3_and_its_proof}
If $a\nequiv b$ mod 3, then $(a,b,a, a,b,a, \ldots, a,b,a)$ does not represent a $P_m(3m)$-polygon.\\

We study the 9 possible combinations of $a$ and $b$ with respect to the residual classes modulo 3 and find that $2a + b$ is divisible by $3$, if and only if $a$ and $b$ belong to the same residual class. The divisibility of $2a + b$ by $3$ must be satisfied, because $u=\dfrac{2a+b}{3}$ as the number of revolutions of a $P_m(3m)$-polygon must be an integer.

\begin{table}[!htp]
\centering
\begin{tabular}{| c | c | c |}
\hline
a mod 3 & b mod 3 & (2a+b) mod 3\\ \hline
0 & 0 & 0 \\
0 & 1 & 1 \\
0 & 2 & 2 \\
1 & 0 & 2 \\
1 & 1 & 0 \\
1 & 2 & 1 \\
2 & 0 & 1 \\
2 & 1 & 2 \\
2 & 2 & 0 \\
\hline
\end{tabular}
\caption{Residual classes modulo 3 for $2a+b$}
\label{tab:residual_classes_modulo_3_for_2a+b}
\end{table}\hfill $\square$\\

\subsubsection{Lemma 4 and its proof}
\label{subsubsec:lemma_4_and_its_proof}
If $a\equiv 0$ mod 3 and $b\equiv 0$ mod 3, then $(a,b,a, a,b,a, \ldots, a,b,a)$ does not represent a $P_m(3m)$-polygon.\\

If one starts the construction of a polygon from the $v_1 = 0$ vertex and has only side lengths that are multiples of 3, then the resulting figure only includes vertices which are also multiples of 3, i.e. not all $n$ vertices will be included. Thus, the resulting figure is not a $P_m(3m)$-polygon, but at most one with a smaller number of vertices. \newline
If we start from the $v_1 = 1$ vertex, only vertices whose number belong to the residual class 1 modulo 3 will be involved in the construction, if we start from the $v_1 = 2$ vertex only vertices whose number belong to the residual class 2 modulo 3, will be involved in the construction, and so on. \newline
However we start, we are not able to construct a $P_m(3m)$-polygon, if $a$ and $b$ are both divisible by 3. 

\hfill $\square$\\

So until now we have shown that when $(a,b,a, a,b,a, \ldots, a,b,a)$ represents a $P_m(3m)$-polygon, either $a$ and $b$ both belong to the residual class 1 modulo 3 or both belong to the residual class 2.
\newpage
\subsubsection{Lemma 5 and its proof}
\label{subsubsec:lemma_5_and_its_proof}
Let $(a,b,a, a,b,a, \ldots, a,b,a)$ be a representation of a $P_m(3m)$-polygon with $a\equiv 2$ mod 3 and $b\equiv 2$ mod 3. Then there are two integers $a'$ and $b'$, which are both \newline congruent 1 mod 3, so that $(a',b',a', a',b',a', \ldots, a',b',a')$ represents the same polygon.\\

\textbf{Example for Lemma 5: $n=12$}
\begin{figure}[!htp]
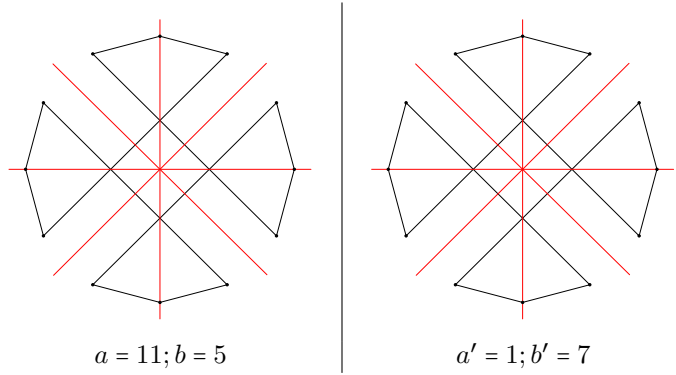

\centering
\begin{tabular}{c | c }
\includegraphics[width=0.3\textwidth]{12_01_07.pdf} & \includegraphics[width=0.3\textwidth]{12_01_07.pdf}\\
$a = 11; b = 5$ & $a' = 1; b' = 7$
\end{tabular}
\caption{$n=12$: The two representations of the same polygon}
\end{figure}

Starting from vertex $v_1 = 0$, each polygon can be traversed in two ways. If the one direction is defined by the $n$-tuple of the sides $(a,b,a, a,b,a, \ldots, a,b,a)$, the $n$-tuple \newline $(n-a,n-b,n-a, n-a,n-b,n-a, \ldots, n-a,n-b,n-a)$ belongs to the opposite direction. Because $n\equiv$ 0 mod 3 and $a\equiv$ 2 mod 3, $(n-a)\equiv$ 1 mod 3. Because $n\equiv$ 0 mod 3 and $b\equiv$ 2 mod 3, $(n-b)\equiv$ 1 mod 3. Thus, $a'=n-a$ and $b'=n-b$ belong to the residual class 1 and the associated $n$-tuple of sides $(a',b',a', a',b',a', \ldots, a',b',a')$ represents the same polygon. \hfill $\square$\\

In the final proof of main-theorem 1, we put $a \equiv$ 1 mod 3 and $b \equiv$ 1 mod 3 without loss of generality and define the following symbol:\newline $(a,b,a, a,b,a, \ldots, a,b,a) \hat{=} P_m(3m):=$ The $n$-tuple of sides $(a,b,a, a,b,a, \ldots, a,b,a)$ represents a $P_m(3m)$-polygon. \\

We now prove the main theorem 1 in the following version:
Let $m>2$, $n=3m$, $a \equiv$ 1 mod 3, $b \equiv$ 1 mod 3, $a \neq b$, $1 \leq a \leq n-2$, $1 \leq b \leq n-2$ and $gcd\left(2a+b,n \right) :=d$. Then follows:\newline
$(a,b,a, a,b,a, \ldots, a,b,a) \hat{=} P_m(3m)\Leftrightarrow d=3$.\\

\subsubsection{Proof of the main-theorem 1$\Rightarrow$:}
\label{subsubsec:proof_of_the_main-theorem_1_rightarrow}
We proof, $d \neq 3\Rightarrow (a,b,a, a,b,a, \ldots, a,b,a) \hat{\neq} P_m(3m)$.\\

$d<3$ is not possible, since $2a+b$ and $n$ are divisible by 3. $\dfrac{2a+b}{3}=u$ is as the number of revolutions a integer. If $d\neq 3$ then $d$ can only accept the values 6, 9, 12 $\ldots$.\newline
So let $d>3$. We show that one of the sums $s_1, s_2, \ldots, s_{n-1}$, is divisible by $n$ and thus not all $n$ vertices are included in the construction.\newline
From $gcd\left(2a+b,n \right)=d>3$ follows the equation $2a+b=g \cdot d$ for an integer $g$ with the property $g<2a+b$.\newline
We choose $i:=\dfrac{3n}{d}<n$. Because $gcd\left(2a+b,n \right)=d$, $\dfrac{n}{d}$ must be an integer and therefore \newline $i=\dfrac{3n}{d}\equiv$ 0 mod 3. Based on Lemma 2 follows: \newline $s_i= \dfrac{3n}{d}\cdot \dfrac{2a+b}{3}=\dfrac{n}{d}\cdot (2a+b)$ and by substituting $2a+b$ further:\newline $s_i=\dfrac{n}{d}\cdot g \cdot d=n \cdot g$. This last equation shows, that $s_i$  is divisible by $n$, while $i$ is less than $n$. \hfill $\Rightarrow\square$\\
\newpage
\subsubsection{Proof of the main-theorem 1$\Leftarrow$:}
\label{subsubsec:proof_of_the_main-theorem_1_leftarrow}

We proof, $(a,b,a, a,b,a, \ldots, a,b,a) \hat{=} P_m(3m)\Leftarrow d=3$\\

That means, that we have to prove: If $d=gcd \left(2a+b,n \right)=3$ no sum of sides except the sum of all $n$ sides is divisible by $n$:\\

\textbf{Part 1:} \newline
The sum of all sides is divisible by $n$, because $s_n=\dfrac{2a+b}{3} \cdot n$ is divisible by $n$, because $\dfrac{2a+b}{3}=u$ is an integer.\\

\textbf{Part 2:} \newline
We want now to prove, that if $d=3$ all sums $s_i$ with $i<n$ and $i\equiv$ 0 mod 3 are not divisible by $n$.\\

We do this prove by showing, that if there exists a index $i$ with $i<n$ and $i\equiv$ 0 mod 3 with the property, that $s_i$ is divisible by $n$, then $d>3$.\\

We choose $i:=\dfrac{3n}{d}$. Because $gcd(2a+b,n)=d$ the number $n$ is divisible by $d$ and therefore the index $i\equiv$ 0 mod 3. We set $i<n$, which includes, that $d>3$. For this index $i$ we show the property, that $s_i$ is divisible by $n$.\\

From Lemma 2 and by substituting $i$ we get: $s_i=\dfrac{i}{3}\cdot(2a+b)=\dfrac{n\cdot(2a+b)}{d}$.\newline Because $gcd(2a+b,n)=d$ we see, that $2a+b$ is divisible by $d$ and $\dfrac{2a+b}{d}:=h$ must be an integer.\newline
We conclude, that $s_i=n \cdot h$, which means, that $s_i$ is divisible by $n$.\\

\textbf{Part 3:} \newline
We want now to prove, that if $d=3$ all sums $s_i$ with $i\leq n-2$ and $i\equiv$ 1 mod 3 are not divisible by $n$.\\

Let $i\equiv$ 1 mod 3 and $i\leq n-2$. As shown in Lemma 2: $s_i=\dfrac{2i+1}{3} \cdot a+\dfrac{i-1}{3} \cdot b$. We can transpose to $s_i=\dfrac{2a+b}{3} \cdot i+\dfrac{a-b}{3}$. We have to show, that in all possible cases $s_i \nequiv$ 0 mod 3.\newline
$s_i=\dfrac{2a+b}{3} \cdot i+\dfrac{a-b}{3}$ \newline
$s_i+\dfrac{2a+b}{3}=\dfrac{2a+b}{3} \cdot i+\dfrac{a-b}{3}+\dfrac{2a+b}{3}$\newline
$s_i=\dfrac{2a+b}{3} \cdot (i-1)+a$\\

Because $i-1 \equiv$ 0 mod 3 and $a \equiv$ 1 mod 3 it doesn't matter to which residual class modulo 3 the factor $\dfrac{2a+b}{3}$ belongs. In every case $s_i \equiv 0+1 \equiv$ 1 mod 3. Therefore we conclude, that if $i\equiv$ 1 mod 3 the sums $s_i$ will never be divisible by 3. On the other hand $n$ is divisible by 3, which together shows, that no sum $s_i$ is divisible by $n$, if $i\equiv$ 1 mod 3.\\

\textbf{Part 4:} \newline
We want now to prove, that if $d=3$ all sums $s_i$ with $i\leq n-1$ and $i\equiv$ 2 mod 3 are not divisible by $n$.\\

Let $i\equiv$ 2 mod 3 and $i\leq n-1$. As shown in Lemma 2: $s_i=\dfrac{2i-1}{3} \cdot a+\dfrac{i+1}{3} \cdot b$. We can transpose to $s_i=\dfrac{2a+b}{3} \cdot i-\dfrac{a-b}{3}$. We have to show, that in all possible cases $s_i \nequiv$ 0 mod 3.\newline
$s_i=\dfrac{2a+b}{3} \cdot i-\dfrac{a-b}{3}$ \newline
$s_i+\dfrac{2a+b}{3}=\dfrac{2a+b}{3} \cdot i -\dfrac{a-b}{3}+\dfrac{2a+b}{3}$\newline
$s_i=\dfrac{2a+b}{3} \cdot (i-1)+\dfrac{a+2b}{3}$\newline
$s_i+\dfrac{2a+b}{3}=\dfrac{2a+b}{3} \cdot (i-1)+\dfrac{a+2b}{3}+\dfrac{2a+b}{3}$\newline
$s_i=\dfrac{2a+b}{3} \cdot (i-2)+a+b$\\

Because $i-2 \equiv$ 0 mod 3 and $a \equiv$ 1 mod 3 and $b \equiv$ 1 mod 3 it doesn't matter to which residual class modulo 3 the factor $\dfrac{2a+b}{3}$ belongs. In every case $s_i \equiv 0+2 \equiv$ 2 mod 3. Therefore we conclude, that if $i\equiv$ 2 mod 3 the sums $s_i$ will never be divisible by 3. On the other hand $n$ is divisible by 3, which together shows, that no sum $s_i$ is divisible by $n$, if $i\equiv$ 2 mod 3.\hfill $\Leftarrow\square$\\

\subsubsection{Conclusion}
\label{conclusion}
Let $m>2$, $n=3m$, $a \equiv$ 1 mod 3, $b \equiv$ 1 mod 3, $a \neq b$, $1 \leq a \leq n-2$, $1 \leq b \leq n-2$, $gcd\left(2a+b,n \right)=3$\\

We obtain a complete set of representatives of the different equivalence classes of the $m$-axial $3m$- polygons $P_m(3m)$ by varying the pairs $(a, b)$ over all allowed values.\\

Obviously: $gcd\left(2a+b,n \right)=3 \Leftrightarrow gcd\left(\dfrac{2a+b}{3},m \right)=1 \Leftrightarrow gcd\left(u,m \right)=1$.\\

Thus, we recognize the validity of Theorem 1 formulated at the beginning:\\

Let $m>2$ be an integer and $n=3m$.\newline
The different equivalence classes of $n$-polygons with $m$ axes are represented by the $n$-tuples $(a, b, a, a, b, a, \ldots, a, b, a)$ of their sides, if $a$ and $b$ have the following six properties:
\begin{enumerate}
\item $a\in \mathbb{N}$ with $a\equiv $1 mod 3,
\item $b\in \mathbb{N}$ with $b\equiv $1 mod 3,
\item $gcd\left(u,m\right)=1$,
\item $1 \leq a \leq n-2$,
\item $1 \leq b \leq n-2$,
\item $a \neq b$.
\end{enumerate}

The question of the number of equivalence classes $\vert P_m(3m) \vert$ will now be answered:

\subsubsection{Proof of the formula for $\vert P_m(3m) \vert$}
\label{proof_of_the_formula_for_|P_m(3m)|}

Let $u$ and $m$ both be integers. Without proof we state that $gcd(u,m)=gcd(u$ mod $m,m)$. We will need this equation below.\\

Let $m>2$ be an integer and $n=3m$. Let $a\in \mathbb{N}$ with $a\equiv $1 mod 3 and $b\in \mathbb{N}$ with $b\equiv $1 mod 3.\newline Let $1 \leq a \leq n-2$ and $1 \leq b \leq n-2$.\\

We set up tables by varying over all possible combinations of $a$ and $b$. Because $a$ and $b$ vary both between 1 and $3m-2$ and are both equivalent to 1 mod 3, we need $m$ tables with $m$ rows. So for $i$ from 1 to $m$ we consider the $i$th table:
\newpage

\begin{table}[!htp]
\centering
\begin{tabular}{| c | c | c | c | c | c |}
\hline
$j$ & $a$ & $b$ & $a$ & $2a+b$ & $u$\\ \hline
$1$ & $3i-2$ & $1$ & $3i-2$ & $6i-3$ & $2i-1$\\
$2$ & $3i-2$ & $4$ & $3i-2$ & $6i$ & $2i$\\
$3$ & $3i-2$ & $7$ & $3i-2$ & $6i+3$ & $2i+1$\\
$\cdots$ & $\cdots$ & $\cdots$ & $\cdots$ & $\cdots$ & $\cdots$\\
$j$ & $3i-2$ & $3j-2$ & $3i-2$ & $6i+3j-6$ & $2i+j-2$\\
$\cdots$ & $\cdots$ & $\cdots$ & $\cdots$ & $\cdots$ & $\cdots$\\
$m$ & $3i-2$ & $3m-2$ & $3i-2$ & $6i+3m-2$ & $2i+m-2$\\ \hline
\end{tabular}
\caption{$i$-th table of the combinations of $a$ and $b$}
\label{tab:i-th_table_of_the_combinations_of_a_and_b}
\end{table}

We recognize in the $u$-column of each table an arithmetic series with the difference 1, which begins at $2i-1$ and ends with $2i+m-2$, i.e:\\

In the first table $u$ takes the values 1, 2, 3, $\cdots, m$. There are $\varphi(m)$ rows with $gcd(u,m)=1$.\newline
In the second table $u$ takes the values 3, 4, 5, $\cdots, m, m+1, m+2$. There are also $\varphi(m)$ rows with $gcd(u,m)=1$, because, as mentioned above, $gcd(u,m)=gcd(u$ mod $m,m)$, and so on.\\

In each of the $m$ tables there are exactly $\varphi(m)$ rows, in which the $gcd(u,m)=1$. Therefore in all $m$ tables there are exactly $m \cdot \varphi(m)$ rows, in which the $gcd(u,m)=1$.\\

In this counting of the rows, in which $u$ is prime to $m$, those combinations of $a$ and $b$ with $a=b$ are also taken into account. To get $\vert P_m(3m)\vert $ we must exclude the completely regular $3m$-polygons, which appear in the tables. Because all completely regular $3m$-polygons appear exactly once in the tables, we have to exclude $\dfrac{\varphi(3m)}{2}$ rows, and we get the stated result:
\begin{center}
$\vert P_m(3m) \vert = m \cdot \varphi(m)-\dfrac{\varphi(3m)}{2}$.
\end{center}\hfill $\square$\\

\subsubsection{Conclusion $\vert P_p(3p) \vert$}
\label{conclusion_p_p(3p)}
Let $p$ be prime and $p>3$. \newline $\vert P_p(3p) \vert=p \cdot \varphi(p)-\dfrac{\varphi(3p)}{2}=p \cdot \varphi(p)-\dfrac{\varphi(3)\cdot\varphi(p)}{2}=p \cdot \varphi(p)-\varphi(p)=\varphi(p)\cdot(p-1)=\underline{\underline{(p-1)^2}}$ \hfill $\square$\\

\newpage
\subsection{Proof of Main-theorem 2}
\label{subsec:proof_of_main_theorem_2}

We start with an investigation upon the sums of sides and the number $u$ of revolutions of the $m$-circular polygon $Q_m(3m)$:

\subsubsection{Lemma 6 and its proof}
\label{subsubsec:lemma_6_and_its_proof}

\textbf{The sums of the sides and the revolutions of the $m$-circular $3m$-polygon $Q_m(3m)$:}\\

Let $m>2$ be an integer and $n=3m$ and let $Q_m(3m)$ be a $m$-circular $3m$-polygon.\newline
$Q_m(n) = (a, b, c, a, b, c, \ldots, a, b, c)$ has the following sums of sides:

\begin{table}[!htp]
\centering
\begin{tabular}{| c | c |}
\hline
$i$ & $ s _i$\\ \hline
1 & $a$ \\
2 & $a+b$\\
3 & $a+b+c$\\
4 & $2a+b+c$\\
5 & $2a+2b+c$\\
6 & $2a+2b+2c=2(a+b+c)$\\
7 & $3a+2b+2c$\\
8 & $3a+3b+2c$\\
9 & $3a+3b+3c=3(a+b+c)$\\
$\cdots$ & $\cdots$\\
$i-2$ & $\dfrac{i+2}{3}\cdot a+\dfrac{i-1}{3}\cdot b+\dfrac{i-1}{3}\cdot c$\\
$i-1$ & $\dfrac{i+1}{3}\cdot a+\dfrac{i+1}{3}\cdot b+\dfrac{i-2}{3}\cdot c$\\
$i$ & $\dfrac{i}{3}\cdot a+\dfrac{i}{3}\cdot b+\dfrac{i}{3}\cdot c=\dfrac{i}{3}\cdot(a+b+c)$\\
$\cdots$ & $\cdots$\\
$n-2$ & $\dfrac{n+2}{3}\cdot a+\dfrac{n-1}{3}\cdot b+\dfrac{n-1}{3}\cdot c$\\
$n-1$ & $\dfrac{n+1}{3}\cdot a+\dfrac{n+1}{3}\cdot b+\dfrac{n-2}{3}\cdot c$\\
$n$ & $\dfrac{n}{3}\cdot a+\dfrac{n}{3}\cdot b+\dfrac{n}{3}\cdot c=\dfrac{n}{3}\cdot(a+b+c)$\\
\hline
\end{tabular}
\caption{Sums of the sides of a $m$-circular $3m$-polygon}
\label{tab:sums_of_the_sides_of_a_a_m-circular_3m_-polygon}
\end{table}
Thus, we have proved that for a $m$\textbf{-circular} $3m$-polygon, the sums $s_i$ of the sides and \newline the number $u$ of revolutions of the polygon are given as follows:
\begin{enumerate}
\item For $i\equiv $ -1 mod 3: $s_i=\dfrac{i+1}{3}\cdot a+\dfrac{i+1}{3}\cdot b+\dfrac{i-2}{3}\cdot c$
\item For $i\equiv $ 0 mod 3: $s_i=\dfrac{i}{3}\cdot a+\dfrac{i}{3}\cdot b+\dfrac{i}{3}\cdot c=\dfrac{i}{3}\cdot(a+b+c)$
\item For $i\equiv $ 1 mod 3: $s_i=\dfrac{i+2}{3}\cdot a+\dfrac{i-1}{3}\cdot b+\dfrac{i-1}{3}\cdot c$
\item $u=\dfrac{a+b+c}{3}$.
\end{enumerate}\hfill $\square$

\subsubsection{Lemma 7 and its proof}
\label{subsubsec:lemma_7_and_its_proof}
Let $m>2$ be an integer and $n=3m$ and let $Q_m(3m)$ be a $m$-circular $3m$-polygon. Then $a\nequiv $ 0 mod 3 and $b\nequiv $ 0 mod 3 and $c\nequiv $ 0 mod 3.\\

Suppose, that $a\equiv $ 0 mod 3, i.e. There exists a integer $k\leq m$, such that $a=3 \cdot k$. We can set a first edge by connecting the vertex $v_1=0$ with the vertex $v_2=3k$, i.e. $e_1=\overline{0\smallsmile 3k}$. Because the construction is to give a $m$-circular $3m$-polygon, a second edge must connect the vertices $3k$ and $6k$, i.e. $e_2=\overline{3k\smallsmile 6k}$. \newline We get $m$ edges $\overline{0\smallsmile 3k},\overline{3k\smallsmile 6k},\ldots,\overline{(m-1)3k\smallsmile m3k}=\overline{(m-1)3k\smallsmile nk}=\overline{(m-1)3k\smallsmile0}$, one after the other and see that the polygon closes prematurely by going through at most $m$ vertices. Therefore our assumption, that $a\equiv $ 0 mod 3 is false. The same argument is valid for $b$ and $c$. \hfill $\square$

\subsubsection{Lemma 8 and its proof}
\label{subsubsec:lemma_8_and_its_proof}
Now let $a\nequiv$ 0 mod 3 and $b\nequiv$ 0 mod 3 and $c\nequiv$ 0 mod 3. \newline If $a\nequiv b$ mod 3 or $b\nequiv c$ mod 3 or $c\nequiv a$ mod 3, then $(a,b,c, a,b,c, \ldots, a,b,c)$ does not represent a $Q_m(3m)$-polygon.\\

We study the 8 possible combinations of $a, b$ and $c$ with respect to the residual classes modulo 3 and find that $a + b+c$ is divisible by $3$, if and only if $a, b$ and $c$ belong to the same residual class. The divisibility of $a + b+c$ by $3$ must be satisfied, because $u=\dfrac{a+b+c}{3}$ as the number of revolutions of a $Q_m(3m)$-polygon must be an integer.

\begin{table}[!htp]
\centering
\begin{tabular}{| c | c | c | c |}
\hline
a mod 3 & b mod 3 & c mod 3 &(a+b+c) mod 3\\
\hline
1 & 1 & 1 & 0  \\
1 & 1 & 2 & 1  \\
1 & 2 & 1 & 1  \\
2 & 1 & 1 & 1  \\
2 & 2 & 2 & 0  \\
2 & 2 & 1 & 2  \\
2 & 1 & 2 & 2  \\
1 & 2 & 2 & 2  \\
\hline
\end{tabular}
\caption{Residual classes modulo 3 for $a+b+c$}
\label{tab:residual_classes_modulo_3_for_a+b+c}
\end{table}\hfill $\square$ \\

So until now we have shown that when $(a,b,c, a,b,c, \ldots, a,b,c)$ represents a $Q_m(3m)$-polygon, either $a,b$ and $c$ all belong to the residual class 1 modulo 3 or all belong to the residual class 2.

\subsubsection{Lemma 9 and its proof}
\label{subsubsec:lemma_9_and_its_proof}
Let $(a,b,c, a,b,c, \ldots, a,b,c)$ be a representation of a $Q_m(3m)$-polygon with $a\equiv 2$ mod 3, $b\equiv 2$ mod 3 and $c\equiv 2$ mod 3. Then there exist three integers $a'$, $b'$ and $c'$ which are all \newline congruent 1 mod 3, so that $(a',b',c', a',b',c', \ldots, a',b',c')$ represents the same polygon.\\

Starting from vertex $v_1 = 0$, each polygon can be traversed in two ways. If the one direction is defined by the $n$-tuple of the sides $(a,b,c, a,b,c, \ldots, a,b,c)$, the $n$-tuple \newline $(n-a,n-b,n-c, n-a,n-b,n-c, \ldots, n-a,n-b,n-c)$ belongs to the opposite direction. Because $n\equiv$ 0 mod 3 and $a\equiv$ 2 mod 3, $(n-a)\equiv$ 1 mod 3. Because $n\equiv$ 0 mod 3 and $b\equiv$ 2 mod 3, $(n-b)\equiv$ 1 mod 3. And because $n\equiv$ 0 mod 3 and $c\equiv$ 2 mod 3, $(n-c)\equiv$ 1 mod 3.Thus, $a':=n-a$, $b':=n-b$ and $c':=n-c$ belong to the residual class 1 and the associated $n$-tuple of sides $(a',b',c', a',b',c', \ldots, a',b',c')$ represents the same polygon. \hfill $\square$\\

In the final proof of main-theorem 2, we put $a \equiv$ 1 mod 3, $b \equiv$ 1 mod 3 and $c \equiv$ 1 mod 3 without loss of generality and define the following symbol:\newline $(a,b,c, a,b,c, \ldots, a,b,c) \hat{=} Q_m(3m):=$ The $n$-tuple of sides $(a,b,c, a,b,c, \ldots, a,b,c)$ represents a $Q_m(3m)$-polygon. \\

We now prove the main theorem 2 in the following version:
Let $m>2$, $n=3m$, $a \equiv$ 1 mod 3, $b \equiv$ 1 mod 3, $c \equiv$ 1 mod 3, $a \neq b$, $b \neq c$, $c \neq a$, $1 \leq a \leq n-2$, $1 \leq b \leq n-2$, $1 \leq c \leq n-2$ and $gcd\left(a+b+c,n \right) :=d$. Then follows:\newline
$(a,b,c, a,b,c, \ldots, a,b,c) \hat{=} Q_m(3m)\Leftrightarrow d=3$.\\

\subsubsection{Proof of the main-theorem 2$\Rightarrow$:}
\label{subsubsec:proof_of_the_main-theorem_2_rightarrow}
We proof, $d \neq 3\Rightarrow (a,b,c, a,b,c, \ldots, a,b,c) \hat{\neq} Q_m(3m)$.\\

$d<3$ is not possible, since $a+b+c$ and $n$ are divisible by 3. $\dfrac{a+b+c}{3}=u$ is as the number of revolutions a integer. If $d\neq 3$ then $d$ can only accept the values 6, 9, 12 $\ldots$.\\

So let be $d>3$. We show that one of the sums $s_1, s_2, \ldots, s_{n-1}$, is divisible by $n$ and thus not all $n$ vertices are included in the construction.\\

From $gcd\left(a+b+c,n \right) =d>3$ follows the equation $a+b+c=g \cdot d$ for an integer $g$ with the property $g< a+b+c$.\\

We choose $i:=\dfrac{3n}{d}<n$. Because $gcd\left(a+b+c,n \right) =d$, the number $\dfrac{n}{d}$ must be an integer and therefore $i=\dfrac{3n}{d}\equiv$ 0 mod 3.\\

 Based on Lemma 6 follows: $s_i= \dfrac{3n}{d}\cdot \dfrac{a+b+c}{3}=\dfrac{n}{d}\cdot (a+b+c)$ and by substituting $a+b+c$ further:\newline $s_i=\dfrac{n}{d}\cdot g \cdot d=n \cdot g$. This last equation shows, that $s_i$  is divisible by $n$, while $i$ is less than $n$. \hfill $\Rightarrow\square$\\

\subsubsection{Proof of the main-theorem 2$\Leftarrow$:}
\label{subsubsec:proof_of_the_main-theorem_1_leftarrow}

We proof, $(a,b,c, a,b,c, \ldots, a,b,c) \hat{=} Q_m(3m)\Leftarrow d=3$\\

That means, that we have to prove: If $d=gcd \left(a+b+c,n \right)=3$ no sum of sides except the sum of all $n$ sides is divisible by $n$:\\

\textbf{Part 1:} \newline
The sum of all sides is divisible by $n$, because $s_n=\dfrac{a+b+c}{3} \cdot n$ is divisible by $n$, because $\dfrac{a+b+c}{3}=u$ is an integer.\\

\textbf{Part 2:} \newline
We want now to prove, that if $d=3$ all sums $s_i$ with $i<n$ and $i\equiv$ 0 mod 3 are not divisible by $n$.\\

We do this prove by showing, that if there exists a index $i$ with $i<n$ and $i\equiv$ 0 mod 3 with the property, that $s_i$ is divisible by $n$, then $d>3$.\\

We choose $i:=\dfrac{3n}{d}$. Because $gcd(a+b+c,n)=d$ the number $n$ is divisible by $d$ and therefore the index $i\equiv$ 0 mod 3. We set $i<n$, which includes, that $d>3$. For this index $i$ we show the property, that $s_i$ is divisible by $n$.\\

From Lemma 6 and by substituting $i$ we get: $s_i=\dfrac{i}{3}\cdot(a+b+c)=\dfrac{n\cdot(a+b+c)}{d}$.\newline Because $gcd(a+b+c,n)=d$ we see, that $a+b+c$ is divisible by $d$ and $\dfrac{a+b+c}{d}:=h$ must be an integer.\newline
We conclude, that $s_i=n \cdot h$, which means, that $s_i$ is divisible by $n$.\\

\textbf{Part 3:} \newline
We want now to prove, that if $d=3$ all sums $s_i$ with $i\leq n-2$ and $i\equiv$ 1 mod 3 are not divisible by $n$.\\

Let $i\equiv$ 1 mod 3 and $i\leq n-2$. As shown in Lemma 6: $s_i=\dfrac{i+2}{3} \cdot a+\dfrac{i-1}{3} \cdot b+\dfrac{i-1}{3}\cdot c$. We can transpose to $s_i=\dfrac{a+b+c}{3} \cdot i+\dfrac{2a-b-c}{3}$. We have to show, that in all possible cases \newline $s_i \nequiv$ 0 mod 3 becomes true.\newline
$s_i=\dfrac{a+b+c}{3} \cdot i+\dfrac{2a-b-c}{3}$ \newline
$s_i+\dfrac{a+b+c}{3}=\dfrac{a+b+c}{3} \cdot i+\dfrac{2a-b-c}{3}+\dfrac{a+b+c}{3}$\newline
$s_i=\dfrac{a+b+c}{3} \cdot (i-1)+a$\\

Because $i-1 \equiv$ 0 mod 3 and $a \equiv$ 1 mod 3 it doesn't matter to which residual class modulo 3 the factor $\dfrac{a+b+c}{3}$ belongs. In every case $s_i \equiv 0+1 \equiv$ 1 mod 3. Therefore we conclude, that if $i\equiv$ 1 mod 3 the sums $s_i$ will never be divisible by 3. On the other hand $n$ is divisible by 3, which together shows, that no sum $s_i$ is divisible by $n$, if $i\equiv$ 1 mod 3.\\

\textbf{Part 4:} \newline
We want now to prove, that if $d=3$ all sums $s_i$ with $i\leq n-1$ and $i\equiv$ 2 mod 3 are not divisible by $n$.\\

Let $i\equiv$ 2 mod 3 and $i\leq n-1$. As shown in Lemma 6: $s_i=\dfrac{i+1}{3} \cdot a+\dfrac{i+1}{3} \cdot b+\dfrac{i-2}{3}\cdot c$. We can transpose to $s_i=\dfrac{a+b+c}{3} \cdot i+\dfrac{a+b-2c}{3}$. We have to show, that in all possible cases $s_i \nequiv$ 0 mod 3.\newline
$s_i=\dfrac{a+b+c}{3} \cdot i+\dfrac{a+b-2c}{3}$ \newline
$s_i+\dfrac{a+b+c}{3}=\dfrac{a+b+c}{3} \cdot i +\dfrac{2a+2b-c}{3}$\newline
$s_i=\dfrac{a+b+c}{3} \cdot (i-1)+\dfrac{2a+2b-c}{3}$\newline
$s_i+\dfrac{a+b+c}{3}=\dfrac{a+b+c}{3} \cdot (i-1)+\dfrac{3a+3b}{3}$\newline
$s_i=\dfrac{2a+b}{3} \cdot (i-2)+a+b$\\

Because $i-2 \equiv$ 0 mod 3 and $a \equiv$ 1 mod 3 and $b \equiv$ 1 mod 3 it doesn't matter to which residual class modulo 3 the factor $\dfrac{a+b+c}{3}$ belongs. In every case $s_i \equiv 0+2 \equiv$ 2 mod 3. Therefore we conclude, that if $i\equiv$ 2 mod 3 the sums $s_i$ will never be divisible by 3. On the other hand $n$ is divisible by 3, which together shows, that no sum $s_i$ is divisible by $n$, if $i\equiv$ 2 mod 3.\hfill $\Leftarrow\square$\\

\subsubsection{Conclusion}
\label{conclusion_2}
Let $m>2$, $n=3m$, $a \equiv$ 1 mod 3, $b \equiv$ 1 mod 3, $c \equiv$ 1 mod 3, $a \neq b$, $b \neq c$, $c \neq a$, $1 \leq a \leq n-2$, $1 \leq b \leq n-2$, $1 \leq c \leq n-2$, $gcd\left(2a+b,n \right)=3$\\

We obtain a complete set of representatives of the equivalence classes of the $m$-circular $3m$- polygons $Q_m(3m)$ by varying the three numbers $a,b$ and $c$ over all allowed values.\\

Obviously: $gcd\left(a+b+c,n \right)=3 \Leftrightarrow gcd\left(\dfrac{a+b+c}{3},m \right)=1 \Leftrightarrow gcd\left(u,m \right)=1$.\\

Thus, we recognize the validity of Theorem 2 formulated at the beginning:\\

Let $m>2$ be an integer and $n=3m$.\newline
The different equivalence classes of $m$-circular $3m$-polygons are represented by the $n$-tuples $(a, b, c, a, b, c, \ldots, a, b, c)$ of their sides, if $a$, $b$ and $c$ have the following ten properties:
\begin{enumerate}
\item $a\in \mathbb{N}$ with $a\equiv $1 mod 3,
\item $b\in \mathbb{N}$ with $b\equiv $1 mod 3,
\item $c\in \mathbb{N}$ with $c\equiv $1 mod 3,
\item $gcd\left(u,m \right)=1$,
\item $1 \leq a\leq n-2$,
\item $1 \leq b\leq n-2$,
\item $1 \leq c\leq n-2$,
\item $a \neq b$,
\item $b \neq c$,
\item $c \neq a$.
\end{enumerate}

The question of the number of equivalence classes $\vert Q_m(3m) \vert$ will now be answered:

\subsubsection{Proof of the formula for $\vert Q_m(3m) \vert$}
\label{proof_of_the_formula_for_|Q_m(3m)|}

Let $u$ and $m$ both be integers. Without proof we state that $gcd(u,m)=gcd(u$ mod $m,m)$. We will need this equation below.\\

Let $m>2$ be an integer and $n=3m$. Let $a\in \mathbb{N}$ with $a\equiv $1 mod 3, $b\in \mathbb{N}$ with $b\equiv $1 mod 3 and $c\in \mathbb{N}$ with $c\equiv $1 mod 3. And let $1 \leq a \leq n-2$, $1 \leq b \leq n-2$ and $1 \leq c \leq n-2$.\\

We set up tables by varying over all possible combinations of $a$, $b$ and $c$. Because $a$, $b$ and $c$ vary all between 1 and $3m-2$ and are all equivalent to 1 mod 3, we need $m^2$ tables with $m$ rows in each table. So for $k$ from 1 to $m$ and $j$ from 1 to $m$:

\begin{table}[!htp]
\centering
\begin{tabular}{| c | c | c | c | c | c |}
\hline
$i$ & $a$ & $b$ & $c$ & $a+b+c$ & $u$\\ \hline
$1$ & $3k-2$ & $3j-2$ & $1$ & $3k+3j-3$ & $k+j-1$\\
$2$ & $3k-2$ & $3j-2$ & $4$ & $3k+3j$ & $k+j$\\
$3$ & $3k-2$ & $3j-2$ & $7$ & $3k+3j+3$ & $k+j+1$\\
$\cdots$ & $\cdots$ & $\cdots$ & $\cdots$ & $\cdots$ & $\cdots$\\
$i$ & $3k-2$ & $3j-2$ & $3i-2$ & $3k+3j+3i-6$ & $k+j+i-2$\\
$\cdots$ & $\cdots$ & $\cdots$ & $\cdots$ & $\cdots$ & $\cdots$\\
$m$ & $3k-2$ & $3j-2$ & $3m-2$ & $3k+3j+3m-6$ & $k+j+m-2$\\ \hline
\end{tabular}
\caption{Table of the combinations of $a$, $b$ and $c$}
\label{tab:table_of_the_combinations_of_a_b_c}
\end{table}

We recognize in the $u$-column of each table an arithmetic series with the difference 1, which begins at $k+j-1$ and ends with $k+j+m-2$. Because, as mentioned above, $gcd(u,m)=gcd(u$ mod $m,m)$ there exist in each table exactly $\varphi(m)$ rows with $gcd(u,m)=1$. In all tables together we obtain $m^2 \cdot\varphi(m)$ rows with $gcd(u,m)=1$.\\

Now consider all the tables together, but only those $m^2\cdot\varphi(m)$ rows, where $u$ is prime to $m$. In these rows occur:
\begin{enumerate}
\item Each equivalence-class of the $m$-circular  $3m$-polygons $Q_m(3m)$ three times. The rows with $a, b, c$ and  $b,c,a$ and $c,a,b$ belong to $m$-circular $3m$-polygons of the same equivalence-class.
\item Each equivalence-class of the $m$-axial $3m$-polygons $P_m(3m)$ three times. The rows with $a,a,b$ and $b,a,a$ and $a,b,a$ belong to $m$-axial $3m$-polygons of the same equivalence-class.
\item Each of the fully-regular $3m$-polygons once. The rows with $a,a,a$ and $gcd(a,m)=1$.
\end{enumerate}

This gives the equation:
\begin{center}
$m^2\cdot \varphi(m)=3 \cdot |Q_m(3m)|+3\cdot |P_m(3m)|+\dfrac{\varphi(3m)}{2}$.
\end{center}

By solving for $|Q_m(3m)|$ and using the already proven term for $|P_m(3m)|$, we obtain the number $|Q_m(3m)|$ of equivalence classes of the $m$-circular $3m$-polygons $Q_m(3m)$ for $m>2$:

\begin{center}
$|Q_m(3m)|=\dfrac{m \cdot \varphi(m) \cdot (m-3)+\varphi(3m)}{3}$.
\end{center}\hfill $\square$

\subsubsection{Conclusion $\vert Q_p(3p) \vert$}
\label{conclusion_q_p_3p}

Let $p$ be prime and $p>3$.\\

$|Q_p(3p)| = \dfrac{p \cdot (p-1) \cdot (p-3)+2\cdot(p-1)}{3}=\underline{\underline{\dfrac{(p-2) \cdot (p-1)\cdot (p+1)}{3}}}$

\listoffigures
\listoftables
\bibliographystyle{alpha}
\bibliography{Literatur}
\end{document}